\documentclass{article}
\usepackage[dvipsnames,svgnames,x11names]{xcolor}
\usepackage{arxiv}
\usepackage{amsmath}
\numberwithin{equation}{section}
\usepackage[utf8]{inputenc} 
\usepackage{hyperref}       
\usepackage{url}            
\usepackage{booktabs}       
\usepackage{amsfonts}       
\usepackage{nicefrac}       
\usepackage{microtype}      
\usepackage{cleveref}       
\usepackage{lipsum}         
\usepackage{graphicx}
\usepackage{natbib}
\usepackage{doi}
\usepackage{tikz}
\usepackage{amsthm}
\usepackage{esint}
\tikzset{>=stealth}
\usepackage{bm}
\usepackage{yhmath}
\usetikzlibrary{patterns}
\usepackage{amssymb}
\usepackage{threeparttable}
\usepackage{makecell}

\title{Complex variable solution on over-/under-break shallow tunnelling in gravitational geomaterial with reasonable far-field displacement}

\date{\today}
\author{
  \href{https://orcid.org/0000-0002-5143-2714}
  {\includegraphics[scale=0.06]{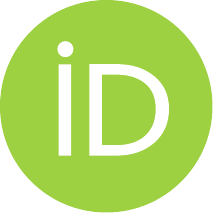}
    \bf {Luo-bin Lin}} \\
  Fujian Provincial Key Laboratory of Advanced Technology and Informatization in Civil Engineering\\
  College of Civil Engineering\\
  Fujian University of Technology\\
  No. 69 Xueyuan Road, Shangjie University Town, Fuzhou, 350118, Fujian, China \\
  \texttt{luobin\_lin@fjut.edu.cn} \\
  \href{https://orcid.org/0000-0002-5583-3734}
  {\includegraphics[scale=0.06]{orcid.pdf}
    \bf {Fu-quan Chen}} \\
  College of Civil Engineering\\
  Fuzhou University\\
  No. 2 Xueyuan Road, Shangjie University Town, Fuzhou, 350108, Fujian, China \\
  \texttt{phdchen@fzu.edu.cn}\\
  \href{https://orcid.org/0000-0002-1302-5296}
  {\includegraphics[scale=0.06]{orcid.pdf}
    \bf {Jin-ping Zhuang}}\\
  Fujian Provincial Key Laboratory of Advanced Technology and Informatization in Civil Engineering\\
  College of Civil Engineering\\
  Fujian University of Technology\\
  No. 69 Xueyuan Road, Shangjie University Town, Fuzhou, 350118, Fujian, China \\
  \texttt{19781208@fjut.edu.cn} \\
}

\renewcommand{\shorttitle}{}

\hypersetup{
  pdftitle={},
  pdfauthor={Luo-bin Lin, Fu-quan Chen, Jin-ping Zhuang},
  pdfkeywords={},
}

\begin{document}
\maketitle

\begin{abstract}
  Over-/under-break excavation is a common phenomenon in shallow tunnelling, which is nonetheless not generally considered in existing complex variable solutions. In this paper, a new equilibrium mechanical model on over-/under-break shallow tunnelling in gravitational geomaterial is established by fixing far-field ground surface to form a corresponding mixed boundary problem. With integration of a newly proposed bidirectional composite conformal mapping using Charge Simulation Method, a complex variable solution of infinite complex potential series is subsequently derived using analytic continuation to tranform the mixed boundaries into a homogenerous Riemann-Hilbert problem, which is iteratively solved to obtain the stress and displacement in geomaterial. The infinite complex potential series of the complex variable solution are truncated to obtain numerical results, which is rectified by Lanczos filtering to reduce the oscillation of Gibbs phenomena. The bidirectional conformal mapping is discussed and validated via several numerical cases, and the subsequent complex variable solution is verified by examining the Lanczos filtering and solution convergence, and comparing with corresponding finite element solution and existing analytical solution. Further discussions are made to disclose possible defects of the proposed solution for objectivity.
\end{abstract}


\section{Introduction}
\label{sec:introduction}

The complex variable method \cite[]{Muskhelishvili1966} is an important analytical approach to study the mechanical behaviour of geomaterial in tunnel engineering. For a deep tunnel, no ground surface or gravitational gradient need to be considered, and the mechanical models of a noncircular cavity excavation in an infinite plane constrained by an orthogonal in-situ stress field can be established \cite[]{Lu2015,Lu2014,Wang2019deep_non-circular,fanghuangcheng2021,zhou2023stress}, where a simply-connected exterior region is conformally mapped onto an infinite mapping plane outside of a unit circle \cite[]{Muskhelishvili1966,ma2022numerical}.

For shallow tunnels, by contrast, the ground surface would further geometrically and mechanically constrain the stress and displacement in geomaterial, and more complicated mechanical models of a circular cavity excavation in an infinite half plane constrained by gravitational gradient should be constructed \cite[]{Strack2002phdthesis,Strack_Verruijt2002buoyancy,Verruijt_Strack2008buoyancy,fang15:_compl,Lu2016,Lu2019new_solution}, where a doubly-connected lower half plane containing a circular cavity is conformally mapped onto a unit annulus with certain inner radius using Verruijt's conformal mapping \cite[]{Verruijt1997displacement,Verruijt1997traction}. Extensively, a new conformal mapping that maps a lower half plane containing an axisymmetrical and noncircular cavity onto a unit annulus is developed by Zeng et al. \cite[]{Zengguisen2019} by adding finite coupled series items. Such an extensive conformal mapping technique is soon applied to several shallow axisymmetrical and noncircular cavity excavation problems \cite[]{lu2021complex,zeng2022analytical,zeng2023analytical,fan2024analytical,zhou2024analytical}. 

All these complex variable solutions mentioned above assume that the cavity boundary is ideally excavated to perfectly match the designed cavity cross section. However, an over-/under-break excavation is a common phenomenon in tunnel engineering  \cite[]{singh2005causes,mandal2009evaluating,mottahedi2018development,fodera2020factors,chen2021experimental}. An over-/under-break tunnelling can be geometrically catagorized as asymmetrical and noncircular cavity excavation, for which the Verruijt's or Zeng's extensive conformal mapping would fail, and no complex variable solution can be subsequently established for shallow tunnelling. Therefore, a robust and relatively accurate conformal mapping technique should be developed to map a lower half plane containing an asymmetrical and noncircular cavity.

Additionally, the displacement singularity at infinity \cite[]{Self2020JEM} that jeopardizes the reliability of the complex variable method is less considered in the above mentioned shallow tunnelling solutions. To eliminate the displacement singularity at infinity, new mechanical models that constrain far-field ground surface to equilibrate the local unbalanced resultant caused by gravitational geomaterial removal are constructed \cite[]{LIN2024106008,LIN2024396}, and the corresponding complex variable solutions \cite[]{LIN2024396,LIN2024106008} apply the analytic continuation principle to turn the mixed boundaries into homogenerous Riemann-Hilbert problems. However, the common over-/under-break phenomena are also not considered in these two newly proposed solutions \cite[]{LIN2024396,LIN2024106008}.

In this paper, the common over-/under-break cavity excavation phenomena would be considered in shallow tunnelling in gravitational geomaterial, and corresponding complex variable solution with consideration of eliminating the far-field displacement at infinity would be proposed by integrating suitable conformal mapping that is capable of dealing with a lower half plane containing an asymmetrical and noncircular cavity. To be more specific, this paper extends our previous solution in Ref \cite[]{LIN2024106008} by eliminating its defects.

\section{Problem definition}
\label{sec:problem-definition}

\subsection{Mechanical model of over/under-break excavation in gravitational geomaterial}
\label{sec:mech-model-mult}

\begin{figure}[htb]
  \centering
  \begin{tikzpicture}
    \tikzstyle{every node} = [scale = 0.7]
    \fill [gray!30] (0,0) rectangle (8,-6);
    \node at (0,-6) [above right] {$ {\bm{\varOmega}} $};
    \draw [cyan, line width = 1pt] (0,0) -- (8,0);
    \node at (4,0) [below,cyan] {Infinite free ground surface $ {\bm{C}}_{1} $};
    \draw [->] (4,0) -- (4.5,0) node [above] {$ x $};
    \draw [->] (4,0) -- (4,0.5) node [right] {$ y $};
    \node at (4,0) [above right] {$ O $};
    \draw [line width = 1pt, ->] (1,-4.5) -- (1,-5) node [left] {$ \gamma $};
    \foreach \x in {1,2,...,12} \draw [->] ({0-\x/12*0.5},{-\x/12*6}) -- (0,{-\x/12*6});
    \foreach \x in {1,2,...,12} \draw [->] ({8+\x/12*0.5},{-\x/12*6}) -- (8,{-\x/12*6});
    \foreach \x in {0,1,2,...,16} \draw [->] ({\x/16*8},-6.5) -- ({\x/16*8},-6);
    \node at (-0.5,-6) [below] {$ k_{0} \gamma y $};
    \node at (8.5,-6) [below] {$ k_{0} \gamma y $};
    \node at (4,-6.5) [below] {$ \gamma y $};
    \draw (0,0) -- (-0.5,-6);
    \draw (8,0) -- (8.5,-6);
    \draw (0,-6.5) -- (8,-6.5);
    \fill [gray!50] (4,-4) arc [start angle = -90, end angle = 90, x radius = 1.5, y radius = 1] arc [start angle = 90, end angle = 270, x radius = 1.2, y radius = 1]; 
    \draw [dashed, violet, line width = 1pt] (4,-4) arc [start angle = -90, end angle = 90, x radius = 1.5, y radius = 1] arc [start angle = 90, end angle = 270, x radius = 1.2, y radius = 1]; 
    \foreach \x in {1,2,3,4,5} \draw [<-] ({4+1.5*cos(-90+180/6*\x)},{-3+1*sin(-90+180/6*\x)}) -- ({4+1.5*cos(-90+180/6*\x)-0.2},{-3+1*sin(-90+180/6*\x)});
    \foreach \x in {1,2,3,4,5} \draw [<-] ({4+1.2*cos(90+180/6*\x)},{-3+1*sin(90+180/6*\x)}) -- ({4+1.2*cos(90+180/6*\x)+0.2},{-3+1*sin(90+180/6*\x)});
    \foreach \x in {1,2,3} \draw [<-] ({4+1.5*cos(180/6*\x)},{-3+1*sin(180/6*\x)}) -- ({4+1.5*cos(180/6*\x)},{-3+1*sin(180/6*\x)-0.2});
    \foreach \x in {3,4,5} \draw [<-] ({4+1.2*cos(180/6*\x)},{-3+1*sin(180/6*\x)}) -- ({4+1.2*cos(180/6*\x)},{-3+1*sin(180/6*\x)-0.2});
    \foreach \x in {1,2,3} \draw [<-] ({4+1.2*cos(180+180/6*\x)},{-3+1*sin(180+180/6*\x)}) -- ({4+1.2*cos(180+180/6*\x)},{-3+1*sin(180+180/6*\x)+0.2});
    \foreach \x in {3,4,5} \draw [<-] ({4+1.5*cos(180+180/6*\x)},{-3+1*sin(180+180/6*\x)}) -- ({4+1.5*cos(180+180/6*\x)},{-3+1*sin(180+180/6*\x)+0.2});
    \node at (4,-4) [below,violet] {$ {\bm{C}}_{2} $};
    \node at (4,-3) {$ {\bm{D}} $};
    \node at (4,-7) [below] {(a) Geomaterial $ {\bm{\varOmega}} $ with geostatic equilibrium and tractions along cavity periphery};

    \node at (9,-3) [left] {\LARGE{\texttt{+}}};

    \fill [gray!30] (9,0) rectangle (17,-6);
    \node at (9,-6) [above right] {$ {\bm{\varOmega}} $};
    \draw [cyan, line width = 1pt] (9,0) -- (17,0);
    \draw [->] (13,0) -- (13.5,0) node [above] {$ x $};
    \draw [->] (13,0) -- (13,0.5) node [right] {$ y $};
    \node at (13,0) [above right] {$ O $};
    \node at (13,0) [below,cyan] {Infinite free ground surface $ {\bm{C}}_{1} $};
    \fill [white] (13,-4) arc [start angle = -90, end angle = 90, x radius = 1.5, y radius = 1] arc [start angle = 90, end angle = 270, x radius = 1.2, y radius = 1]; 
    \draw [violet, line width = 1pt] (13,-4) arc [start angle = -90, end angle = 90, x radius = 1.5, y radius = 1] arc [start angle = 90, end angle = 270, x radius = 1.2, y radius = 1]; 
    \foreach \x in {1,2,3,4,5} \draw [->] ({13+1.5*cos(-90+180/6*\x)},{-3+1*sin(-90+180/6*\x)}) -- ({13+1.5*cos(-90+180/6*\x)-0.2},{-3+1*sin(-90+180/6*\x)});
    \foreach \x in {1,2,3,4,5} \draw [->] ({13+1.2*cos(90+180/6*\x)},{-3+1*sin(90+180/6*\x)}) -- ({13+1.2*cos(90+180/6*\x)+0.2},{-3+1*sin(90+180/6*\x)});
    \foreach \x in {1,2,3} \draw [->] ({13+1.5*cos(180/6*\x)},{-3+1*sin(180/6*\x)}) -- ({13+1.5*cos(180/6*\x)},{-3+1*sin(180/6*\x)-0.2});
    \foreach \x in {3,4,5} \draw [->] ({13+1.2*cos(180/6*\x)},{-3+1*sin(180/6*\x)}) -- ({13+1.2*cos(180/6*\x)},{-3+1*sin(180/6*\x)-0.2});
    \foreach \x in {1,2,3} \draw [->] ({13+1.2*cos(180+180/6*\x)},{-3+1*sin(180+180/6*\x)}) -- ({13+1.2*cos(180+180/6*\x)},{-3+1*sin(180+180/6*\x)+0.2});
    \foreach \x in {3,4,5} \draw [->] ({13+1.5*cos(180+180/6*\x)},{-3+1*sin(180+180/6*\x)}) -- ({13+1.5*cos(180+180/6*\x)},{-3+1*sin(180+180/6*\x)+0.2});
    \node at (13,-4) [below,violet] {$ {\bm{C}}_{2} $};
    \node at (13,-7) [below] {(b) Inverse tractions along cavity periphery to simulate excavation};
  \end{tikzpicture}
  \caption{Shallow asymmetrical cavity excavation in gravitational geomaterial}
  \label{fig:1}
\end{figure}
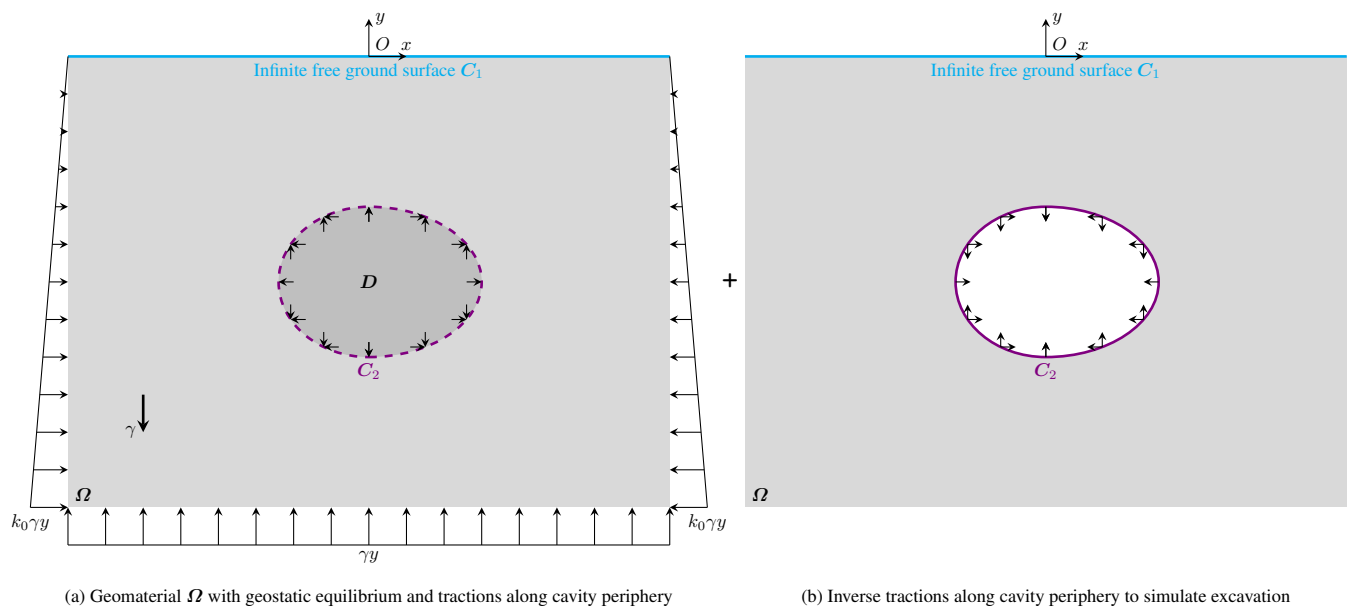

We may start from a geomaterial of gravitational gradient $ \gamma $ occupying the whole lower half complex plane, as shown in Fig. \ref{fig:1}a, and the lateral coefficient is $ k_{0} $. A global and complex coordinate system $ z = x + {\rm{i}} y $ is set ($ {\rm{i}} $ denotes the imaginarty unit with $ {\rm{i}}^{2} = -1 $), and the initial stress field in the geomaterial can be expressed as:
\begin{equation}
  \label{eq:2.1}
  \left\{
    \begin{aligned}
      & \sigma_{x}^{0} = k_{0} \gamma y \\
      & \sigma_{y}^{0} = \gamma y \\
      & \tau_{xy}^{0} = 0
    \end{aligned}
  \right.
\end{equation}
where $ \sigma_{x}^{0} $, $ \sigma_{y}^{0} $, and $ \tau_{xy}^{0} $ denote horizontal, vertical, and shear stress components of the initial stress field, respectively. With the initial stress field, the geomaterial has reached geostatic equilibrium. The ground surface is free, and can be denoted by $ {\bm{C}}_{1} $.

An over/under-break excavavation is subsequently to be conducted in the gravitational geomaterial, as shown in Fig. \ref{fig:1}a. The cavity boundary is denoted by $ {\bm{C}}_{2} $, and an arbitrary point along $ {\bm{C}}_{2} $ is denoted by $ S $. The open region within boundary $ {\bm{C}}_{2} $ is denoted by $ {\bm{D}} $, while the rest open region is denoted by $ {\bm{\varOmega}} $. To be concise, we denote $ \overline{{\bm{D}}} = {\bm{D}} \cup {\bm{C}}_{2} $, while $ \overline{\bm{\varOmega}} = {\bm{C}}_{2} \cup {\bm{\varOmega}} \cup {\bm{C}}_{1} $. Some certain point $ z_{c} = x_{c} + {\rm{i}}y_{c} $ within boundary $ {\bm{C}}_{2} $ is selected as the center of closed region $ \overline{\bm{D}} $. The induced orientation of boundary $ {\bm{C}}_{2} $ is specified to be clockwise to always keep the remaining open region $ {\bm{\varOmega}} $ on the left side of boundary $ {\bm{C}}_{2} $, while the outward normal direction of boundary $ {\bm{C}}_{2} $ is denoted by $ \vec{n} $, which always points to the inside of boundary $ {\bm{C}}_{2} $. The exacavtion is conducted in Fig. \ref{fig:1}b, and the following tractions are applied along boundary $ {\bm{C}}_{2} $ to equilibrate the initial stress field:
\begin{equation}
  \label{eq:2.2}
  \left\{
    \begin{aligned}
      & X(S) = - X^{0}(S) = - \sigma_{x}^{0}(S) \cdot \cos\langle \vec{n}, x \rangle - \tau_{xy}^{0}(S) \cdot \cos\langle \vec{n}, y \rangle \\
      & Y(S) = - Y^{0}(S) = - \tau_{xy}^{0}(S) \cdot \cos\langle \vec{n}, x \rangle - \sigma_{y}^{0}(S) \cdot \cos\langle \vec{n}, y \rangle \\
    \end{aligned}
  \right.
\end{equation}
where $ X(S) $ and $ Y(S) $ denote horizontal and vertical inverse surface tractions of arbitrary point $ S $ along boundary $ {\bm{C}}_{2} $ in Fig \ref{fig:1}b, $ X^{0}(S) $ and $ Y^{0}(S) $ denote horizontal and vertical surface tractions of arbitrary point $ S $ along boundary $ {\bm{C}}_{2} $ in Fig \ref{fig:1}a, $ \sigma_{x}^{0}(S) $, $ \sigma_{y}^{0}(S) $, and $ \tau_{xy}^{0}(S) $ denote horizontal, vertical, and shear stress components of the initial stress field of arbitrary point $ S $ along boundary $ {\bm{C}}_{2} $, $ \langle \vec{n}, x \rangle $ and $ \langle \vec{n}, y \rangle $ denote the angles between the outward normal direction $ \vec{n} $ of boundary $ {\bm{C}}_{2} $ and the positive directions of $x$ and $y$ axes, respectively. If the path differential segment of boundary $ {\bm{C}}_{2} $ is denoted by $ {\rm{d}}S $, the following mathematical expression $ \cos\langle \vec{n}, x \rangle = \frac{{\rm{d}}y}{{\rm{d}}S} $ and $ \cos\langle \vec{n}, y \rangle = - \frac{{\rm{d}}x}{{\rm{d}}S} $ always stand, according the induced orientation and outward normal direction specified above.

The resultant along boundary $ {\bm{C}}_{2} $ due to cavity excavation can then be expressed as:
\begin{equation}
  \label{eq:2.3}
  \left\{
    \begin{aligned}
      & R_{x} = \varointclockwise_{{\bm{C}}_{2}} X(S) |{\rm{d}}S| = k_{0} \gamma \varointclockwise_{{\bm{C}}_{2}} y {\rm{d}}y = - k_{0} \gamma \ointctrclockwise_{{\bm{C}}_{2}} \left( 0{\rm{d}}x + y {\rm{d}}y \right) = k_{0}\gamma \iint_{{\bm{D}}} 0 \cdot {\rm{d}}x{\rm{d}}y = 0 \\
      & R_{y} = \varointclockwise_{{\bm{C}}_{2}} Y(S) |{\rm{d}}S| = \gamma \varointclockwise_{{\bm{C}}_{2}} y {\rm{d}}x = - \gamma \ointctrclockwise_{{\bm{C}}_{2}} \left( y{\rm{d}}x + 0 {\rm{d}}y \right) = \gamma \iint_{{\bm{D}}} 1 \cdot {\rm{d}}x{\rm{d}}y = 2 \pi W \gamma
    \end{aligned}
  \right.
\end{equation}
where $ R_{x} $ and $ R_{y} $ denote the horizontal and vertical traction resultant components for boundary $ {\bm{C}}_{2} $, respectively, $ |{\rm{d}}S| = -{\rm{d}}S $ for clockwise integration, and the fourth equation applies the Green theorem for the simply connected region $ \overline{{\bm{D}}} $ with the induced orientation (counter-clockwise), and $ 2\pi W $ denotes the area of the excavated region $ \overline{{\bm{D}}} $. Eq. (\ref{eq:2.3}) indicates that the traction resultant along boundary $ {\bm{C}}_{2} $ only contains a vertically upward component and free from any horizontal one for whatever cavity shape.

\begin{figure}[htb]
  \centering
  \begin{tikzpicture}
    \tikzstyle{every node} = [scale = 0.7]
    \fill [gray!30] (0,0) rectangle (8,-6);
    \node at (0,-6) [above right] {$ {\bm{\varOmega}} $};
    \draw [cyan, line width = 1pt] (0,0) -- (8,0);
    \draw [->] (4,0) -- (4.5,0) node [above] {$ x $};
    \draw [->] (4,0) -- (4,0.5) node [right] {$ y $};
    \node at (4,0) [above right] {$ O $};
    \node at (4,0) [below,cyan] {Infinite free ground surface $ {\bm{C}}_{1} $};
    \fill [white] (4,-4) arc [start angle = -90, end angle = 90, x radius = 1.5, y radius = 1] arc [start angle = 90, end angle = 270, x radius = 1.2, y radius = 1]; 
    \draw [violet, line width = 1pt] (4,-4) arc [start angle = -90, end angle = 90, x radius = 1.5, y radius = 1] arc [start angle = 90, end angle = 270, x radius = 1.2, y radius = 1]; 
    \draw [red, ->] (4,-3) -- (4,-2.5) node [right] {$ R_{y} $};
    \fill [violet] (4,-3) circle [radius = 0.05];
    \node at (4,-3) [right, violet] {$ z_{c} $};
    \node at (4,-4) [below,violet] {$ {\bm{C}}_{2} $};
    \node at (4,-6) [below] {(a) Nonzero resultant along cavity boundary};

    \fill [gray!30] (9,0) rectangle (17,-6);
    \node at (9,-6) [above right] {$ {\bm{\varOmega}} $};
    \draw [cyan, line width = 1pt] (9.5,0) -- (16.5,0);
    \draw [line width = 1pt] (9,0) -- (9.5,0);
    \fill [pattern = north east lines] (9,0) rectangle (9.5,0.2);
    \draw [line width = 1pt] (17,0) -- (16.5,0);
    \fill [pattern = north east lines] (17,0) rectangle (16.5,0.2);
    \node at (13,0.7) [above] {Infinite constrained far-field ground surface $ {\bm{C}}_{1c} $};
    \draw [<-] (9.5,0.2) -- (10.5,0.7);
    \draw [<-] (16.5,0.2) -- (15.5,0.7); 
    \draw [->] (13,0) -- (13.5,0) node [above] {$ x $};
    \draw [->] (13,0) -- (13,0.5) node [right] {$ y $};
    \node at (13,0) [above right] {$ O $};
    \node at (13,0) [below,cyan] {Finite free ground surface $ {\bm{C}}_{1f} $};
    \fill [white] (13,-4) arc [start angle = -90, end angle = 90, x radius = 1.5, y radius = 1] arc [start angle = 90, end angle = 270, x radius = 1.2, y radius = 1]; 
    \draw [violet, line width = 1pt] (13,-4) arc [start angle = -90, end angle = 90, x radius = 1.5, y radius = 1] arc [start angle = 90, end angle = 270, x radius = 1.2, y radius = 1]; 
    \draw [red, ->] (13,-3) -- (13,-2.5) node [right] {$ R_{y} $};
    \fill [violet] (13,-3) circle [radius = 0.05];
    \node at (13,-3) [right, violet] {$ z_{c} $};
    \node at (13,-4) [below,violet] {$ {\bm{C}}_{2} $};
    \fill [red] (9.5,0) circle [radius = 0.05];
    \fill [red] (16.5,0) circle [radius = 0.05];
    \node at (9.5,0) [below, red] {$ T_{1} $};
    \node at (16.5,0) [below, red] {$ T_{2} $};
    \node at (13,-6) [below] {(b) Geomaterial with constrained far-field ground surface};
  \end{tikzpicture}
  \caption{New mechanical model construction with constrained far-field ground surface}
  \label{fig:2}
\end{figure}
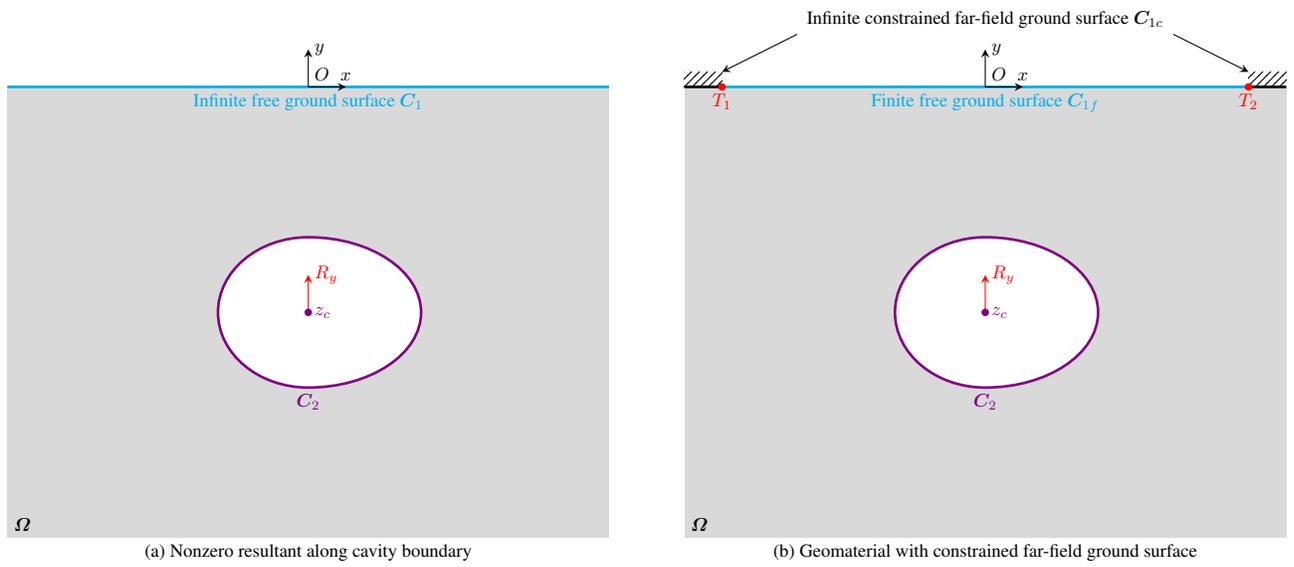

Eq. (\ref{eq:2.3}) indicates that the traction resultant along the asymmetrical cavity boundary is equal to the opposite gravitational weight of the excavated geomaterial, as shown in Fig. \ref{fig:2}a. The location point of the vertical nonzero resultant along boundary $ {\bm{C}}_{2} $ is set to the center $ z_{c} = x_{c} + {\rm{i}} y_{c} $ of closure region $ \overline{{\bm{D}}} $. The nonzero resultant $ Y $ indicates that no static equilibrium is established in the mechanical model in Fig. \ref{fig:2}a.

The following complex potentials $ \varphi(z) $ and $ \psi(z) $ can be established based on the mechanical model in Fig. \ref{fig:2}a according to Strack \cite[]{Strack2002phdthesis}:
\begin{subequations}
  \label{eq:2.4}
  \begin{equation}
    \label{eq:2.4a}
    \varphi(z) = -\frac{{\rm{i}} \kappa R_{y}} {2\pi(1+\kappa)} \ln (z-\overline{z}_{c}) - \frac{{\rm{i}}R_{y}} {2\pi(1+\kappa)} \ln(z-z_{c}) + \varphi^{\ast}_{0}(z), \quad z \in \overline{{\bm{\varOmega}}}
  \end{equation}
  \begin{equation}
    \label{eq:2.4b}
    \psi(z) = -\frac{{\rm{i}} R_{y}} {2\pi(1+\kappa)} \ln (z-\overline{z}_{c}) - \frac{{\rm{i}} \kappa R_{y}} {2\pi(1+\kappa)} \ln(z-z_{c}) + \psi^{\ast}_{0}(z), \quad z \in \overline{{\bm{\varOmega}}}
  \end{equation}
\end{subequations}
where $ \kappa $ is the Kolosov parameter with $ \kappa = 3-4\nu $ for plane strain and $ \kappa = \frac{3-\nu}{1-\nu} $ for plane stress, $ \varphi^{\ast}_{0}(z) $ and $ \psi^{\ast}_{0}(z) $ denote the single-valued items. The above complex potentials can be further transformed as
\begin{equation}
  \label{eq:2.4a'}
  \tag{2.4a'}
  \varphi(z) = -\frac{{\rm{i}}R_{y}} {2\pi} \ln (z-\overline{z}_{c}) + \frac{{\rm{i}}R_{y}} {2\pi(1+\kappa)} \ln\frac{z-\overline{z}_{c}}{z-z_{c}} + \varphi^{\ast}_{0}(z), \quad z \in \overline{{\bm{\varOmega}}}
\end{equation}
\begin{equation}
  \label{eq:2.4b'}
  \tag{2.4b'}
  \psi(z) = -\frac{{\rm{i}}R_{y}} {2\pi} \ln (z-\overline{z}_{c}) + \frac{{\rm{i}} \kappa R_{y}} {2\pi(1+\kappa)} \ln\frac{z-\overline{z}_{c}}{z-z_{c}} + \psi^{\ast}_{0}(z), \quad z \in \overline{{\bm{\varOmega}}}
\end{equation}

The nonzero resultant specified by Eqs. (\ref{eq:2.4a'}) and (\ref{eq:2.4b'}) would further lead to a unique displacement singularity. The first items in Eqs. (\ref{eq:2.4a'}) and (\ref{eq:2.4b'}) indicate that an upward concentrated force $ R_{y} $ is located at the point $ \overline{z}_{c} $ along the boundary of the infinitely small circle centered at point $ \overline{z}_{c} $. Since the circle is infinitely small, the point along the boundary of this circle and the center of this circle itself overlap to each other. Each of the second items in Eqs. (\ref{eq:2.4a'}) and (\ref{eq:2.4b'}) indicates that a pair of concentrated forces of the same magnitude $ R_{y} $ the located points $ \overline{z}_{c} $ and $ z_{c} $ are in the same vertical line but of opposite directions, which can be cancelled by each other. Traction continuation has been implicitly used in Eqs. (\ref{eq:2.4a'}) and (\ref{eq:2.4b'}), and similar analyses can be found in the prior studies \cite[]{LIN2024396,LIN2024106008}. With the complex potentials, the stress and displacement in geomaterial can be computed as \cite[]{Muskhelishvili1966}
\begin{subequations}
  \label{eq:2.5}
  \begin{equation}
    \label{eq:2.5a}
    \left\{
      \begin{aligned}
        & \sigma_{y} + \sigma_{x} = 2 [\varphi^{\prime}(z) + \overline{\varphi^{\prime}(z)}] \\
        & \sigma_{y} - \sigma_{x} + 2{\rm{i}}\tau_{xy} = 2[\overline{z}\varphi^{\prime\prime}(z) + \psi^{\prime}(z)] \\
      \end{aligned}
    \right.
    , \quad z \in \overline{{\bm{\varOmega}}}
  \end{equation}
  \begin{equation}
    \label{eq:2.5b}
    2G(u+{\rm{i}}v) = \kappa \varphi(z) - z\overline{\varphi^{\prime}(z)} - \overline{\psi(z)}, \quad z \in \overline{{\bm{\varOmega}}}
  \end{equation}
\end{subequations}
where $ \sigma_{x} $, $ \sigma_{y} $, and $ \tau_{xy} $ denote horizontal, vertical, and shear stress components in the remaining geomaterial $ \overline{{\bm{\varOmega}}} $ due to cavity exacavtion, respectively; $ u $ and $ v $ denote horizontal and vertical displacement components in the remaining geomaterial $ \overline{{\bm{\varOmega}}} $ due to cavity exacavtion, respectivley; $ G = \frac{E}{2(1+\nu)} $ denotes shear modulus of the geomaterial; $ \bullet^{\prime} $ and $ \bullet^{\prime\prime} $ denote the first and second derivatives of $ \bullet $, respectively; $ \overline{\bullet} $ denotes conjugate of $ \bullet $. Apparently, when substituting Eqs. (\ref{eq:2.4a'}) and (\ref{eq:2.4b'}) into Eq. (\ref{eq:2.5}), a unique displacement singularity at infinity ($ z \rightarrow \infty $) would exist, and similar deduction of this unique displacement singularity has been reported in Refs \cite[]{LIN2024396,LIN2024106008}, and would not be repeated here.

In summary, the above analyses indicate that the mechanical model in Fig. \ref{fig:2}a contain a nonzero principal vector, and would further lead to a unique displacement at infinity, which is not identical to asymmetrical cavity excavation in real-world engineering, and should be eliminated accordingly. To eliminate the displacement singularity, we only need to constrain the far-field ground surface to equilibrate the nonzero principal vector, as shown in Fig. \ref{fig:2}b, and the displacement singularity at infinity would vanish spontaneously. Correspondingly, the ground surface is separated into the still free but finite segment above the multiple cavities $ {\bm{C}}_{1f} $ and the constrained and infinite segment of far-field surface $ {\bm{C}}_{1c} $. The joint points of these two segments are denoted by $ T_{1} $ and $ T_{2} $, respectively. Similar models can be found in Refs \cite[]{LIN2024396,LIN2024106008}. No matter how the nonzero resultant change in the mechanical model in Fig. \ref{fig:2}b, the constrained far-field ground surface would generate necessary constraining forces to satisfy static equilibrium. Therefore, the former unbalanced mechanical model turns to a new and equilibrium one.

\subsection{Mixed boundary conditions of new mechanical model}
\label{sec:boundary-conditions}

In the new and equilibrium mechanical model in Fig. \ref{fig:2}b, the ground surface is separated into two segments $ {\bm{C}}_{1c} $ and $ {\bm{C}}_{1f} $. The former is constrained, indicating no displacement should occur along boundary $ {\bm{C}}_{1c} $; and the latter is free, indicating no surface traction should be applied along boundary $ {\bm{C}}_{1f} $. Thus, the following explicit boundary conditions should be applied along ground surface:
\begin{subequations}
  \label{eq:2.6}
  \begin{equation}
    \label{eq:2.6a}
    u(T) + {\rm{i}}v(T) = 0, \quad T \in {\bm{C}_{1c}}
  \end{equation}
  \begin{equation}
    \label{eq:2.6b}
    X(T) + {\rm{i}}Y(T) = 0, \quad T \in {\bm{C}_{1f}}
  \end{equation}
\end{subequations}
where $ X(T) $ and $ Y(T) $ denote horizontal and vertical tractions along ground surface, respectively; $ u(T) $ and $ v(T) $ denote horizontal and vertical displacement along ground surface, respectively; $ T $ denotes some point along ground surface. The boundary condition in Eq. (\ref{eq:2.6a}) introduces an infinite cut in the lower half complex plane, which starts from point $ T_{1} $, and extends horizontally and left towards infinity, and then further extends horizontally and left towards point $ T_{2} $. Apparently, Eq. (\ref{eq:2.6}) contains mixed boundary conditions.

Another explicit boundary conditions in the mechanical model in Fig. \ref{fig:2} is the inverse traction to simulate excavation in Eq. (\ref{eq:2.2}), which can be rewritten as
\begin{equation}
  \label{eq:2.7}
  X(S) + {\rm{i}}Y(S) = - X^{0}(S) - {\rm{i}}Y^{0}(S), \quad S \in {\bm{C}}_{2}
\end{equation}
The boundary condition in Eq. (\ref{eq:2.7}) is first-kind.

The static equilibrium in the remaining geomaterial $ \overline{{\bm{\varOmega}}} $ in the mechanical model in Fig. \ref{fig:2}b requires that the unbalanced resultants of the inverse traction along the cavity boundary should be equilibrated by the constrained far-field ground surface, and the following implicit equilibrium should be satisfied \cite[]{LIN2024396,LIN2024106008}:
\begin{equation}
  \label{eq:2.8}
  \int_{{\bm{C}}_{1c}} \left[ X(T) + {\rm{i}} Y(T) \right] {\rm{d}}T = - R_{x} - {\rm{i}}R_{y} = -{\rm{i}}R_{y}
\end{equation}
The mechanical model in Fig. \ref{fig:2}b further requires the single-valuedness of displacement in the remaining geomaterial $ \overline{{\bm{\varOmega}}} $, which cannot be expressed in formula for now and will be discussed below.

In summary, we should solve the stress and displacement solution in the remaining geomaterial $ \overline{{\bm{\varOmega}}} $ using the explicit boundary conditions in Eqs. (\ref{eq:2.6}) and (\ref{eq:2.7}), as well as the implicit equilibrium in Eq. (\ref{eq:2.8}) and the displacement single-valuedness.

\section{Bidirectional composite conformal mapping and geomaterial simulation}
\label{sec:appr-comp-conf}

\subsection{Bidirectional composite conformal mapping}
\label{sec:bidir-comp-conf-1}

As shown in Figs. \ref{fig:1} and \ref{fig:2}, a shallow over/under-break excavtion in geomaterial would probably cause an asymmetrical cavity in an infinite lower half plane. To solve the stress and displacement in the remaining geomaterial $ \overline{{\bm{\varOmega}}} $ using the complex variable method, a proper conformal mapping that bidirectionally maps the remaining geomaterial $ \overline{{\bm{\varOmega}}} $ and a corresponding unit annulus onto each other should be constructed. However, as far as we know, no such direct conformal mapping has been constructed. To overcome such a difficulty, a two-step composite conformal mapping consisting of two existing conformal mappings with certain constraint is applied to approximately simulate the necessary bidirectional conformal mapping.

The first step of the composite conformal mapping is to bidirectionally map the geomaterial occupying the remaining infinite lower half plane containing an asymmetrical cavity $ \overline{\bm{\varOmega}} $ in the complex plane $ z = x + {\rm{i}} y $ and a quasi lower half plane containing a unit circular cavity $ \overline{\bm{\varOmega}}^{\prime} $ in complex plane $ w = \Re w + {\rm{i}}\Im w $ onto each other, as shown in Figs. \ref{fig:3}a-1 and a-2, and the forward and backward conformal mappings are respectively denoted as
\begin{subequations}
  \label{eq:3.1}
  \begin{equation}
    \label{eq:3.1a}
    w = w(z)
  \end{equation}
  \begin{equation}
    \label{eq:3.1b}
    z = z(w)
  \end{equation}
\end{subequations}
where the detailed expressions can be found in Appendix \ref{sec:first-step-composite}. In Fig. \ref{fig:3}a-1, the blue dashed line indicates the expected excavation cavity boundary, with which an expected axisymmetrical line could be established. However, the real excavation cavity boundary is asymmetrical, thus, such axisymmetry would no longer exist. Via the conformal mapping in Eq. (\ref{eq:3.1}), the cavity boundary $ {\bm{C}}_{2} $, upper boundaries $ {\bm{C}}_{1c} $ and $ {\bm{C}}_{1f} $, and joint points $ T_{1} $ and $ T_{2} $ in physical plane $ z = x + {\rm{i}}y $ in Fig. \ref{fig:3}a-1 are bidirectionally mapped onto the cavity boundary $ {\bm{C}}_{2}^{\prime} $, upper boundaries $ {\bm{C}}_{1c}^{\prime} $ and $ {\bm{C}}_{1f}^{\prime} $, and joint points $ T_{1}^{\prime} $ and $ T_{2}^{\prime} $ in mapping plane $ w = \Re w + {\rm{i}}\Im w $ in Fig. \ref{fig:3}a-2, respectively. It should addressed that the first-step bidirectional conformal mapping only focuses on the cavity boundary mapping, and the upper boundary of the mapping geomaterial in Fig. \ref{fig:3}a-2 is no longer a straight line, but is slightly curved due to the forward conformal mapping $ w = w(z) $, as long as the following geometrical constraint is applied that the depth of the real cavity should {\em{not}} be too close to the upper boundary. Such a geometrical constraint is discussed in the numerical cases in Section .

The second step of the composite conformal mapping is to further bidirectionally map the geomaterial occupying the quasi lower half plane $ \overline{\bm{\varOmega}}^{\prime} $ in complex plane $ w = \Re w + {\rm{i}}\Im w $ and a quasi unit annulus $ \overline{\bm{\omega}} $ in complex plane $ \zeta = \rho \cdot {\rm{e}}^{{\rm{i}}\theta} $ onto each other, as shown in Figs. \ref{fig:3}a-2 and a-3, and the forward and backward conformal mappings are respectively denoted as
\begin{subequations}
  \label{eq:3.2}
  \begin{equation}
    \label{eq:3.2a}
    \zeta = \zeta(w)
  \end{equation}
  \begin{equation}
    \label{eq:3.2b}
    w = w(\zeta)
  \end{equation}
\end{subequations}
where the detailed expressions can be found in Appendix \ref{sec:second-step-comp}. Via the conformal mapping in Eq. (\ref{eq:3.2}), the cavity boundary $ {\bm{C}}_{2}^{\prime} $, upper boundaries $ {\bm{C}}_{1c}^{\prime} $ and $ {\bm{C}}_{1f}^{\prime} $, and joint points $ T_{1}^{\prime} $ and $ T_{2}^{\prime} $ in mapping plane $ w = \Re w + {\rm{i}}\Im w $ in Fig. \ref{fig:3}a-2 are bidirectionally mapped onto the inner periphery $ {\bm{c}}_{2} $, outer peripheries $ {\bm{c}}_{1c} $ and $ {\bm{c}}_{1f} $, and joint points $ t_{1} $ and $ t_{2} $ in mapping plane $ \zeta = \rho \cdot {\rm{e}}^{{\rm{i}}\theta} $ in Fig. \ref{fig:3}a-3, respectively. The polar angles of joint points $ t_{1} $ and $ t_{2} $ are denoted by $ \theta_{1} $ and $ \theta_{2} $, respectively. It should be addressed that the slightly curved upper boundary in Fig. \ref{fig:3}a-2 would cause that (1) the cavity depth $ h $ in Fig. \ref{fig:3}a-2 should be recalculated according to the imaginary part of joint point $ T_{2} $; (2) the outer periphery of the quasi unit annulus in Fig. \ref{fig:3}a-3 is a quasi unit circle.

Therefore, the forward and backward conformal mappings of the bidrectional two-step composite conformal mapping can be respectively expressed as
\begin{subequations}
  \label{eq:3.3}
  \begin{equation}
    \label{eq:3.3a}
    \zeta(z) = \zeta[w(z)]
  \end{equation}
  \begin{equation}
    \label{eq:3.3b}
    z(\zeta) = z[w(\zeta)]
  \end{equation}
\end{subequations}
Via the bidirectional two-step composite conformal mapping in Eq. (\ref{eq:3.3}), the geomaterial in a lower half plane containing an asymmetrical cavity in physical plane $ z = x + {\rm{i}}y $ and a quasi unit annulus in mapping plane $ \zeta = \rho \cdot {\rm{e}}^{{\rm{i}}\theta} $ can be mutually mapped onto each other.

\subsection{Geomaterial simulation}
\label{sec:geom-resel}

The purpose of the bidirectional conformal mapping in Eq. (\ref{eq:3.3}) is to facilitate the construction of Muskhelishvili's complex potentials in Laurent series formation \cite[]{Muskhelishvili1966} in the mapped geomaterial $ \overline{\bm{\omega}} $ in Fig. \ref{fig:3}a-3 to further solve stress and displacement in geomaterial via Eq. (\ref{eq:2.5}), as long as region $ \overline{\bm{\omega}} $ is a real unit annulus. However, region $ \overline{\bm{\omega}} $ is not a real unit annulus, and the outer periphery is just a quasi unit circle. Consequently and conceptually, the Muskhelishvili's complex potentials can not be expanded in Laurent series formation in the mapped geomaterial in Fig. \ref{fig:3}a-3. Given that the Laurent series expanding formation is very efficient to solve the Muskhelishvili's complex potential, certain modification and approximation of the quasi annulus may be conducted. 

We may recall that a certain geometrical constraint has been set in the first step of composite conformal mapping that the cavity should {\emph{not}} be too close to the upper boundary in Fig. \ref{fig:3}a-1, so that after the composite foward conformal mapping in Eq. (\ref{eq:3.3a}), the outer periphery in Fig. \ref{fig:3}a-3 should be approximate to a unit circle, and the `radius` should be approximate to unit. Thus, it may be possible to reselect the mapped geomaterial in the quasi unit annulus region $ \overline{\bm{\omega}} $ in Fig. \ref{fig:3}a-3 to a real unit annulus $ \overline{\bm{\omega}}_{0} $ in Fig. \ref{fig:3}b-3. In other words, we attempt to use the unit annulus region $ \overline{\bm{\omega}}_{0} $ in Fig. \ref{fig:3}b-3 to geometrically simulate the mapping quasi unit annulus region $ \overline{\bm{\omega}} $ in Fig. \ref{fig:3}a-3. Correspondingly, the inner periphery $ {\bm{c}}_{2} $, outer peripheries $ {\bm{c}}_{1c} $ and $ {\bm{c}}_{1f} $, and joint points $ t_{1} $ and $ t_{2} $ in the mapping quasi unit annulus in Fig. \ref{fig:3}a-3 are geometrically simulated by the inner periphery $ {\bm{c}}_{20} $, outer peripheries $ {\bm{c}}_{1c0} $ and $ {\bm{c}}_{1f0} $, and joint points $ t_{10} $ and $ t_{20} $ in the mapping quasi unit annulus in Fig. \ref{fig:3}b-3. It should be addressed that the polar angles of joint points $ t_{10} $ and $ t_{20} $ in Fig. \ref{fig:3}b-3 are exactly the same to those of joint points $ t_{1} $ and $ t_{2} $ in Fig. \ref{fig:3}a-3.

Theoretically, such a geomaterial simulation scheme is mathematically bold, since it \emph{changes the definition domain of the mapped geomaterial}. As mentioned above that the `radius` of the outer periphery of the quasi unit annulus $ \overline{\bm{\omega}} $ is approximate to unit, thus, such a unit annulus geomaterial simulation scheme would not greatly alter the geometry of the geomaterial after the composite backward conformal mapping in Eq. (\ref{eq:3.3b}) (detailed discussion can be found in Section ). Fow now, we may simply assume that such a geomaterial simulation scheme is reasonable, and we should consequently expect that the second-step backward conformal mapping in Eq. (\ref{eq:3.2b}) would map the simulating unit annulus $ \overline{\bm{\omega}}_{0} $ in mapping plane $ \zeta = \rho \cdot {\rm{e}}^{{\rm{i}}\theta} $ in Fig. \ref{fig:3}b-3 onto a lower half plane $ \overline{\bm{\varOmega}}_{0}^{\prime} $ containing a unit circular cavity with a fully straight upper boundary in mapping plane $ w = \Re w + {\rm{i}}\Im w $ in Fig. \ref{fig:3}b-2. Correspondingly, the inner periphery $ {\bm{c}}_{20} $, outer peripheries $ {\bm{c}}_{1c0} $ and $ {\bm{c}}_{1f0} $, and joint points $ t_{10} $ and $ t_{20} $ in the mapping quasi unit annulus in Fig. \ref{fig:3}b-3 would be backwardly mapped onto the cavity boundary $ {\bm{C}}_{20}^{\prime} $, upper boundaries $ {\bm{C}}_{1c0}^{\prime} $ and $ {\bm{C}}_{1f0}^{\prime} $, and joint points $ T_{10}^{\prime} $ and $ T_{20}^{\prime} $ in mapping plane $ w = \Re w + {\rm{i}}\Im w $ in Fig. \ref{fig:3}b-2, respectively. All the components in Fig. \ref{fig:3}b-2 are simulations of corresponding ones in Fig. \ref{fig:3}a-2.

The simulating mapping geomaterial $ \overline{\bm{\varOmega}}_{0}^{\prime} $ in Fig. \ref{fig:3}b-2 could further be mapped by the first-step backward conformal mapping in Eq. (\ref{eq:3.1b}) onto a quasi lower half plane $ \overline{\bm{\varOmega}}_{0} $ containing an asymmetrical cavity with a slightly curved upper boundary in the physical plane $ z = x + {\rm{i}}y $. Correspondingly, the cavity boundary $ {\bm{C}}_{20}^{\prime} $, upper boundaries $ {\bm{C}}_{1c0}^{\prime} $ and $ {\bm{C}}_{1f0}^{\prime} $, and joint points $ T_{10}^{\prime} $ and $ T_{20}^{\prime} $ in mapping plane $ w = \Re w + {\rm{i}}\Im w $ in Fig. \ref{fig:3}b-2 are mapped onto the cavity boundary $ {\bm{C}}_{20} $, upper boundaries $ {\bm{C}}_{1c0} $ and $ {\bm{C}}_{1f0} $, and joint points $ T_{10} $ and $ T_{20} $ in physical plane $ z = x + {\rm{i}}y $ in Fig. \ref{fig:3}b-1, respectively. All the components in Fig. \ref{fig:3}b-1 are simulations of corresponding ones in Fig. \ref{fig:3}b-1. The shape of cavity boundary $ {\bm{C}}_{20} $ in Fig. \ref{fig:3}b-1 should be almost the same to the original cavity boundary $ {\bm{C}}_{2} $ in Fig. \ref{fig:3}a-1, and the upward resultant $ R_{y} $ can be thereby applied at point $ z_{c0} $ in Fig. \ref{fig:3}b-1, which is the simulation of resultant location $ z_{c} $ in Fig. \ref{fig:3}a-1. The upper boundaries $ {\bm{C}}_{1f0} $ and $ {\bm{C}}_{1c0} $ in Fig. \ref{fig:3}b-1 should be approximate to the upper boundaries $ {\bm{C}}_{1f} $ and $ {\bm{C}}_{1c} $ in Fig. \ref{fig:3}a-1.

With the geomaterial simulation scheme above, the Muskhelishvili's complex potentials in Laurent series formation can be constructed in the simulating geomaterial region $ \overline{\bm{\omega}}_{0} $, and can be further operated using the two-step backward conformal mapping.

\subsection{Simulation of boundary conditions and complex potentials}
\label{sec:simul-bound-cond}

Since the original geomaterial in Fig. \ref{fig:3}a-1 is simulated by the one in Fig. \ref{fig:3}b-1, the new mechanical model in Fig. \ref{fig:2}b would be replaced by the one in Fig. \ref{fig:3}b-1 accordingly. Correspondingly, the boundary conditions in Section \ref{sec:boundary-conditions} should also be simulated. To be specific, the mixed boundary conditions along ground surface $ {\bm{C}}_{1c} $ and $ {\bm{C}}_{1f} $ in Eq. (\ref{eq:2.6}) are respectively simulated by
\begin{equation}
  \label{eq:2.6a'}
  \tag{2.6a'}
  u_{0}(T_{0}) + {\rm{i}}v_{0}(T_{0}) = 0, \quad T_{0} \in {\bm{C}}_{1c0}
\end{equation}
\begin{equation}
  \label{eq:2.6b'}
  \tag{2.6b'}
  X_{0}(T_{0}) + {\rm{i}}Y_{0}(T_{0}) = 0, \quad T_{0} \in {\bm{C}}_{1f0}
\end{equation}
where $ u_{0}(T_{0}) $ and $ v_{0}(T_{0}) $ denote horizontal and vertical displacement along the simulated ground surface segment $ {\bm{C}}_{1c0} $ in Fig. \ref{fig:3}b-1, respectively; $ X_{0}(T_{0}) $ and $ Y_{0}(T_{0}) $ denote horizontal and vertical surface traction components along the simulated ground surface segment $ {\bm{C}}_{1f0} $ in Fig. \ref{fig:3}b-1, respectively; $ T_{0} $ denotes some point along the simulated ground surface in Fig. \ref{fig:3}b-1.

Similarly, Eqs. (\ref{eq:2.7}) and (\ref{eq:2.8}) can be simulated as
\begin{equation}
  \label{eq:2.7'}
  \tag{2.7'}
  X_{0}(S_{0}) + {\rm{i}}Y_{0}(S_{0}) = - X^{0}(S_{0}) - {\rm{i}}Y^{0}(S_{0}) = - X^{0}(S) - {\rm{i}}Y^{0}(S)
\end{equation}
\begin{equation}
  \label{eq:2.8'}
  \tag{2.8'}
  \int_{{\bm{C}}_{1c0}} [X_{0}(T_{0})+{\rm{i}}Y_{0}(T_{0})]{\rm{d}}T_{0} = -{\rm{i}}R_{y}
\end{equation}
where $ X_{0}(S_{0}) $ and $ Y_{0}(S_{0}) $ denote horizontal and vertical surface traction components along the simulated cavity boundary $ {\bm{C}}_{20} $ in Fig. \ref{fig:3}b-1, respectively; $ S_{0} $ denotes some point along the simulated cavity boundary $ {\bm{C}}_{20} $ in Fig. \ref{fig:3}b-1. The reason for the second equality of Eq. (\ref{eq:2.7'}) is that the shape of the simulated cavity boundary $ {\bm{C}}_{20} $ in Fig. \ref{fig:3}b-1 is almost the same to that of the original cavity boundary $ {\bm{C}}_{2} $ in Fig. \ref{fig:3}a-1, thus, each simulating point $ S_{0} $ along cavity boundary $ {\bm{C}}_{20} $ can be replaced by the original point $ S $ along cavity boundary $ {\bm{C}}_{2} $. The displacement single-valuedness in Fig. \ref{fig:2}b mentioned in Section \ref{sec:boundary-conditions} should also be inherited by Fig. \ref{fig:3}b-1.

Owing to the simulation of boundary conditions, the complex potentials in Eq. (\ref{eq:2.5}) should also be simulated as
\begin{equation}
  \label{eq:3.4}
  \left\{
    \begin{aligned}
      & \varphi_{0}(z) \sim \varphi(z) \\
      & \psi_{0}(z) \sim \psi(z)
    \end{aligned}
  \right.
\end{equation}
where $ \sim $ denotes the latter can be simulated by the former. The stress and displacement solution in the simulated geomaterial region $ \overline{\bm{\varOmega}}_{0} $ can be obtained by replacing the simulated complex potentials in Eq. (\ref{eq:3.4}) into Eq. (\ref{eq:2.5}). Therefore, the solution in this paper is a simulation one with tolerant errors.

\section{Mixed boundary problem transformation}
\label{sec:probl-transf}

With the simulating model in Fig. \ref{fig:3}b-1 and the two-step composite backward conformal mapping in Eq. (\ref{eq:3.3b}), the simulating mixed boundary conditions in Eq. (\ref{eq:2.6a'}), (\ref{eq:2.6b'}), and (\ref{eq:2.7'}) can be mapped from physical plane $ z $ onto mapping plane $ \zeta $ as
\begin{subequations}
  \label{eq:4.1}
  \begin{equation}
    \label{eq:4.1a}
    \left.[u_{0}(\zeta) + {\rm{i}} v_{0}(\zeta)]\right|_{\zeta \rightarrow t_{0}} = u_{0}(T_{0}) + {\rm{i}} v_{0}(T_{0}) = 0, \quad t_{0} \in {\bm{c}}_{1c0}
  \end{equation}
  \begin{equation}
    \label{eq:4.1b}
    \left. \frac{\zeta z^{\prime}(\zeta)}{|\zeta z^{\prime}(\zeta)|} [\sigma_{\rho0}(\zeta) + {\rm{i}}\tau_{\rho\theta0}(\zeta)] \right|_{\zeta \rightarrow t_{0}} = X_{0}(T_{0}) + {\rm{i}} Y_{0}(T_{0}) = 0, \quad t_{0} \in {\bm{c}}_{1f0} 
  \end{equation}
  \begin{equation}
    \label{eq:4.1c}
    \left. \frac{\zeta z^{\prime}(\zeta)}{|\zeta z^{\prime}(\zeta)|} [\sigma_{\rho0}(\zeta) + {\rm{i}}\tau_{\rho\theta0}(\zeta)] \right|_{\zeta \rightarrow s_{0}} = -k_{0}\gamma y \frac{{\rm{d}}y}{{\rm{d}}S_{0}} + {\rm{i}}\gamma y \frac{{\rm{d}}x}{{\rm{d}}S_{0}}, \quad s_{0} \in {\bm{c}}_{20}
  \end{equation}
\end{subequations}
where $ u_{0}(\zeta) $ and $ v_{0}(\zeta) $ denote horizontal and vertical displacement components in the simulating geomaterial $ \overline{\bm{\varOmega}}_{0} $ mapped onto the simulating unit annulus $ \overline{\bm{\omega}}_{0} $ in mapping plane $ \zeta $, respectively; $ \sigma_{\rho0}(\zeta) $ and $ \tau_{\rho\theta0}(\zeta) $ denote the curvilinear radial and shear stress components in the simulating geomaterial $ \overline{\bm{\varOmega}}_{0} $ mapped onto the simulating unit annulus $ \overline{\bm{\omega}}_{0} $ in mapping plane $ \zeta $, respectively; $ t_{0} $ and $ s_{0} $ denote arbitrary points along boundaries $ {\bm{c}}_{10} ({\bm{c}}_{10} = {\bm{c}}_{1c0} \cup {\bm{c}}_{1f0}) $ and $ {\bm{c}}_{20} $ in the simulating unit annulus $ \overline{\bm{\omega}}_{0} $ in mapping plane $ \zeta $, and are bijectively to $ T_{0} $ and $ S_{0} $, respectively. Considering the two-step composite backward conformal mapping in Eq. (\ref{eq:3.3b}), we have $ z^{\prime}(\zeta) = z^{\prime}(w) \cdot w^{\prime}(\zeta) $. For completion, the complementary expressions of Eqs. (\ref{eq:4.1a}) and (\ref{eq:4.1b}) are also given as
\begin{equation}
  \label{eq:4.1a'}
  \tag{4.1a'}
  \left.[u_{0}(\zeta) + {\rm{i}} v_{0}(\zeta)]\right|_{\zeta \rightarrow t_{0}} = u_{0}(T_{0}) + {\rm{i}} v_{0}(T_{0}) \neq 0, \quad t_{0} \in {\bm{c}}_{1f0}
\end{equation}
\begin{equation}
  \label{eq:4.1b'}
  \tag{4.1b'}
  \left. \frac{\zeta z^{\prime}(\zeta)}{|\zeta z^{\prime}(\zeta)|} [\sigma_{\rho0}(\zeta) + {\rm{i}}\tau_{\rho\theta0}(\zeta)] \right|_{\zeta \rightarrow t_{0}} = X_{0}(T_{0}) + {\rm{i}} Y_{0}(T_{0}) \neq 0, \quad t_{0} \in {\bm{c}}_{1c0}
\end{equation}

The following curvilinear complex potentials defined in the simulating unit annulus $ \overline{\bm{\omega}}_{0} $ are introduced to obtain the curvilinear stress and displacement components and to solve the mixed boundary problem in Eq. (\ref{eq:4.1}) as
\begin{subequations}
  \label{eq:4.2}
  \begin{equation}
    \label{eq:4.2a}
    \sigma_{\theta0}(\zeta) + \sigma_{\rho0}(\zeta)  = 2\left[ \varPhi_{0}(\zeta) + \overline{\varPhi_{0}(\zeta)} \right], \quad \zeta \in \overline{\bm{\omega}}_{0}
  \end{equation}
  \begin{equation}
    \label{eq:4.2b}
    \sigma_{\rho0}(\zeta) + {\rm{i}}\tau_{\rho\theta0}(\zeta) = \left[ \varPhi_{0}(\zeta) + \overline{\varPhi_{0}(\zeta)} \right] - \frac{\overline{\zeta}}{\zeta} \left[ \frac{z(\zeta)}{z^{\prime}(\zeta)} \overline{\varPhi^{\prime}_{0}(\zeta)} + \frac{\overline{z^{\prime}(\zeta)}}{z^{\prime}(\zeta)} \overline{\varPsi_{0}(\zeta)} \right], \quad \zeta \in \overline{\bm{\omega}}_{0}
  \end{equation}
  \begin{equation}
    \label{eq:4.2c}
    g_{0}(\zeta) = 2G[u_{0}(\zeta) + {\rm{i}}v_{0}(\zeta)] = \kappa \varphi_{0}(\zeta) - z(\zeta) \overline{\varPhi_{0}(\zeta)} - \overline{\psi_{0}(\zeta)}, \quad \zeta \in \overline{\bm{\omega}}_{0}
  \end{equation}
\end{subequations}
where
\begin{equation}
  \label{eq:4.3}
  \left\{
    \begin{aligned}
      & \varPhi_{0}(\zeta) = \frac{\varphi_{0}^{\prime}(\zeta)}{z^{\prime}(\zeta)} \\
      & \varPsi_{0}(\zeta) = \frac{\psi_{0}^{\prime}(\zeta)}{z^{\prime}(\zeta)}
    \end{aligned}
  \right.
\end{equation}

Using analytic continuation, the definition domain of $ \varPhi_{0}(\zeta) $ can be expanded from $ \overline{\bm{\omega}}_{0} $ to $ \overline{\bm{\omega}}_{0} \cup {\bm{\omega}}_{0}^{-} $ \cite[]{Muskhelishvili1966,LIN2024396,LIN2024106008}, where $ {\bm{\omega}}_{0}^{-} $ denotes the region outside of boundary $ {\bm{c}}_{10} $ of the simulating unit annulus, as shown in Fig. \ref{fig:3}b-3. Similarly, the open region $ {\bm{\omega}}_{0}^{+} $ inside of the boundary $ {\bm{c}}_{10} $ is defined as $ \overline{\bm{\omega}}_{0} = {\bm{\omega}}_{0}^{+} \cup {\bm{c}}_{10} $. Correspondingly, the other complex potential $ \psi_{0}^{\prime}(\zeta) $ can be expressed as \cite[]{LIN2024106008}
\begin{equation}
  \label{eq:4.4}
  \psi_{0}^{\prime}(\zeta) = \zeta^{-2} \overline{\varphi_{0}^{\prime}}(\zeta^{-1}) - \left[ \frac{\overline{z}(\zeta^{-1})}{z^{\prime}(\zeta)}\varphi_{0}^{\prime}(\zeta) \right]^{\prime}, \quad \zeta \in \overline{\bm{\omega}}_{0}
\end{equation}
Subsequently, Eqs. (\ref{eq:4.2b}) and (\ref{eq:4.2c}) can be equivalently rewritten \cite[]{Muskhelishvili1966,LIN2024396,LIN2024106008} as
\begin{equation}
  \label{eq:4.2b'}
  \tag{4.2b'}
  \begin{aligned}
    \sigma_{\rho0}(\zeta) + {\rm i} \tau_{\rho\theta0}(\zeta) = & \left[ \overline{\zeta}^{2} \frac{z\left( \overline{\zeta}^{-1} \right)} {z^{\prime}\left( \overline{\zeta}^{-1} \right)} - \frac{\overline{\zeta}}{\zeta} \frac{z(\zeta)} {z^{\prime}(\zeta)} \right] \overline{\varPhi_{0}^{\prime}(\zeta)} + \left[ \overline{\zeta}^{2} \frac{\overline{z^{\prime}(\zeta)}} {z^{\prime}\left( \overline{\zeta}^{-1} \right) } - \frac{\overline{\zeta}}{\zeta} \frac{\overline{z^{\prime}(\zeta)}} {z^{\prime}(\zeta)} \right] \overline{\varPsi_{0}(\zeta)} \\
    & + \varPhi_{0}(\zeta) - \varPhi_{0} \left( \overline{\zeta}^{-1} \right), \quad \zeta = \rho \cdot \sigma \in {\bm{\omega}}_{0}^{+}
  \end{aligned}
\end{equation}
\begin{equation}
  \label{eq:4.2c'}
  \tag{4.2c'}
  \begin{aligned}
    \frac{{\rm d}g_{0}(\zeta)}{{\rm d}\zeta} = & \left[ \frac{\overline{\zeta}}{\zeta} z(\zeta) - \overline{\zeta}^{2} \frac{z^{\prime}(\zeta)} {z^{\prime}\left( \overline{\zeta}^{-1} \right)} z\left( \overline{\zeta}^{-1} \right) \right] \overline{\varPhi_{0}^{\prime}(\zeta)} + \overline{z^{\prime}(\zeta)} \left[ \frac{\overline{\zeta}}{\zeta} - \overline{\zeta}^{2} \frac{z^{\prime}(\zeta)} {z^{\prime}\left( \overline{\zeta}^{-1} \right)} \right] \overline{\varPsi_{0}(\zeta)} \\
                                                                                                              & + z^{\prime}(\zeta) \left[\kappa  \varPhi_{0}(\zeta) + \varPhi_{0}\left(\overline{\zeta}^{-1}\right) \right], \quad \zeta = \rho \cdot \sigma \in {\bm \omega}_{0}^{+}
  \end{aligned}
\end{equation}
where $ \sigma = {\rm{e}}^{{\rm{i}}\theta} $ for simplicity in expressions.

When $ \rho \rightarrow 1 $, Eqs. (\ref{eq:4.2b'}) and (\ref{eq:4.2c'}) can be simplified by substitution of Eq. (\ref{eq:4.3}) as \cite[]{LIN2024396,LIN2024106008}
\begin{subequations}
  \label{eq:4.5}
  \begin{equation}
    \label{eq:4.5a}
    z^{\prime}(\zeta) \left[ \sigma_{\rho0}(\zeta) + {\rm{i}} \tau_{\rho\theta0}(\zeta) \right]|_{\rho \rightarrow 1} = \varphi_{0}^{\prime+}(\sigma) - \varphi_{0}^{\prime-}(\sigma) = 0, \quad \sigma \in {\bm{c}}_{1f0}
  \end{equation}
  \begin{equation}
    \label{eq:4.5b}
    \left. \frac{{\rm{d}}g_{0}(\zeta)}{{\rm{d}}\zeta} \right.|_{\rho \rightarrow 1} = \kappa \varphi_{0}^{\prime+}(\sigma) + \varphi_{0}^{\prime-}(\sigma) = 0, \quad \sigma \in {\bm{c}}_{1c0}
  \end{equation}
\end{subequations}
where $ \varphi_{0}^{\prime+}(\sigma) $ and $ \varphi_{0}^{\prime-}(\sigma) $ denote the values of $ \varphi_{0}^{\prime}(\zeta) = z^{\prime}(\zeta)\varPhi_{0}(\zeta) $ approaching boundary $ {\bm{c}}_{10} $ from $ {\bm{\omega}}_{0}^{+} $ and $ {\bm{\omega}}_{0}^{-} $ sides, respectively. Eq. (\ref{eq:4.5}) forms a homogenerous Riemann-Hilbert problem with extra constraint in Eq. (\ref{eq:4.1c}]) \cite[]{LIN2024396,LIN2024106008}. Similar to Eqs. (\ref{eq:4.1a'}) and (\ref{eq:4.1b'}), the complementary expressions of Eqs. (\ref{eq:4.5a}) and (\ref{eq:4.5b}) can be respectively given as
\begin{equation}
  \label{eq:4.5a'}
  \tag{4.5a'}
  z^{\prime}(\zeta) \left[ \sigma_{\rho0}(\zeta) + {\rm{i}} \tau_{\rho\theta0}(\zeta) \right]|_{\rho \rightarrow 1} = \varphi_{0}^{\prime+}(\sigma) - \varphi_{0}^{\prime-}(\sigma) \neq 0, \quad \sigma \in {\bm{c}}_{1c0}
\end{equation}
\begin{equation}
  \label{eq:4.5b'}
  \tag{4.5b'}
  \left. \frac{{\rm{d}}g_{0}(\zeta)}{{\rm{d}}\zeta} \right.|_{\rho \rightarrow 1} = \kappa \varphi_{0}^{\prime+}(\sigma) + \varphi_{0}^{\prime-}(\sigma) \neq 0, \quad \sigma \in {\bm{c}}_{1f0}
\end{equation}

\section{Problem solution}
\label{sec:problem-solution}

\subsection{Preparation of complex potentials}
\label{sec:complex-potentials}

Via Plemelj formula \cite[]{Muskhelishvili1966}, the solution of Eq. (\ref{eq:4.5}) can be expressed as
\begin{equation}
  \label{eq:5.1}
  \varphi_{0}^{\prime}(\zeta) = X(\zeta) \sum\limits_{n=-\infty}^{\infty} d_{n}\zeta^{n}, \quad \zeta \in \overline{\bm{\omega}}_{0} \cup {\bm{\omega}}_{0}^{-}
\end{equation}
where
\begin{equation}
  \label{eq:5.2}
  X(\zeta) = (\zeta-t_{10})^{-\frac{1}{2}-{\rm{i}}\lambda} (\zeta-t_{20})^{-\frac{1}{2}+{\rm{i}}\lambda}, \quad \lambda = \frac{\ln\kappa}{2\pi}, t_{10} = {\rm{e}}^{{\rm{i}}\theta_{1}}, t_{20} = {\rm{e}}^{{\rm{i}}\theta_{2}}
\end{equation}
$ d_{n} $ denotes the complex coefficients to be determined. Eq. (\ref{eq:5.2}) can be respectively expanded into rational series near the origin and infinity in the mapping plane $ \zeta $ as
\begin{subequations}
  \label{eq:5.3}
  \begin{equation}
    \label{eq:5.3a}
    X(\zeta) = \sum\limits_{k = 0}^{\infty} \alpha_{k}\zeta^{k}, \quad \zeta \in {\bm{\omega}}_{0}^{+}
  \end{equation}
  \begin{equation}
    \label{eq:5.3b}
    X(\zeta) = \sum\limits_{k = 0}^{\infty} \beta_{k}\zeta^{-k}, \quad \zeta \in {\bm{\omega}}_{0}^{-}
  \end{equation}
\end{subequations}
where
\begin{subequations}
  \label{eq:5.4}
  \begin{equation}
    \label{eq:5.4a}
    \left\{
      \begin{aligned}
        \alpha_{0} = & \; - t_{10}^{-\frac{1}{2}-{\rm{i}}\lambda} t_{20}^{-\frac{1}{2}+{\rm{i}}\lambda} \\
        \alpha_{k} = & \; - t_{10}^{-\frac{1}{2}-{\rm{i}}\lambda} t_{20}^{-\frac{1}{2}+{\rm{i}}\lambda} \cdot (-1)^{k} \left( c_{k}t_{10}^{-k} + \overline{c}_{k}t_{20}^{-k} + \sum\limits_{l=1}^{k-1} c_{l}\overline{c}_{k-l}\cdot t_{10}^{-l}t_{20}^{-k+l} \right), \quad k \geq 1 \\
      \end{aligned}
    \right.
  \end{equation}
  \begin{equation}
    \label{eq:5.4b}
    \left\{
      \begin{aligned}
        \beta_{0} = & \; 0 \\
        \beta_{1} = & \; 1 \\
        \beta_{k} = & \; (-1)^{k-1} \left( c_{k-1}t_{10}^{k-1} + \overline{c}_{k-1}t_{20}^{k-1} + \sum\limits_{l=1}^{k-2} c_{l}\overline{c}_{k-l-1}\cdot t_{10}^{l}t_{20}^{k-l-1} \right), \quad k \geq 2\\
      \end{aligned}
    \right.
  \end{equation}
  with
  \begin{equation}
    \label{eq:5.4c}
    c_{k} = \prod\limits_{l=1}^{k} \left( \frac{1}{2}-{\rm{i}}\lambda-l \right)/{k!}
  \end{equation}
\end{subequations}
The branch of Eq. (\ref{eq:5.3}) can be easily verfied by $ \lim\limits_{\zeta\rightarrow\infty} \zeta X(\zeta) = 1 $.

Substituting Eq. (\ref{eq:5.3}) into Eq. (\ref{eq:5.1}) yields
\begin{subequations}
  \label{eq:5.5}
  \begin{equation}
    \label{eq:5.5a}
    \varphi_{0}^{\prime}(\zeta) = \sum\limits_{k=-\infty}^{\infty} A_{k} \zeta^{k}, \quad A_{k} = \sum\limits_{n=-\infty}^{k} \alpha_{k-n}d_{n}, \quad \zeta \in {\bm{\omega}}_{0}^{+}
  \end{equation}
  \begin{equation}
    \label{eq:5.5b}
    \varphi_{0}^{\prime}(\zeta) = \sum\limits_{k=-\infty}^{\infty} B_{k} \zeta^{k}, \quad B_{k} = \sum\limits_{n=k}^{\infty} \beta_{-k+n}d_{n}, \quad \zeta \in {\bm{\omega}}_{0}^{-}
  \end{equation}
\end{subequations}
Substituting Eq. (\ref{eq:5.5a}) into Eq. (\ref{eq:4.4}) yields
\begin{equation}
  \label{eq:5.6}
  \psi_{0}^{\prime}(\zeta) = \sum\limits_{k=-\infty}^{\infty} \overline{B}_{-k-2}\zeta^{k} - \left[ \frac{\overline{z}(\zeta^{-1})}{z^{\prime}(\zeta)}\varphi_{0}^{\prime}(\zeta) \right]^{\prime}, \quad \zeta \in {\bm{\omega}}_{0}^{+}
\end{equation}
Eqs. (\ref{eq:5.5}) and (\ref{eq:5.6}) complete the preparation of complex potentials.

\subsection{Static equilibrium and displacement single-valuedness}
\label{sec:result-equil-displ}

The simulating static equilibrium in Eq. (\ref{eq:2.8'}) should be considered. Substituting Eq. (\ref{eq:4.1b'}) into the left-hand side of Eq. (\ref{eq:2.8'}) with consideration of Eq. (\ref{eq:4.1b}) and backward conformal mapping in Eq. (\ref{eq:3.3b}) yields
\begin{equation}
  \label{eq:5.7}
  \int_{{\bm{C}}_{1c0}} [X_{0}(T_{0})+{\rm{i}}Y_{0}(T_{0})]|{\rm{d}}T_{0}| = \int_{{\bm{C}}_{10}} [X_{0}(T_{0})+{\rm{i}}Y_{0}(T_{0})]|{\rm{d}}T_{0}| = \ointctrclockwise_{{\bm{c}}_{10}} t_{0} \cdot z^{\prime}(t_{0}) [\sigma_{\rho0}(t_{0}) + {\rm{i}}\tau_{\rho\theta0}(t_{0})]{\rm{d}}\theta, \quad t_{0} = {\rm{e}}^{{\rm{i}}\theta}
\end{equation}
where $ |{\rm{d}}T_{0}| = |z^{\prime}(t_{0})| \cdot |{\rm{i}}t_{0}|{\rm{d}}\theta = |z^{\prime}(t_{0})|{\rm{d}}\theta $. Substituting Eqs. (\ref{eq:4.5a'}) and (\ref{eq:5.5}) into the last integration of Eq. (\ref{eq:5.7}) yields
\begin{equation}
  \label{eq:5.7'}
  \tag{5.7'}
  -{\rm{i}}\ointctrclockwise_{{\bm{c}}_{10}} \sum\limits_{k=-\infty}^{\infty} (A_{k}-B_{k})t_{0}^{k} {\rm{d}}t_{0} = \sum\limits_{k=-\infty}^{\infty} (A_{k}-B_{k}) \int_{0}^{2\pi} {\rm{e}}^{{\rm{i}}k\theta} {\rm{d}}\theta = 2\pi(A_{-1}-B_{-1})
\end{equation}
Replacing the left-hand side of Eq. (\ref{eq:2.8'}) with Eq. (\ref{eq:5.7'}) yields
\begin{equation}
  \label{eq:5.8}
  A_{-1} - B_{-1} = \frac{R_{y}}{2\pi} = -{\rm{i}}W\gamma
\end{equation}

The displacement single-valuedness in the simulating geomaterial $ \overline{\bm{\varOmega}}_{0} $ should also be considered. Substituting Eqs. (\ref{eq:5.5a}), (\ref{eq:5.6}), and (\ref{eq:4.3}) into Eq. (\ref{eq:4.2c}) yields
\begin{equation}
  \label{eq:5.9}
  \begin{aligned}
    g_{0}(\zeta)
    = & \; \kappa \int \sum\limits_{k=-\infty}^{\infty} A_{k}\zeta^{k} {\rm{d}}\zeta - \int \sum\limits_{k=-\infty}^{\infty} B_{-k-2} \overline{\zeta}^{k} {\rm{d}}\overline{\zeta} + \frac{z(\overline{\zeta}^{-1}) - z(\zeta)} {\overline{z^{\prime}[w(\zeta)]} \cdot \overline{w^{\prime}(\zeta)}} {\overline{\varphi_{0}^{\prime}(\zeta)}}  \\
    = & \; \sum\limits_{\substack{k=-\infty \\ k \neq -1}}^{\infty} \left( \kappa A_{k} \frac{\zeta^{k+1}}{k+1} - B_{-k-2} \frac{\overline{\zeta}^{k+1}}{k+1} \right) + \frac{z(\overline{\zeta}^{-1}) - z(\zeta)} {\overline{z^{\prime}[w(\zeta)]} \cdot \overline{w^{\prime}(\zeta)}} {\overline{\varphi_{0}^{\prime}(\zeta)}} \\
      & \; + (\kappa A_{-1} - B_{-1})\ln\rho + C_{0} + (\kappa A_{-1}+B_{-1}){\rm{Ln}}\sigma \\
  \end{aligned}
\end{equation}
where $ {\rm{Ln}} $ denotes the multi-valued logarithm, $ C_{0} $ denotes the complex constant of rigid-body displacement translation, which should be determined in Eq. (\ref{eq:6.5}). To guarantee displacement single-valuedness, the last item of Eq. (\ref{eq:5.9}) should spontaneously vanish, thus, we have
\begin{equation}
  \label{eq:5.10}
  \kappa A_{-1} + B_{-1} = 0
\end{equation}
Eqs. (\ref{eq:5.8}) and (\ref{eq:5.10}) give
\begin{equation}
  \label{eq:5.11}
  \left\{
    \begin{aligned}
      & A_{-1} = \frac{-{\rm{i}}W\gamma}{1+\kappa} \\
      & B_{-1} = \frac{{\rm{i}}\kappa W\gamma}{1+\kappa} \\
    \end{aligned}
  \right.
\end{equation}
Eqs. (\ref{eq:5.11}) indicates that the values of $ A_{-1} $ and $ B_{-1} $ are only related the traction resultant along the inner boundary and the Kolosov parameter, and are free from conformal mapping.

\subsection{Traction along inner boundary}
\label{sec:traction-along-inner}

The traction equilibrium along inner boundary $ {\bm{c}}_{20} $ in Eq. (\ref{eq:4.1c}) can be equivalently modified to its length integration formation \cite[]{LIN2024396,LIN2024106008} as
\begin{equation}
  \label{eq:5.12}
  -{\rm{i}}\int \frac{s_{0}}{|s_{0}|} \frac{z^{\prime}(s_{0})}{|z^{\prime}(s_{0})|} [\sigma_{\rho0}(s_{0}) + {\rm{i}}\tau_{\rho\theta0}(s_{0})] |{\rm{d}}S_{0}| = -{\rm{i}}\int \left[ -k_{0}\gamma y(S_{0})\frac{{\rm{d}}y(S_{0})}{{\rm{d}}S_{0}} + {\rm{i}}\gamma y(S_{0})\frac{{\rm{d}}x(S_{0})}{{\rm{d}}S_{0}} \right] |{\rm{d}}S_{0}|
\end{equation}
Considering $ |{\rm{d}}S_{0}| = -{\rm{d}}S_{0} $ and $ |{\rm{d}}S_{0}| = |z^{\prime}(s_{0})|\cdot|{\rm{d}}s_{0}| = |z^{\prime}(s_{0})|\cdot\alpha|{\rm{d}}\theta| = -|z^{\prime}(s_{0})|\cdot\alpha{\rm{d}}\theta $, Eq. (\ref{eq:5.12}) can be simplified as
\begin{equation}
  \label{eq:5.12'}
  \tag{5.12'}
  \int z^{\prime}(s_{0})[\sigma_{\rho0}(s_{0})+{\rm{i}}\tau_{\rho\theta0}(s_{0})]{\rm{d}}s_{0} = -\gamma\int {y(S_{0})}\left( k_{0}\frac{{\rm{i}}{\rm{d}}y(S_{0})}{{\rm{d}}\sigma} + \frac{{\rm{d}}x(S_{0})}{{\rm{d}}\sigma} \right) {\rm{d}}\sigma
\end{equation}

Substituting Eqs. (\ref{eq:4.2b}), (\ref{eq:5.5a}), and (\ref{eq:5.6}) into the left-hand side of Eq. (\ref{eq:5.12'}) yields
\begin{equation}
  \label{eq:5.13}
  \begin{aligned}
    & \int z^{\prime}(s_{0})[\sigma_{\rho0}(s_{0})+{\rm i}\tau_{\rho\theta0}(s_{0})] {\rm{d}}s_{0} \\
    = & \; \int \varphi_{0}^{\prime}(s_{0}){\rm{d}}s_{0} + \frac{z(s_{0})}{\overline{z^{\prime}(s_{0})}}\overline{\varphi_{0}^{\prime}(s_{0})} + \int \overline{\psi_{0}^{\prime}(s_{0})}{\rm{d}}\overline{s}_{0} \\
    = & \; \sum\limits_{\substack{k=-\infty \\ k \neq 0}}^{\infty} A_{k-1} \frac{\alpha^{k}}{k}\sigma^{k} + \sum\limits_{\substack{k=-\infty \\ k \neq 0}}^{\infty} B_{-k-1} \frac{\alpha^{k}}{k}\sigma^{-k} + \frac{z[w(\alpha\sigma)]-z[w(\alpha^{-1}\sigma)]}{\overline{z^{\prime}[w(\alpha\sigma)]} \cdot \overline{w^{\prime}(\alpha\sigma)}} \sum\limits_{k=-\infty}^{\infty} \overline{A}_{k} \alpha^{k}\sigma^{-k} \\
    & \; + (A_{-1}+B_{-1})\ln\alpha + C_{a} + (A_{-1}-B_{-1}){\rm{Ln}}\sigma
  \end{aligned}
\end{equation}
where $ C_{a} $ denotes the integral constant, and the following expansion can be conducted:
\begin{equation}
  \label{eq:5.14}
   \sum\limits_{k=-\infty}^{\infty} f_{k}(\alpha) \sigma^{k} = \frac {z[w(\alpha\sigma)]-z[w(\alpha^{-1}\sigma)]} {\overline{z^{\prime}[w(\alpha\sigma)]} \cdot \overline{w^{\prime}(\alpha\sigma)}}
\end{equation}
where $ f_{k}(\alpha) $ denotes the simulating coefficients when $ \rho = \alpha $. The expanding method of Eq. (\ref{eq:5.14}) is a numerical one with truncation and sample points, and the details can be found in Appendix B of Ref \cite[]{LIN2024106008}, and would not be repeated here. Substituting the right-hand side of Eq. (\ref{eq:5.14}) into Eq. (\ref{eq:5.13}) yields
\begin{equation}
  \label{eq:5.13'}
  \tag{5.13'}
  \begin{aligned}
    \int z^{\prime}(s_{0})[\sigma_{\rho0}(s_{0})+{\rm i}\tau_{\rho\theta0}(s_{0})] {\rm d}{s_{0}} 
    = \enspace & \; \sum\limits_{\substack{k=-\infty \\ k \neq 0}}^{\infty} \left( A_{k-1} \frac{\alpha^{k}}{k} - B_{k-1} \frac{\alpha^{-k}}{k} + \sum\limits_{l=-\infty}^{\infty} f_{l}(\alpha) \overline{A}_{l-k} \alpha^{l-k} \right) \sigma^{k} \\
    & \; + \sum\limits_{l=-\infty}^{\infty} f_{l}(\alpha) \overline{A}_{l} \alpha^{l} + (A_{-1}+B_{-1})\ln\alpha + C_{a} + (A_{-1}-B_{-1}){\rm{Ln}}\sigma
  \end{aligned}
\end{equation}

According to the two-step composite conformal mapping in Eq. (\ref{eq:3.3b}), the integrand elements of the right-hand side in Eq. (\ref{eq:5.12'}) can be expanded as
\begin{subequations}
  \label{eq:5.15} 
  \begin{equation}
    \label{eq:5.15a}
    x(S_{0}) = \frac{1}{2} \left\{ z[w(\alpha\sigma)] + \overline{z[w(\alpha\sigma)]} \right\}
  \end{equation}
  \begin{equation}
    \label{eq:5.15b}
    y(S_{0}) = -\frac{{\rm{i}}}{2} \left\{ z[w(\alpha\sigma)] - \overline{z[w(\alpha\sigma)]} \right\}
  \end{equation}
  \begin{equation}
    \label{eq:5.15c}
    \frac{{\rm{d}}x(S_{0})}{{\rm{d}}\sigma} = \frac{\alpha}{2} \left\{ z^{\prime}[w(\alpha\sigma)] \cdot w^{\prime}(\alpha\sigma) - \overline{z^{\prime}[w(\alpha\sigma)]} \cdot \overline{w^{\prime}(\alpha\sigma)} \cdot \sigma^{-2} \right\}
  \end{equation}
  \begin{equation}
    \label{eq:5.15d}
    \frac{{\rm{i}}{\rm{d}}y(S_{0})}{{\rm{d}}\sigma} = \frac{\alpha}{2} \left\{ z^{\prime}[w(\alpha\sigma)] \cdot w^{\prime}(\alpha\sigma) + \overline{z^{\prime}[w(\alpha\sigma)]} \cdot \overline{w^{\prime}(\alpha\sigma)} \cdot \sigma^{-2} \right\}
  \end{equation}
\end{subequations}
Similar to Eq. (\ref{eq:5.14}), the components the right-hand side in Eq.(\ref{eq:5.12'}) can be expanded according to Eq. (\ref{eq:5.15}) as
\begin{equation}
  \label{eq:5.16}
  \sum\limits_{k=-\infty}^{\infty} g_{k} \sigma^{k} = {y(S_{0})}\left( k_{0}\frac{{\rm{i}}{\rm{d}}y(S_{0})}{{\rm{d}}\sigma} + \frac{{\rm{d}}x(S_{0})}{{\rm{d}}\sigma} \right)
\end{equation}
Eq. (\ref{eq:5.16}) applies the same method of Eq. (\ref{eq:5.14}). Substituting the left-hand side of Eq. (\ref{eq:5.16}) into the right-hand side of Eq. (\ref{eq:5.12'}) yields
\begin{equation}
  \label{eq:5.17}
  -\gamma\int {y(S_{0})}\left( k_{0}\frac{{\rm{i}}{\rm{d}}y(S_{0})}{{\rm{d}}\sigma} + \frac{{\rm{d}}x(S_{0})}{{\rm{d}}\sigma} \right) {\rm{d}}\sigma = - \gamma \int \sum\limits_{k=-\infty}^{\infty} g_{k}\sigma^{k} {\rm{d}}\sigma = \sum\limits_{k=1}^{\infty} I_{k} \sigma^{-k} + \sum\limits_{k=1}^{\infty} J_{k} \sigma^{k} + K_{0}{\rm{Ln}}\sigma
\end{equation}
where
\begin{equation}
  \label{eq:5.18}
  \left\{
    \begin{aligned}
      & I_{k} = - \gamma \frac{g_{-k-1}}{-k} \\
      & J_{k} = - \gamma \frac{g_{k-1}}{k} \\
      & K_{0} = - \gamma g_{-1}
    \end{aligned}
  \right.
\end{equation}

\subsection{Iterative solution}
\label{sec:approximate-solution}

According to Eq. (\ref{eq:5.12'}), Eqs. (\ref{eq:5.13'}) and (\ref{eq:5.17}) should be equal to each other for arbitrary polar angle $ \theta $. Comparing the coefficients between Eqs. (\ref{eq:5.13'}) and (\ref{eq:5.17}) gives
\begin{subequations}
  \label{eq:5.19}
  \begin{equation}
    \label{eq:5.19a}
    A_{-k-1} \frac{\alpha^{-k}}{-k} - B_{-k-1} \frac{\alpha^{k}}{-k} + \sum\limits_{l=-\infty}^{\infty} f_{l}(\alpha) \overline{A}_{l+k} \alpha^{l+k} = I_{k}, \quad k \geq 1
  \end{equation}
  \begin{equation}
    \label{eq:5.19b}
    A_{k-1} \frac{\alpha^{k}}{k} - B_{k-1} \frac{\alpha^{-k}}{k} + \sum\limits_{l=-\infty}^{\infty} f_{l}(\alpha) \overline{A}_{l-k} \alpha^{l-k} = J_{k}, \quad k \geq 1
  \end{equation}
  \begin{equation}
    \label{eq:5.19c}
    \sum\limits_{l=-\infty}^{\infty} f_{l}(\alpha) \overline{A}_{l} \alpha^{l} + (A_{-1}+B_{-1})\ln\alpha + C_{a} = 0
  \end{equation}
  \begin{equation}
    \label{eq:5.19d}
    A_{-1} - B_{-1} = K_{0}
  \end{equation}
\end{subequations}
Eqs. (\ref{eq:5.19a}) and (\ref{eq:5.19b}) can be rewritten as
\begin{equation}
  \label{eq:5.19a'}
  \tag{5.19a'}
  A_{-k-1} = - k \alpha^{k} I_{k} + \alpha^{2k}B_{-k-1} + k \alpha^{2k} \sum\limits_{l=-\infty}^{\infty} f_{l}(\alpha) \overline{A}_{l+k} \alpha^{l}, \quad k \geq 1
\end{equation}
\begin{equation}
  \label{eq:5.19b'}
  \tag{5.19b'}
  B_{k-1} = - k \alpha^{k} J_{k} + \alpha^{2k}A_{k-1} + k \sum\limits_{l=-\infty}^{\infty} f_{l}(\alpha) \overline{A}_{l-k} \alpha^{l}, \quad k \geq 1
\end{equation}
With Eq. (\ref{eq:5.11}), Eqs. (\ref{eq:5.19c}) and (\ref{eq:5.19d}) can be rewritten as
\begin{equation}
  \label{eq:5.19c'}
  \tag{5.19c'}
  C_{a} = - \sum\limits_{l=-\infty}^{\infty} f_{l}(\alpha) \overline{A}_{l} \alpha^{l} + {\rm{i}}\frac{\kappa-1}{\kappa+1}W\gamma \ln\alpha 
\end{equation}
\begin{equation}
  \label{eq:5.19d'}
  \tag{5.19d'}
  K_{0} = - {\rm{i}}W \gamma
\end{equation}

Substituting Eqs. (\ref{eq:5.5b}) and (\ref{eq:5.5a}) into the left-hand sides of Eqs. (\ref{eq:5.19a'}) and (\ref{eq:5.19b'}) yields
\begin{subequations}
  \label{eq:5.20} 
  \begin{equation}
    \label{eq:5.20a}
    \sum\limits_{n=k+1}^{\infty} \alpha_{-k-1+n} d_{-n} = - k \alpha^{k} I_{k} + \alpha^{2k}B_{-k-1} + k \alpha^{2k} \sum\limits_{l=-\infty}^{\infty} \alpha^{l} f_{l}(\alpha) \overline{A}_{l+k}, \quad k \geq 1
  \end{equation}
  \begin{equation}
    \label{eq:5.20b}
    \sum\limits_{n=k}^{\infty} \beta_{-k+1+n} d_{n} = - k \alpha^{k} J_{k} + \alpha^{2k}A_{k-1} + k \sum\limits_{l=-\infty}^{\infty} \alpha^{l} f_{l}(\alpha) \overline{A}_{l-k}, \quad k \geq 1
  \end{equation}
\end{subequations}
Eq. (\ref{eq:5.20}) constrains $ d_{n} (n \geq 1) $ and $ d_{-n} (n \geq 2) $, while $ d_{0} $ and $ d_{-1} $ are not constrained. To obtain $ d_{0} $ and $ d_{-1} $, Eq. (\ref{eq:5.11}) can be rewritten with substitution of Eqs. (\ref{eq:5.5}) and (\ref{eq:5.19d'}) as
\begin{subequations}
  \label{eq:5.21}
  \begin{equation}
    \label{eq:5.21a}
    \sum\limits_{n=1}^{\infty} \alpha_{-1+n}d_{-n} = \frac{K_{0}}{1+\kappa}
  \end{equation}
  \begin{equation}
    \label{eq:5.21b}
    \sum\limits_{n=0}^{\infty} \beta_{1+n}d_{n} = \frac{-\kappa K_{0}}{1+\kappa}
  \end{equation}
\end{subequations}
Therefore, Eqs. (\ref{eq:5.20}) and (\ref{eq:5.21}) constrain $ d_{n} (-\infty < n < \infty) $. Eq. (\ref{eq:5.19c'}) is generally trivial, since it just determines an integral constant, which is not necessary to solve. Eq. (\ref{eq:5.19d'}) should be always spontaneously satisfied for a given cavity shape.

To obtain $ d_{n} (-\infty < n < \infty) $, an iterativie method is used. We assume that $ d_{n} (-\infty < n < \infty) $ can be expanded into iteration formation as
\begin{equation}
  \label{eq:5.22}
  d_{n} = \sum\limits_{q=0}^{\infty} d_{n}^{(q)}, \quad -\infty < n < \infty
\end{equation}
where $ q $ denotes iteration reps. With Eq. (\ref{eq:5.22}), Eqs. (\ref{eq:5.20}) and (\ref{eq:5.21}) can be reorganized and rewritten as
\begin{subequations}
  \label{eq:5.23}
  \begin{equation}
    \label{eq:5.23a}
    \left\{
      \begin{aligned}
        & \sum\limits_{n=1}^{\infty} \alpha_{-1+n}\sum\limits_{q=0}^{\infty} d_{-n}^{(q)} = \frac{K_{0}}{1+\kappa} \\
        & \sum\limits_{n=k+1}^{\infty} \alpha_{-k-1+n} \sum\limits_{q=0}^{\infty} d_{-n}^{(q)} = - k \alpha^{k} I_{k} + \alpha^{2k}\sum\limits_{q=0}^{\infty} B_{-k-1}^{(q)} + k \alpha^{2k} \sum\limits_{l=-\infty}^{\infty} \alpha^{l} f_{l}(\alpha) \sum\limits_{q=0}^{\infty} \overline{A}_{l+k}^{(q)}, \quad k \geq 1
      \end{aligned}
    \right.
  \end{equation}
  \begin{equation}
    \label{eq:5.23b}
    \left\{
      \begin{aligned}
        & \sum\limits_{n=0}^{\infty} \beta_{1+n} \sum\limits_{q=0}^{\infty} d_{n}^{(q)} = \frac{-\kappa K_{0}}{1+\kappa} \\
        & \sum\limits_{n=k}^{\infty} \beta_{-k+1+n} \sum\limits_{q=0}^{\infty} d_{n}^{(q)} = - k \alpha^{k} J_{k} + \alpha^{2k} \sum\limits_{q=0}^{\infty} A_{k-1}^{(q)} + k \sum\limits_{l=-\infty}^{\infty} \alpha^{l} f_{l}(\alpha) \sum\limits_{q=0}^{\infty} \overline{A}_{l-k}^{(q)}, \quad k \geq 1
      \end{aligned}
    \right.
  \end{equation}
\end{subequations}
where $ A_{k} $ and $ B_{k} $ are also expanded into iteration formation as
\begin{equation*}
  \left\{
    \begin{aligned}
      & A_{k} = \sum\limits_{q=0}^{\infty} A_{k}^{(q)} \\
      & B_{k} = \sum\limits_{q=0}^{\infty} B_{k}^{(q)} \\
    \end{aligned}
  \right.
\end{equation*}

For iteration $ q = 0 $, certain items in Eq. (\ref{eq:5.23}) can be artificially extracted to obtain $ d_{n}^{(0)} (-\infty < n < \infty) $ as
\begin{subequations}
  \label{eq:5.24}
  \begin{equation}
    \label{eq:5.24a}
    \left\{
      \begin{aligned}
        & \sum\limits_{n=1}^{\infty} \alpha_{-1+n}d_{-n}^{(0)} = \frac{K_{0}}{1+\kappa} \\
        & \sum\limits_{n=k+1}^{\infty} \alpha_{-k-1+n} d_{-n}^{(0)} = - k \alpha^{k} I_{k}, \quad k \geq 1
      \end{aligned}
    \right.
  \end{equation}
  \begin{equation}
    \label{eq:5.24b}
    \left\{
      \begin{aligned}
        & \sum\limits_{n=0}^{\infty} \beta_{1+n}d_{n}^{(0)} = \frac{-\kappa K_{0}}{1+\kappa} \\
        & \sum\limits_{n=k}^{\infty} \beta_{-k+1+n} d_{n}^{(0)} = - k \alpha^{k} J_{k}, \quad k \geq 1 
      \end{aligned}
    \right.
  \end{equation}
\end{subequations}
Eq. (\ref{eq:5.24}) forms two independent linear systems to start the iteration. Then for iteration $ q \geq 0 $, $ A_{k}^{(q)} $ and $ B_{k}^{(q)} $ in Eq. (\ref{eq:5.23}) can be obtained by Eq. (\ref{eq:5.5}) as
\begin{subequations}
  \label{eq:5.25}
  \begin{equation}
    \label{eq:5.25a}
    A_{k}^{(q)} = \sum\limits_{n=-k}^{\infty} \alpha_{k+n}d_{-n}^{(q)}
  \end{equation}
  \begin{equation}
    \label{eq:5.25b}
    B_{k}^{(q)} = \sum\limits_{n=k}^{\infty} \beta_{-k+n}d_{n}^{(q)}
  \end{equation}
\end{subequations}
For the next iteration $ q+1 $, corresponding items in Eq. (\ref{eq:5.23}) can be extracted to obtain $ d_{n}^{(q+1)} (-\infty < n < \infty) $ as
\begin{subequations}
  \label{eq:5.26}
  \begin{equation}
    \label{eq:5.26a}
    \left\{
      \begin{aligned}
        & \sum\limits_{n=1}^{\infty} \alpha_{-1+n}d_{-n}^{(q+1)} = 0 \\
        & \sum\limits_{n=k+1}^{\infty} \alpha_{-k-1+n} d_{-n}^{(q+1)} = \alpha^{2k}B_{-k-1}^{(q)} + k \alpha^{2k} \sum\limits_{l=-\infty}^{\infty} \alpha^{l} f_{l}(\alpha) \overline{A}_{l+k}^{(q)}, \quad k \geq 1
      \end{aligned}
    \right. 
  \end{equation}
  \begin{equation}
    \label{eq:5.26b}
    \left\{
      \begin{aligned}
        & \sum\limits_{n=0}^{\infty} \beta_{1+n}d_{n}^{(q+1)} = 0 \\
        & \sum\limits_{n=k}^{\infty} \beta_{-k+1+n} d_{n}^{(q+1)} = \alpha^{2k}A_{k-1}^{(q)} + k \sum\limits_{l=-\infty}^{\infty} \alpha^{l} f_{l}(\alpha) \overline{A}_{l-k}^{(q)}, \quad k \geq 1
      \end{aligned}
    \right.
  \end{equation}
\end{subequations}
Eq. (\ref{eq:5.26}) also forms two independent linear systems to proceed the iteration. The iteration procedure in Eqs. (\ref{eq:5.24})-(\ref{eq:5.26}) may stop, when the following threshold is reached:
\begin{equation}
  \label{eq:5.27}
  \max|d_{n}^{(Q+1)}| \leq \varepsilon, \quad -\infty < n < \infty
\end{equation}
where $ Q $ denotes the iteration rep to satisfy Eq. (\ref{eq:5.27}), $ \varepsilon $ is a small real constant, and $ \varepsilon = 10^{-16} $ could be used. After the iteration, the results of $ d_{n} $ could be accumulated via Eq. (\ref{eq:5.22}) as
\begin{equation}
  \label{eq:5.22'}
  \tag{5.22'}
  d_{n} = \sum\limits_{q=0}^{Q} d_{n}^{(q)}, \quad -\infty < n < \infty
\end{equation}

\section{Stress and displacement in simulating geomaterial}
\label{sec:stress-displ-simul}

Eq. (\ref{eq:5.13'}) is not only suitable along the inner periphery of the simulating unit annulus $ {\bm{\omega}}_{0} $, but also suitable for arbitrary circle within the simulating unit annulus $ {\bm{\omega}}_{0} $ for a given polar radius $ \rho $ as
\begin{equation}
  \label{eq:6.1}
  \begin{aligned}
    \int z^{\prime}(\rho\sigma)[\sigma_{\rho0}(\rho\sigma)+{\rm i}\tau_{\rho\theta0}(\rho\sigma)] {\rm d}({\rho\sigma}) 
    = \enspace & \; \sum\limits_{\substack{k=-\infty \\ k \neq 0}}^{\infty} \left( A_{k-1} \frac{\rho^{k}}{k} - B_{k-1} \frac{\rho^{-k}}{k} + \sum\limits_{l=-\infty}^{\infty} f_{l}(\rho) \overline{A}_{l-k} \rho^{l-k} \right) \sigma^{k} \\
    & \; + \sum\limits_{l=-\infty}^{\infty} f_{l}(\rho) \overline{A}_{l} \rho^{l} + (A_{-1}+B_{-1})\ln\rho + C_{a} + (A_{-1}-B_{-1}){\rm{Ln}}\sigma
  \end{aligned}
\end{equation}
where the following simulation is used in the same manner of Eq. (\ref{eq:5.14}) as
\begin{equation}
  \label{eq:6.2}
  \sum\limits_{k=-\infty}^{\infty} f_{k}(\rho) \sigma^{k} = \frac {z[w(\rho\sigma)]-z[w(\rho^{-1}\sigma)]} {\overline{z^{\prime}[w(\rho\sigma)]}\cdot\overline{w^{\prime}(\rho\sigma)}}
\end{equation}
With Eq. (\ref{eq:6.1}), the curvilinear radial and shear stress components mapped onto the simulating unit annulus $ {\bm{\omega}}_{0} $ can be obtained in a compact formation as
\begin{subequations}
  \label{eq:6.3}
  \begin{equation}
    \label{eq:6.3a}
    \begin{aligned}
      \sigma_{\rho0}(\rho\sigma)+{\rm i}\tau_{\rho\theta0}(\rho\sigma) = \frac{1}{z^{\prime}[w(\rho\sigma)]\cdot w^{\prime}(\rho\sigma)} \sum\limits_{k=-\infty}^{\infty} \left[ A_{k}\rho^{k} - B_{k}\rho^{-k-2} + (k+1) \sum\limits_{l=-\infty}^{\infty} f_{l}(\rho)\overline{A}_{l-k-1}\rho^{l-k-2} \right] \sigma^{k}
    \end{aligned}
  \end{equation}
  Substituting Eq. (\ref{eq:5.5}) into Eq. (\ref{eq:4.2a}) with substitution of $ \zeta = \rho\sigma $ yields the curvilinear hoop and radial stress components mapped onto the simulating unit annulus $ {\bm{\omega}}_{0} $ as
  \begin{equation}
    \label{eq:6.3b}
    \sigma_{\theta0}(\rho\sigma) + \sigma_{\rho0}(\rho\sigma) = 4\Re \left[ \frac{1}{z^{\prime}[w(\rho\sigma)]\cdot w^{\prime}(\rho\sigma) } \sum\limits_{k=-\infty}^{\infty} A_{k}\rho^{k}\sigma^{k} \right]
  \end{equation}
  Substituting Eq. (\ref{eq:6.2}) into Eq. (\ref{eq:5.9}) with substitution of $ \zeta = \rho\sigma $ gives the rectangular displacement components mapped onto the simulating unit annulus $ {\bm{\omega}}_{0} $ as
  \begin{equation}
    \label{eq:6.3c}
    \begin{aligned}
      g_{0}(\rho\sigma) = 2G[u_{0}(\rho\sigma) + {\rm{i}}v_{0}(\rho\sigma)]
      = & \; \sum\limits_{\substack{k=-\infty \\ k \neq 0}}^{\infty} \left( \kappa A_{k-1} \frac{\rho^{k}}{k} + B_{k-1} \frac{\rho^{-k}}{k} \right) \sigma^{k} + C_{0} \\
        & \; - \sum\limits_{k=-\infty}^{\infty} \sum\limits_{l=-\infty}^{\infty} f_{l}(\rho) \overline{A}_{l-k} \rho^{l-k} \sigma^{k} + (\kappa A_{-1} - B_{-1})\ln{\rho}
    \end{aligned}
  \end{equation}
\end{subequations}

When $ \rho \rightarrow 1 $, Eq. (\ref{eq:6.3}) would result in the stress and displacement components along the outer periphery of the simulating unit annulus $ {\bm{\omega}}_{0} $ in mapping plane $ \zeta = \rho\cdot{\rm{e}}^{{\rm{i}}\theta} $, which is corresponding to the ground surface of the simulating geomaterial $ {\bm{\varOmega}}_{0} $ in physical plane $ z = x + {\rm{i}}y $, as
\begin{subequations}
  \label{eq:6.4}
  \begin{equation}
    \label{eq:6.4a}
    \begin{aligned}
      \sigma_{\rho0}(\sigma)+{\rm i}\tau_{\theta0}(\sigma) = \frac{1}{z^{\prime}[w(\sigma)]\cdot w^{\prime}(\sigma)} \sum\limits_{k=-\infty}^{\infty} \left( A_{k} - B_{k} \right) \sigma^{k}
    \end{aligned}
  \end{equation}
  \begin{equation}
    \label{eq:6.4b}
    \sigma_{\theta0}(\sigma) + \sigma_{\rho0}(\sigma) = 4\Re \left[ \frac{1}{z^{\prime}[w(\sigma)]\cdot w^{\prime}(\sigma) } \sum\limits_{k=-\infty}^{\infty} A_{k}\sigma^{k} \right]
  \end{equation}
  \begin{equation}
    \label{eq:6.4c}
      g_{0}(\sigma) = 2G[u_{0}(\sigma) + {\rm{i}}v_{0}(\sigma)] = \sum\limits_{\substack{k=-\infty \\ k \neq 0}}^{\infty} \left( \kappa A_{k-1} + B_{k-1} \right) \frac{\sigma^{k}}{k} + C_{0}
  \end{equation}
\end{subequations}
Since $ z[w(1)] \rightarrow \infty $ for the two-step composite backward conformal mapping in Eq. (\ref{eq:3.3b}), $ g_{0}(1) = 0 $ should be satisfied in Eq. (\ref{eq:6.4c}) to meet the boundary condition of the mechanical model in Fig. \ref{fig:3}a-1, and we have
\begin{equation}
  \label{eq:6.5}
  C_{0} = - \sum\limits_{\substack{k=-\infty \\ k \neq 0}}^{\infty} (\kappa A_{k-1} + B_{k-1}) \frac{1^{k}}{k}
\end{equation}
where the unit constant is deliberately kept to correctly apply the Lanczos filtering in Eq. (\ref{eq:7.1}) after series truncation. Eq. (\ref{eq:6.5}) determines the constant mentioned in Eq. (\ref{eq:5.9}).

With the curvilinear stress and displacement components caused by excavation in the simulating unit annulus $ {\bm{\omega}}_{0} $ in the mapping plane $ \zeta = \rho\cdot{\rm{e}}^{{\rm{i}}\theta} $ in Eq. (\ref{eq:6.3}), the rectangular stress and displacement components caused by excavation in the simulating geomaterial $ {\bm{\varOmega}}_{0} $ in the physical plane $ z = x + {\rm{i}}y $ can be obtained as
\begin{subequations}
  \label{eq:6.6} 
  \begin{equation}
    \label{eq:6.6a}
    \left\{
      \begin{aligned}
        & \sigma_{y0}(z) + \sigma_{x0}(z) = \sigma_{\theta0}(\zeta) + \sigma_{\rho0}(\zeta) \\
        & \sigma_{y0}(z) - \sigma_{x0}(z) + 2{\rm{i}}\tau_{xy0}(z) = \left[ \sigma_{\theta0}(\zeta) - \sigma_{\rho0}(\zeta) + 2{\rm{i}} \tau_{\rho\theta0}(\zeta) \right] \cdot \frac {\overline{\zeta}} {\zeta} \frac {\overline{z^{\prime}[w(\zeta)]}} {z^{\prime}[w(\zeta)]} \frac{\overline{w^{\prime}(\zeta)}}{w^{\prime}(\zeta)} \\
      \end{aligned}
    \right.
  \end{equation}
  \begin{equation}
    \label{eq:6.6b}
    u_{0}(z) + {\rm{i}}v_{0}(z) = u_{0}(\zeta) + {\rm{i}}v_{0}(\zeta)
  \end{equation}
\end{subequations}
where $ \sigma_{x0}^{\ast} $, $ \sigma_{y0}^{\ast} $, and $ \tau_{xy0}^{\ast} $ denote the horizontal, vertical, and shear stress components caused by excavation in the simulating geomaterial $ {\bm{\varOmega}}_{0} $, respectively. The final stress components considering both excavation and the initial stress field can be obtained by accumulation of Eqs. (\ref{eq:6.6a}) and (\ref{eq:2.1}) as
\begin{equation}
  \label{eq:6.7}
  \left\{
    \begin{aligned}
      & \sigma_{x}^{\ast}(z) = \sigma_{x0}(z) + \sigma_{x}^{0}(z) \\
      & \sigma_{y}^{\ast}(z) = \sigma_{y0}(z) + \sigma_{y}^{0}(z) \\
      & \tau_{xy}^{\ast}(z) = \tau_{xy0}(z) + \tau_{xy}^{0}(z) \\
    \end{aligned}
  \right.
\end{equation}
where $ \sigma_{x}^{\ast} $, $ \sigma_{y}^{\ast} $, and $ \tau_{xy}^{\ast} $ denote the final horizontal, vertical, and shear stress components in the simulating geomaterial $ {\bm{\varOmega}}_{0} $, respectively. Till now, the rectangular stress and displacement components in the simulating geomaterial in Fig. \ref{fig:3}a-1 are solved in Eqs. (\ref{eq:6.7}) and (\ref{eq:6.6b}). Alternatively, Eq. (\ref{eq:6.7}) can be transformed into polar formation in the mapping plane as
\begin{equation}
  \label{eq:6.7'}
  \tag{6.7'}
  \left\{
    \begin{aligned}
      & \sigma_{\theta}^{\ast}(z) + \sigma_{\rho}^{\ast}(z) = \sigma_{y}^{\ast}(z) + \sigma_{x}^{\ast}(z) \\
      & \sigma_{\theta}^{\ast}(z) - \sigma_{\rho}^{\ast}(z) + 2{\rm{i}}\tau_{\rho\theta}^{\ast}(z) = [\sigma_{y}^{\ast}(z) - \sigma_{x}^{\ast}(z) + 2{\rm{i}}\tau_{xy}^{\ast}(z)]\cdot \left. \frac{\zeta}{\overline{\zeta}} \frac{z^{\prime}[w(\zeta)]}{\overline{z^{\prime}[w(\zeta)]}} \frac{w^{\prime}(\zeta)}{\overline{w^{\prime}(\zeta)}} \right|_{\zeta \rightarrow \zeta[w(z)]}
    \end{aligned}
  \right.
\end{equation}
where $ \sigma_{\theta}^{\ast} $, $ \sigma_{\rho}^{\ast} $, and $ \tau_{\rho\theta}^{\ast} $ denote the final curvilinear hoop, radial, and tangential stress components, respectively.

\section{Numerical cases and verification}
\label{sec:numerical-cases}

To pratically obtain numerical results, the infinite series of the solution of the mixed boundary problem in Eq. (\ref{eq:5.1}) should be truncated into $ 2N_{0}+1 $ items ($ -N_{0} \leq n \leq N_{0} $), and Eq. (\ref{eq:5.5}) should be truncated correspondingly as
\begin{equation}
  \label{eq:5.5a'}
  \tag{5.5a'}
  \varphi^{\prime}_{0}(\zeta) = \sum\limits_{k=-N_{0}}^{N_{0}} A_{k} \zeta^{k}, \quad A_{k} = \sum\limits_{n=-k}^{N_{0}} \alpha_{n+k}d_{-n}, \quad \zeta \in {\bm \omega}_{0}^{+}
\end{equation}
\begin{equation}
  \label{eq:5.5b'}
  \tag{5.5b'}
  \varphi^{\prime}_{0}(\zeta) = \sum\limits_{k=-N_{0}}^{N_{0}} B_{k} \zeta^{k}, \quad B_{k} = \sum\limits_{n=k}^{N_{0}} \beta_{n-k}d_{n}, \quad \zeta \in {\bm \omega}_{0}^{-}
\end{equation}
Subsequently, all the series in the iterative solution procedure in Eqs. (\ref{eq:5.22})-(\ref{eq:5.27}) should be truncated into finite items correspondingly. Similarly, the expansions in Eqs. (\ref{eq:5.14}), (\ref{eq:5.16}), and (\ref{eq:6.2}) should be also truncated into $ 2M+1 $ items ($ -M \leq k \leq M $), where $ M $ should be large for good approximation accuracy. The linear systems of the iterative solution procedure in Eqs. (\ref{eq:5.24})-(\ref{eq:5.26}) is numerically stable with small condition numbers \cite[]{LIN2024396}. Finally, the stress and displacement solution in Eqs. (\ref{eq:6.3}), (\ref{eq:6.4}), (\ref{eq:6.6}), and (\ref{eq:6.7}) would be truncated into finite items as well. Therefore, certain errors of the final stress and displacement in Eqs. (\ref{eq:6.7}) and (\ref{eq:6.6b}) would exist. To be specific, the abrupt change of mixed boundary conditions along the ground surface segement $ {\bm{C}}_{1c0} $ and $ {\bm{C}}_{1f0} $ in Eqs. (\ref{eq:2.6a'}) and (\ref{eq:2.6b'}) can not be perfectly simulated, since the mixed boundary conditions require series of infinite items, instead of series of finite items after truncation \cite[]{LIN2024396}. The Gibbs phenomena would thereby occur in the results of stress and displacement components in the simulating geomaterial \cite[]{LIN2024396,LIN2024106008}. To reduce the Gibbs phenomena, the following Lanczos filtering is applied in Eqs. (\ref{eq:6.3}) and (\ref{eq:6.4}) to replace $ \sigma^{k} $ with $ L_{k}\cdot\sigma^{k} $ \cite[]{Lanczos1956,LIN2024396,LIN2024106008} as
\begin{equation}
  \label{eq:7.1}
  L_{k} = \left\{
    \begin{aligned}
      & \; 1, \quad k = 0 \\
      & \; \sin\left( \frac{k}{N_{0}}\pi \right)/\left( \frac{k}{N_{0}}\pi \right), \quad {\rm{otherwise}}
    \end{aligned}
  \right.
\end{equation}
where $ -N_{0} \leq k \leq N_{0} $.

In this section, several numerical cases are conducted to verify the proposed solution. All numerical cases are computed using the programming code {\texttt{FORTRAN}} of compiler {\texttt{GCC}} (version 13.2.1 20230801). All linear systems are solved using the {\texttt{DGESV}} package of {\texttt{LAPACK}} (version 3.11.0). All data figures are drawn using {\texttt{GNUPLOT}} (version 6.0 patchlevel 0). All the codes are open sourced, and can be found in author Luobin Lin's {\texttt{github}} repository: {\underline{github.com/luobinlin987/over-under-break-excavation-shallow-tunnelling}}.

\subsection{Cavity shapes and bidirectional conformal mappings}
\label{sec:cavity-shapes}

The following asymmetrical cavity shape of four different depths are used as the numerical cases to verify the Charge Simulation Method (CSM) in Appendix \ref{sec:first-step-composite} and the bidirectional composite conformal mapping in Eq. (\ref{eq:3.3}):
\begin{equation}
  \label{eq:7.2}
  \left\{
    \begin{aligned}
      & \left( \frac{x}{4} \right)^{2} + \left( \frac{y+H_{j}}{5} \right)^{2} = 1, \quad x \leq 0 \\
      & \left( \frac{x}{6} \right)^{2} + \left( \frac{y+H_{j}}{5} \right)^{2} = 1, \quad x > 0 \\
    \end{aligned}
  \right.
  , \quad j = 1,2,3,4
\end{equation}
where $ H_{1} = 10, H_{2} = 8, H_{3} = 6 $, $ H_{4} = 5.2 $ denote four different cavity depths. The expected cavity boundaries of the above four cases in Eq. (\ref{eq:7.2}) are fully circular as
\begin{equation}
  \label{eq:7.2'}
  \tag{7.2'}
  \left( \frac{x}{5} \right)^{2} + \left( \frac{y+H_{j}}{5} \right)^{2} = 1, \quad j = 1,2,3,4
\end{equation}
Comparing to Eq. (\ref{eq:7.2'}), the left ($ x \leq 0 $) and right ($ x > 0 $) parts of the cavity in Eq. (\ref{eq:7.2}) can be treated as under-break and over-break regions, respectively. The coordinates of joint points $ T_{1} $ and $ T_{2} $ between constraining ground surface $ {\bm{C}}_{1c} $ and free ground surface $ {\bm{C}}_{1f} $ are selected as $ (-10,0) $ and $ (10,0) $, respectively.

\begin{figure}[htb]
  \centering
  \includegraphics{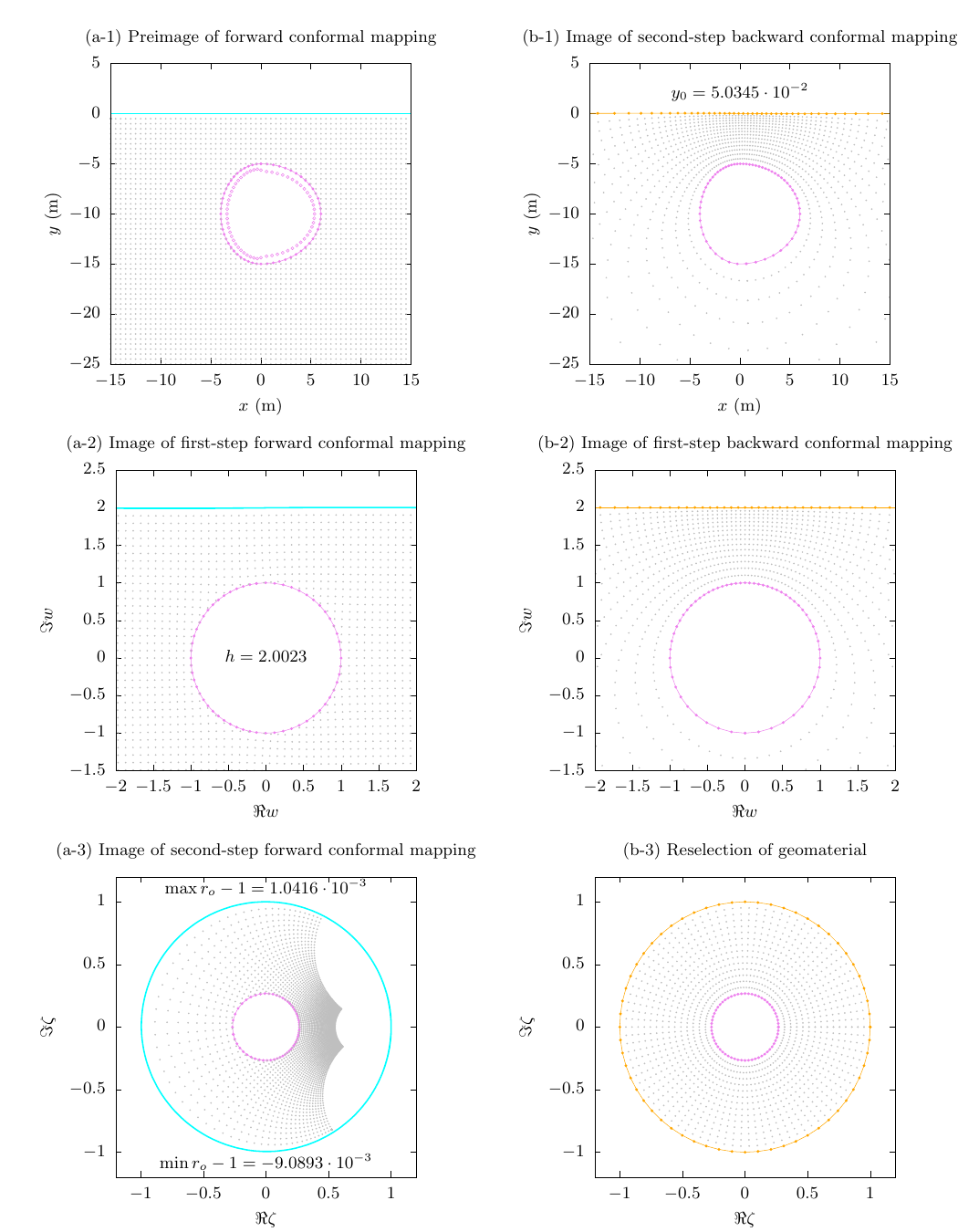}
  \caption{Bidirectional conformal mapping of Case 1}
  \label{fig:4}
\end{figure}

\begin{figure}[htb]
  \centering
  \includegraphics{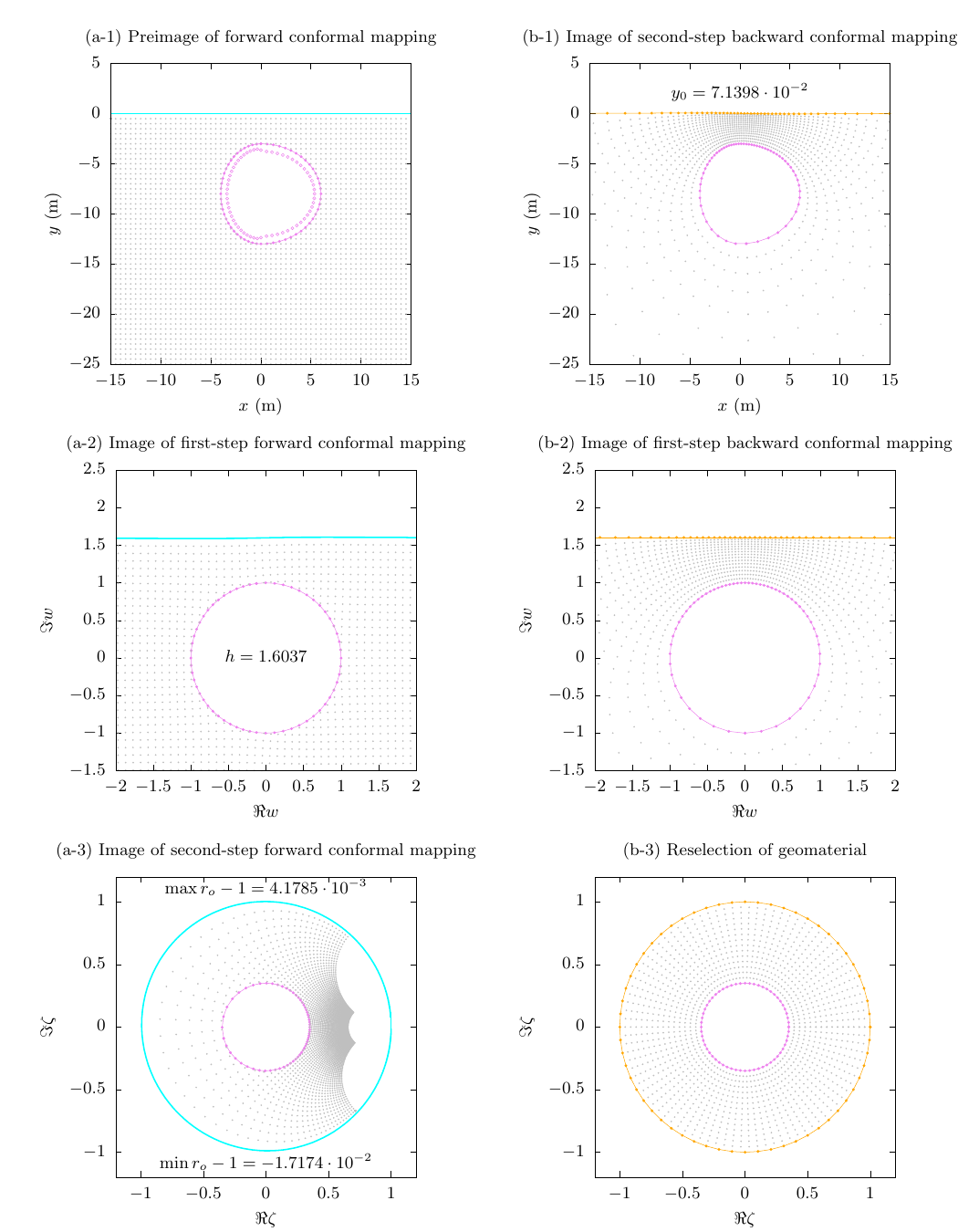}
  \caption{Bidirectional conformal mapping of Case 2}
  \label{fig:5}
\end{figure}

\begin{figure}[htb]
  \centering
  \includegraphics{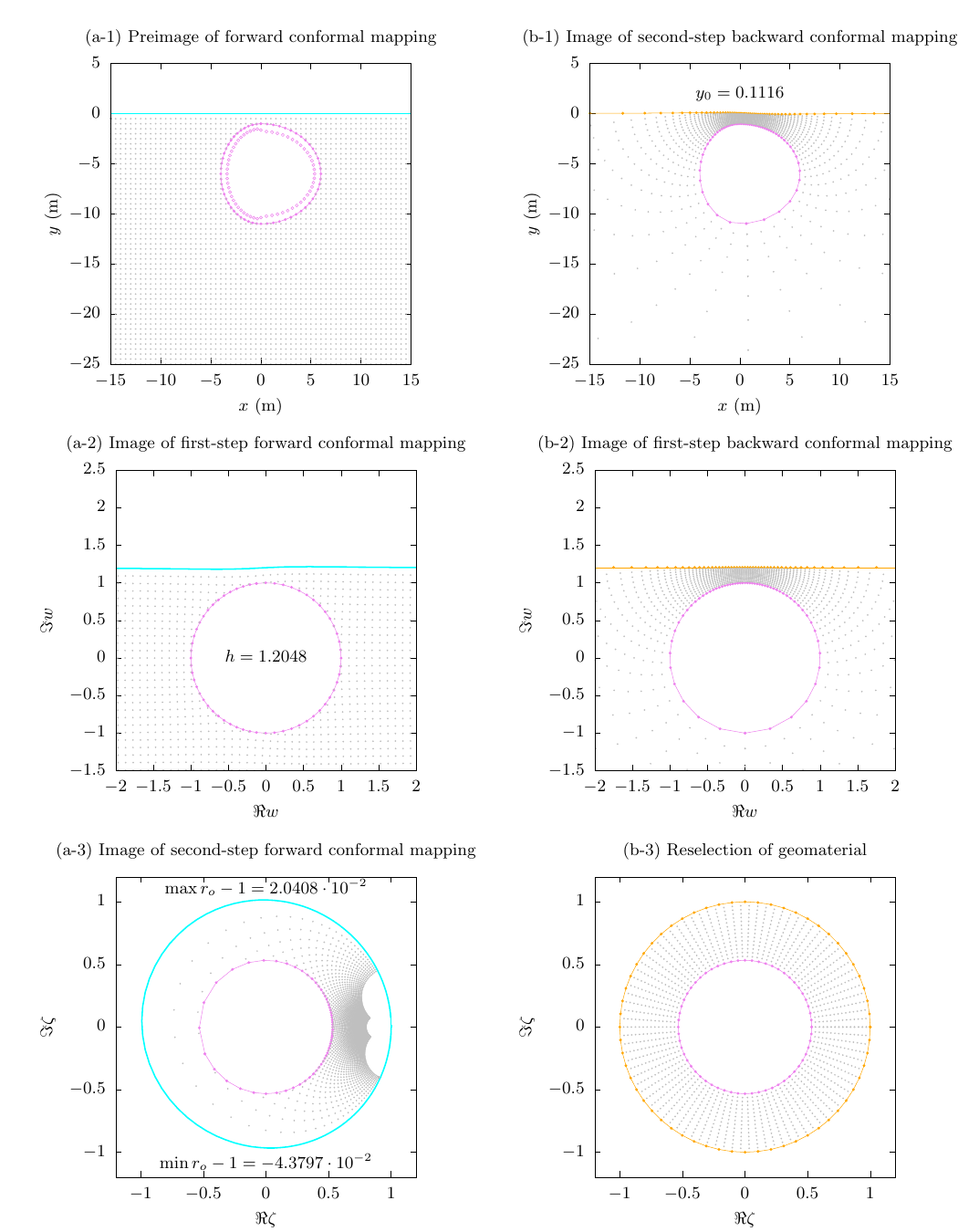}
  \caption{Bidirectional conformal mapping of Case 3}
  \label{fig:6}
\end{figure}

\begin{figure}[htb]
  \centering
  \includegraphics{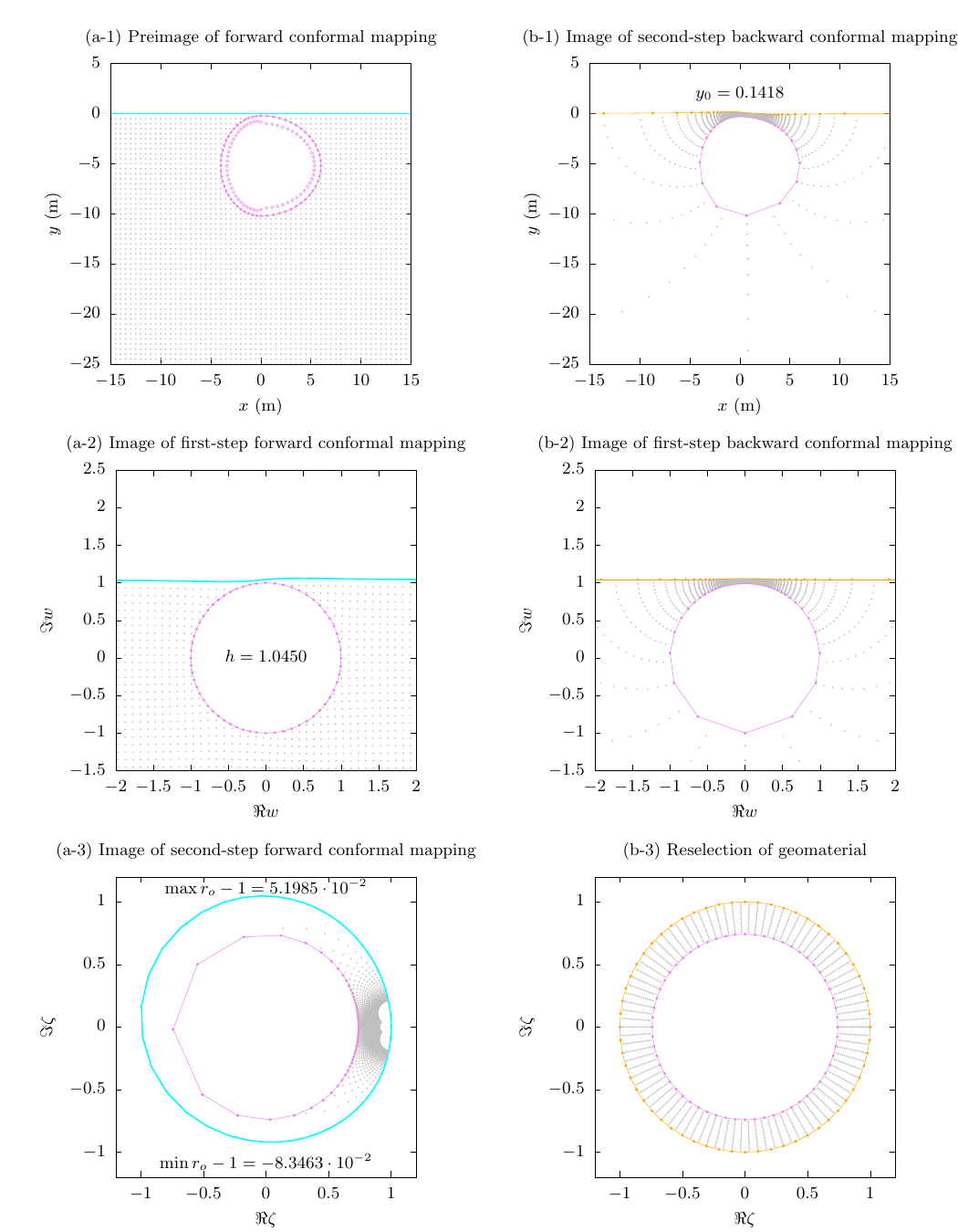}
  \caption{Bidirectional conformal mapping of Case 4}
  \label{fig:7}
\end{figure}

The {\emph{collocation points}} of the four cavity boundaries are selected as
\begin{subequations}
  \label{eq:7.3}
  \begin{equation}
    \label{eq:7.3a}
    \left\{
      \begin{aligned}
        & \Re z_{j,i} = 6\cdot\cos\left(\frac{\pi}{2}-\frac{i-1}{N}\pi\right) \\
        & \Im z_{j,i} = 5\cdot\sin\left(\frac{\pi}{2}-\frac{i-1}{N}\pi\right) - {\rm{i}}H_{j} \\
      \end{aligned}
      ,\quad i = 1,2,3,\cdots,N
    \right.
  \end{equation}
  \begin{equation}
    \label{eq:7.3b}
    \left\{
      \begin{aligned}
        & \Re z_{j,N+i} = 4\cdot\cos\left(-\frac{\pi}{2}-\frac{i-1}{N}\pi\right) \\
        & \Im z_{j,N+i} = 5\cdot\sin\left(-\frac{\pi}{2}-\frac{i-1}{N}\pi\right) - {\rm{i}}H_{j} \\
      \end{aligned}
      ,\quad i = 1,2,3,\cdots,N
    \right.
  \end{equation}
\end{subequations}
where $ N = 30 $ and $ j = 1,2,3,4 $. The {\emph{assignment factor}} in Eq. (\ref{eqa:5}) takes $ k_{2} = 1.2 $. Substituting the {\emph{collocation points}} in Eq. (\ref{eq:7.3}) of the four cavity shapes and $ k_{2} $ into the bidirectional composite conformal mapping in Eq. (\ref{eq:3.3}) and the two-step mapping details in Appendix \ref{sec:bidir-comp-conf} gives the mapping results in Figs. \ref{fig:4}-\ref{fig:7}.

The deployment and logic among subfigures in Figs. \ref{fig:4}-\ref{fig:7} are the same to those in Fig. \ref{fig:3}. In subfigure a-1 of Figs. \ref{fig:4}-\ref{fig:7}, the {\emph{charge points}} are located inside of the cavity boundaries, and the sample points (the light gray points) are filled in the finite region $ -15 \leq x \leq 15, -30 \leq y < 0 $. Correspondingly, only part of subfigure a-3 in Figs. \ref{fig:4}-\ref{fig:7} would be filled with sample points, while the rest region would remain blank as illustrated. In subfigure a-2, the cavity depth $ h = \Im T_{2} $ is presented. In subfigure a-3, the maximum outer radius $ \max r_{o} $ and minimum outer radius $ \min r_{0} $ are illustrated. In subfigure b-3 in Figs. \ref{fig:4}-\ref{fig:7}, the sample points are taken radially and tangentially and are filled within the reselected geomaterial region to distinguish from subfigure a-3. Correspondingly, the sample point distributions in subfigures b-2 and b-1 in Figs. \ref{fig:4}-\ref{fig:7} would be different from subfigures a-2 and a-1, respectively. In subfigure b-1, the maximum absolute vertical coordinate along ground surface $ y_{0} = \max |z|_{x \rightarrow 0} $ is illustrated, which indicates the approximation of the simulating geomaterial in subfigure b-1 to the original geomaterial in subfigure a-1.

Figs. \ref{fig:4}-\ref{fig:7} show that when the cavity depth becomes smaller, the ground surfaces in subfigures a-2 and a-3 would not be a straight line and a circle any longer, respectively. Correspondingly, the reselected geomaterial in subfigure b-3 that {\emph{changes the definition domain of the mapped geomaterial}} causes a non-straight ground surface in subfigure a-1. Such a feature has been emphasized in Section \ref{sec:geom-resel} that the geomaterial in subfigure b-1 is only a simulation to that in subfigure a-1, and now been visually shown in Figs. \ref{fig:4}-\ref{fig:7}.

For Cases 2, 3, and 4, the ground surfaces in subfigure b-1 in Figs. \ref{fig:5}-\ref{fig:7} are obviously curved with relatively large $ y_{0} $, and subfigure b-1 can not simulate subfigure a-1 ideally. For Case 1, the ground surface in subfigure b-1 in Fig. \ref{fig:4} is almost straight with a small $ y_{0} $, which can be treated as a good simulation of subfigure a-1 in Fig. \ref{fig:4}. Subsequently, we may estimate in a rough and intuitive sense that when cavity depth is about two times of cavity size, the two-step bidirectional conformal mapping would provide good simulation. In real-world tunnel and pipeline engineering, such a geometry of cavity depth versus cavity size is commonly seen. Therefore, the four numerical cases indicate that the the bidirectional composite conformal mapping binded by the geomaterial reselection strategy is reasonable and valuable in nonaxisymmetrical cavity excavation problem in geometrical propective for some real-world engineering scenarios.

\subsection{Lanczos filtering}
\label{sec:lanczos-filtering}

\begin{figure}[htb]
  \centering
  \includegraphics{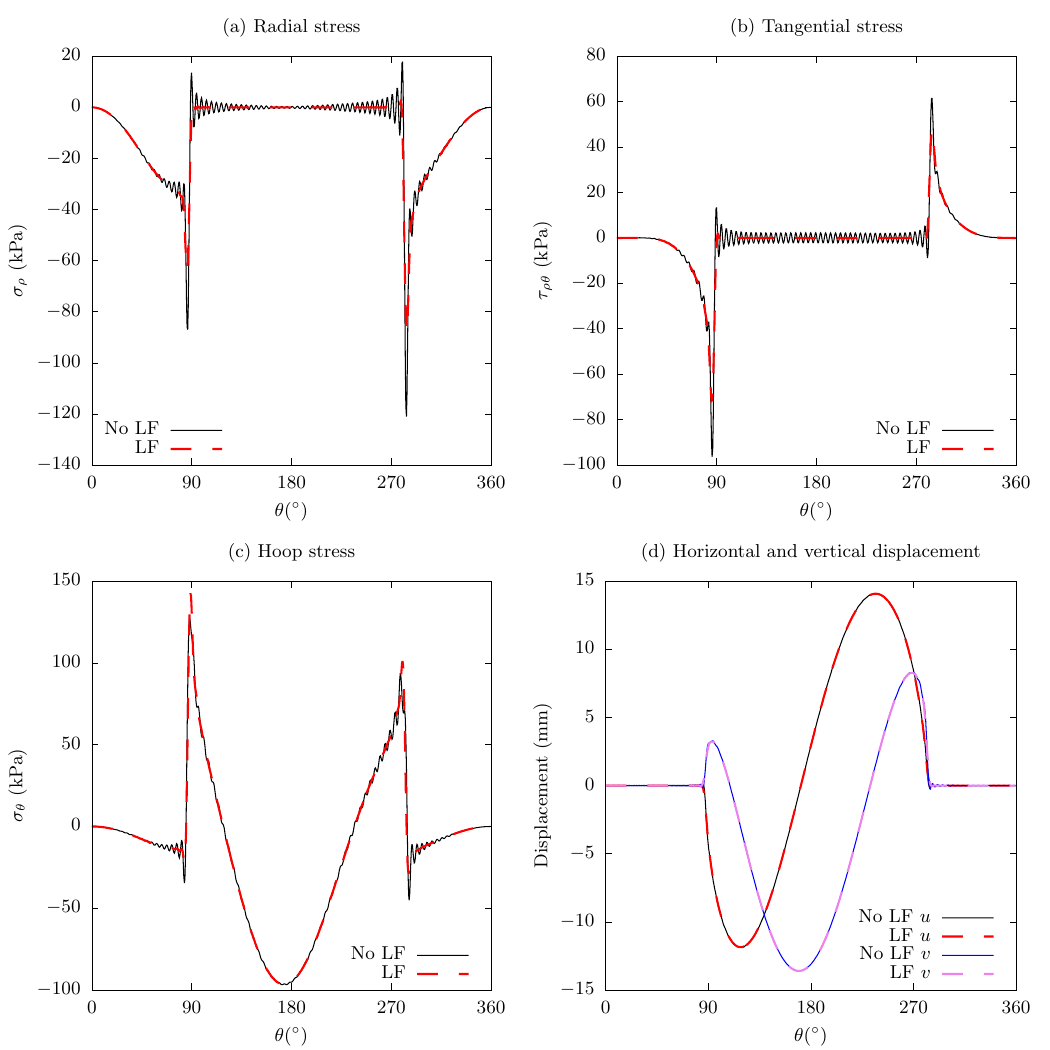}
  \caption{Comparisons of curvilinear stress and displacement components in mapping plane $ \zeta $ of Lanczos filtering (LF) along ground surface of Case 1}
  \label{fig:8}
\end{figure}

\begin{figure}[htb]
  \centering
  \includegraphics{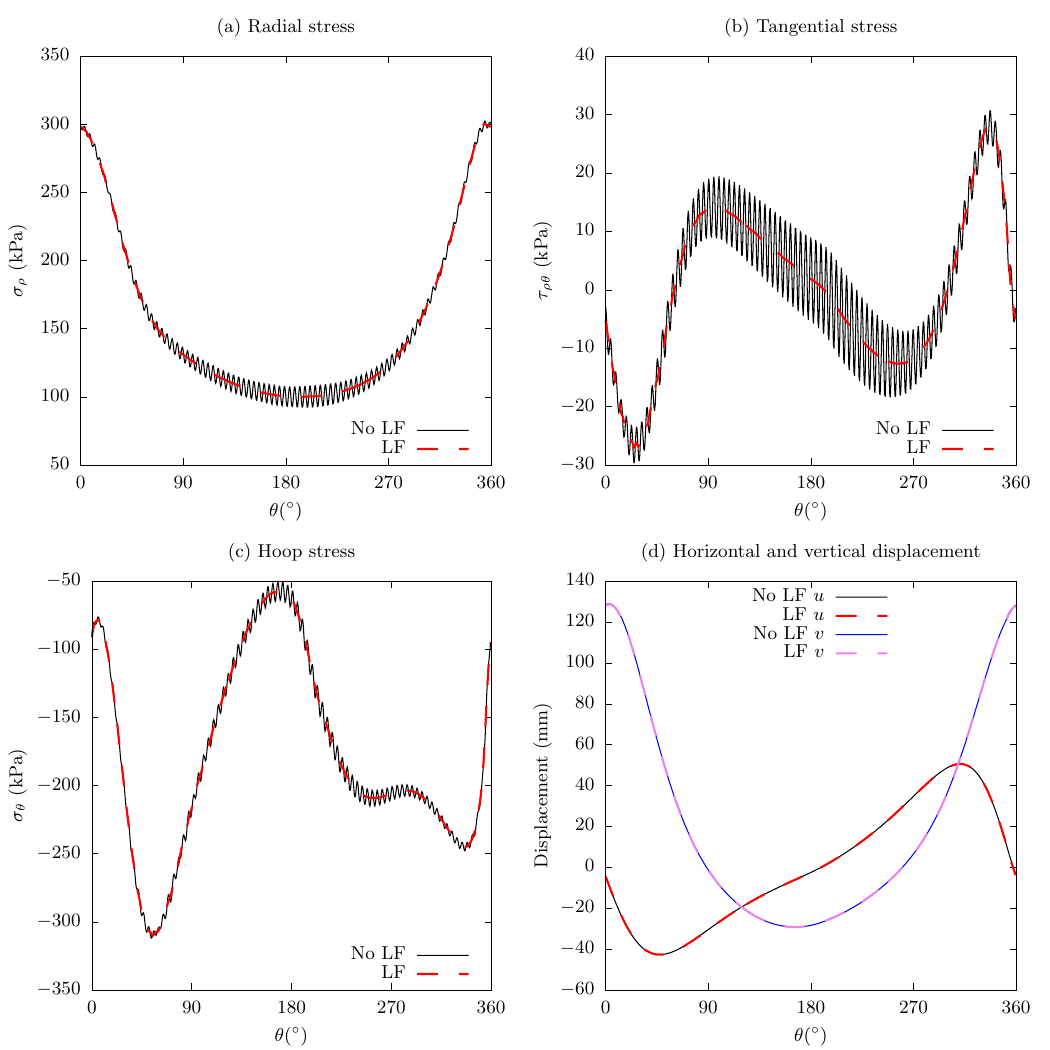}
  \caption{Comparisons of curvilinear stress and displacement components in mapping plane $ \zeta $ of Lanczos filtering (LF) along cavity boundary of Case 1}
  \label{fig:9}
\end{figure}

Last section indicates that only Case 1 is geometrically suitable to proceed further verification of subsequent solution in Sections \ref{sec:probl-transf}-\ref{sec:stress-displ-simul}. The parameters in Table \ref{tab:1} are used in the subsequent solution, and the plane strain condition is considered ($ \kappa = 3-4\nu $). The parameter $ x_{0} = \frac{\Re|T_{1}|}{H_{1}} = \frac{\Re|T_{2}|}{H_{1}} $ is dimensionless, and denotes the normalized range of the free ground surface $ {\bm{C}}_{1f} $. The truncation of the infinite series in Eqs. (\ref{eq:5.14}), (\ref{eq:5.16}), and (\ref{eq:6.2}) take $ M = 360 $ after numerous trial computation by balancing approximation accuracy and consumption time of programming code. The truncation of Eqs. (\ref{eq:5.5a'}) and (\ref{eq:5.5b'}) takes $ N_{0} = 80 $ with good accuracy after trial computations. Substituting the parameters in Table \ref{tab:1} into the subsequent solution gives the curvilinear stress and displacement components without Lanczos filtering in the mapping plane $ \zeta $ along ground surface (Eq. (\ref{eq:6.4})) and cavity boundary (Eq. (\ref{eq:6.3}) when $ \rho = \alpha$), as shown in black lines in Figs. \ref{fig:8} and \ref{fig:9}. After applying the Lanczos filtering in Eq. (\ref{eq:7.1}), the results in red lines can be obtained. The results without Lanczos filtering illustrate clear oscillation with the Gibbs phenomena, while the results with Lanczos filtering is much smoother to greatly reduce the Gibbs phenomena. Similar result comparisons can be found in previous studies \cite[]{LIN2024396,LIN2024106008}. Additionally, Figs. \ref{fig:8} and \ref{fig:9} also indicate that the stress and displacement at infinity in the physical plane $ z $ (corresponding to the curvilinear stress and displacement in the mapping plane $ \zeta $ of $ \theta = 0^{\circ} $ in Figs. \ref{fig:8} and \ref{fig:9}) vanish as expected in the mechanical model in Fig. \ref{fig:2}b.

\begin{table}[htb]
  \centering 
  \caption{Parameters of Case 1}
  \label{tab:1}
  \begin{tabular}{ccccccc}
    \toprule
    $ E $ (MPa) & $ \nu $ & $ \gamma $ (kPa) & $ k_{0} $ & $ x_{0} $ & $ N_{0} $ & $ M $ \\
    \midrule
    20 & 0.3 & 20 & 0.8 & 1 & 80 & 360 \\
    \bottomrule
  \end{tabular}
\end{table}

\subsection{Solution convergence against joint points along ground surface $ x_{0} $ and truncation number $ N_{0} $}
\label{sec:solut-conv-joint}

\begin{figure}[htb]
  \centering
  \includegraphics{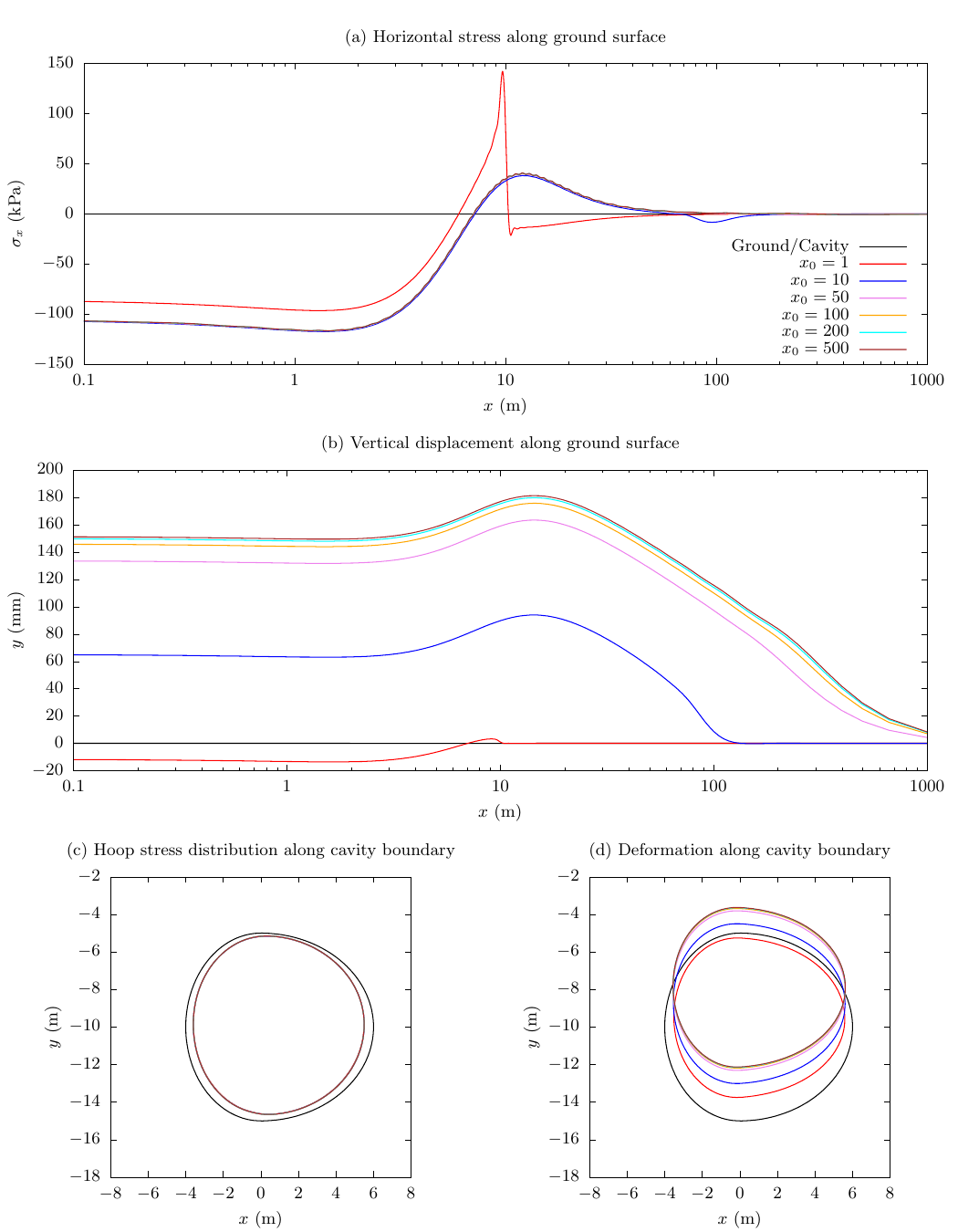}
  \caption{Stress and displacement convergence against normalized free ground surface $ x_{0} $ of Case 1}
  \label{fig:10}
\end{figure}

Logically, we should expect that as the joint points $ T_{1} $ and $ T_{2} $ along ground surface locate far away from the cavity, the stress and displacement components along ground surface and cavity boundary should be gradually convergent. To illustrate the solution convergence of the joint points along ground surface, the parameter $ x_{0} $ of the normalized range of free ground surface takes the values $ x_{0} = 1, 10, 50, 100, 200, 500 $. Substituting the selected $ x_{0} $ and the other parameters in Table \ref{tab:1} into the solution gives Fig. \ref{fig:10}. Fig. \ref{fig:10} illustrates clear convergence of stress and displacement components along ground surface and cavity boundary. Therefore, we should take $ x_{0} = 500 $ in the following computations.

\begin{figure}[htb]
  \centering
  \includegraphics{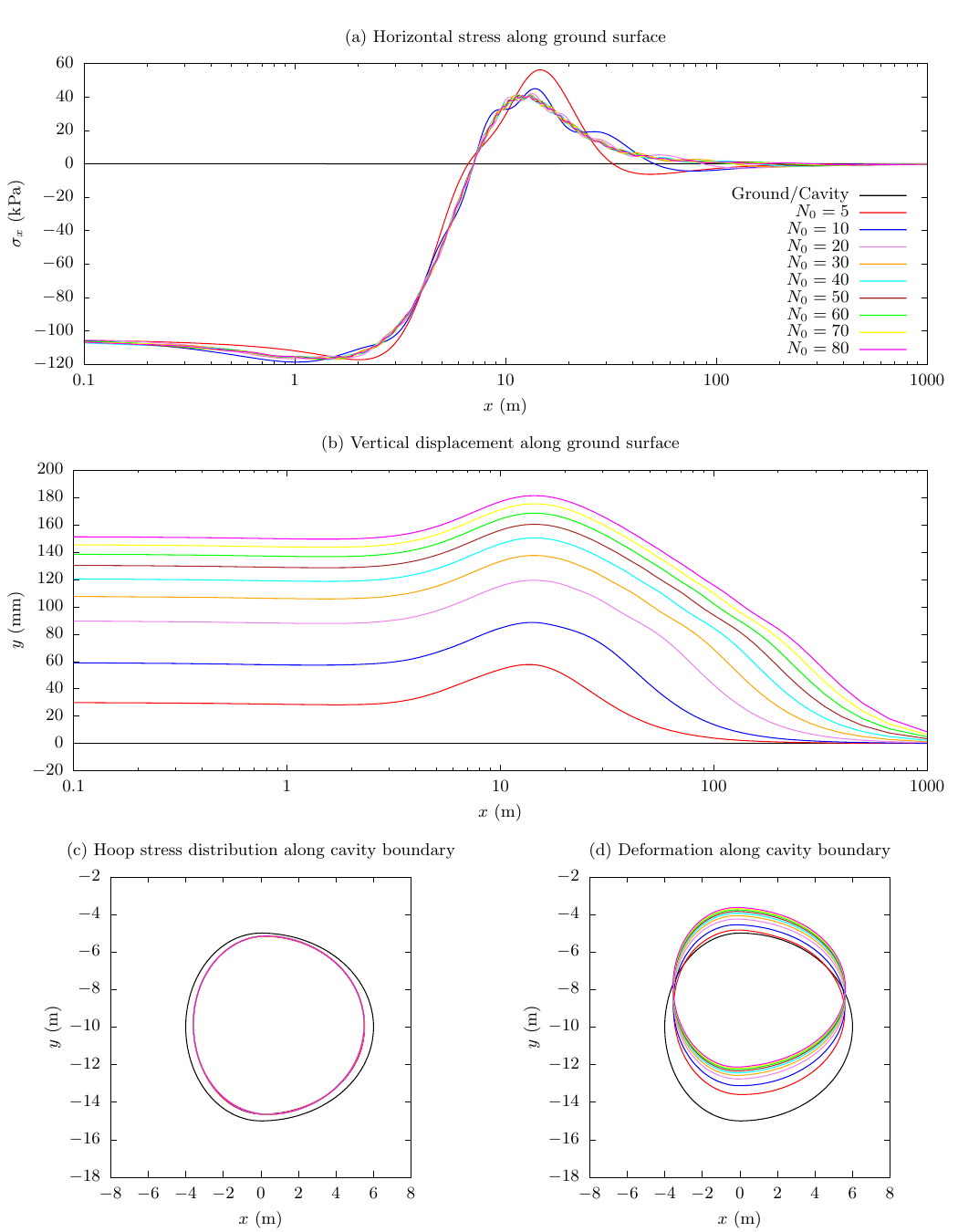}
  \caption{Stress and displacement convergence against truncation number $ N_{0} $ of Case 1}
  \label{fig:11}
\end{figure}

Similarly, the solution should be convergent, when the truncation number $ N_{0} $ in Eqs. (\ref{eq:5.5a'}) and (\ref{eq:5.5b'}) gets larger, but should not exceed the threshold $ 2N_{0}+1 \leq M $ due to Eq. (\ref{eq:5.26}). Correspondingly, the truncation parameter $ N_{0} $ takes the values $ N_{0} = 5, 10, 20, 30, 40, 50, 60, 70, 80 $. Substituting the selected $ N_{0} $ and other parameters in Table \ref{fig:1} ($ x_{0} = 500 $ as discussed above) into the solution gives Fig. \ref{fig:11}. Similar to Fig. \ref{fig:10}, Fig. \ref{fig:11} also illstrates clear convergence of stress and displacement components along ground surface and cavity boundary. Therefore, the value $ N_{0} = 80 $ is suitable (not too large or too small) to meet accuracy requirement.

\subsection{Comparisons with finite element solution}
\label{sec:comp-with-finite}

\begin{figure}[htb]
  \centering
  \begin{tikzpicture}
    \fill [pattern = north east lines] (-6.2,0) rectangle (6.2,-6.2);
    \fill [gray!30] (-6,0) rectangle (6,-6);
    \draw [line width = 1.5pt] (-6,0) -- (-6,-6) -- (6,-6) -- (6,0);
    \draw [line width = 1.5pt, cyan] (-6,0) -- (6,0);
    \draw [dashed, violet] (0,-2) ellipse [x radius = 1.2, y radius = 1];
    \draw [dashed, blue] (0,-2) ellipse [x radius = 0.8, y radius = 1]; 
    \draw [dashed] (-6,-4) -- (6,-4);
    \draw [dashed] (-2,0) -- (-2,-6);
    \draw [dashed] (2,0) -- (2,-6);
    \node at (-4,0) [align = center] {
      \scriptsize{88.82m} \\
      \scriptsize{\textcolor{red}{80 seeds}}
    };
    \node at (-4,-4) [align = center] {
      \scriptsize{\textcolor{red}{80 seeds}} \\
      \scriptsize{88.82m}
    };
    \node at (-4,-6) [align = center] {
      \scriptsize{\textcolor{red}{80 seeds}} \\
      \scriptsize{88.82m}
    };
    \node at (4,-6) [align = center] {
      \scriptsize{\textcolor{red}{80 seeds}} \\
      \scriptsize{88.82m}
    };
    \node at (4,-4) [align = center] {
      \scriptsize{\textcolor{red}{80 seeds}} \\
      \scriptsize{88.82m}
    };
    \node at (4,0) [align = center] {
      \scriptsize{88.82m} \\
      \scriptsize{\textcolor{red}{80 seeds}}
    };
    \node at (-6,-5) [rotate = 90, align = center] {
      \scriptsize{72.36m} \\
      \scriptsize{\textcolor{red}{70 seeds}}
    };
    \node at (-2,-5) [rotate = 90, align = center] {
      \scriptsize{72.36m} \\
      \scriptsize{\textcolor{red}{70 seeds}}
    };
    \node at (2,-5) [rotate = 90, align = center] {
      \scriptsize{\textcolor{red}{70 seeds}} \\
      \scriptsize{72.36m}
    };
    \node at (6,-5) [rotate = 90, align = center] {
      \scriptsize{\textcolor{red}{70 seeds}} \\
      \scriptsize{72.36m}
    };
    \node at (-6,-2) [rotate = 90, align = center] {
      \scriptsize{27.64m} \\
      \scriptsize{\textcolor{red}{30 seeds}}
    };
    \node at (-2,-2) [rotate = 90, align = center] {
      \scriptsize{27.64m} \\
      \scriptsize{\textcolor{red}{30 seeds}}
    };
    \node at (2,-2) [rotate = 90, align = center] {
      \scriptsize{\textcolor{red}{30 seeds}} \\
      \scriptsize{27.64m}
    };
    \node at (6,-2) [rotate = 90, align = center] {
      \scriptsize{\textcolor{red}{30 seeds}} \\
      \scriptsize{27.64m}
    };
    \node at (0,-6) [align = center] {
      \scriptsize{\textcolor{red}{30 seeds}} \\
      \scriptsize{22.36m}
    };
    \node at (0,-4) [align = center] {
      \scriptsize{\textcolor{red}{30 seeds}} \\
      \scriptsize{22.36m}
    };
    \node at (0,0) [above] {\scriptsize{22.36m} \enspace \scriptsize{\textcolor{red}{60 seeds}}};
    \draw [->, blue] (0,-2) -- (0,-1);
    \draw [->, blue] (0,-2) -- (-0.8,-2);
    \draw [->, violet] (0,-2) -- (1.2,-2);
    \node at (-0.4,-2) [above, blue] {\scriptsize{4m}};
    \node at (0.6,-2) [above, violet] {\scriptsize{6m}};
    \node at (0,-1.5) [above, rotate = 90, blue]  {\scriptsize{5m}};
    \fill (0,-2) circle [radius = 0.05];
    \draw [<->, violet] (0,0) -- (0,-1);
    \node at (0,-0.5) [above, rotate = 90, violet] {\scriptsize{5m}};
    \node at (0,-3) [below] {\scriptsize{\textcolor{violet}{120 seeds}} \enspace \scriptsize{\textcolor{blue}{120 seeds}}};
    \node at (0,-2.3) [blue] {\uppercase\expandafter{\romannumeral3}};
    \node at (0.95,-2.3) [violet] {\uppercase\expandafter{\romannumeral2}};
    \node at (-0.95,-2.3) [violet] {\uppercase\expandafter{\romannumeral1}};
  \end{tikzpicture}
  \caption{Schematic diagram of FEM model and meshing strategy}
  \label{fig:12}
\end{figure}
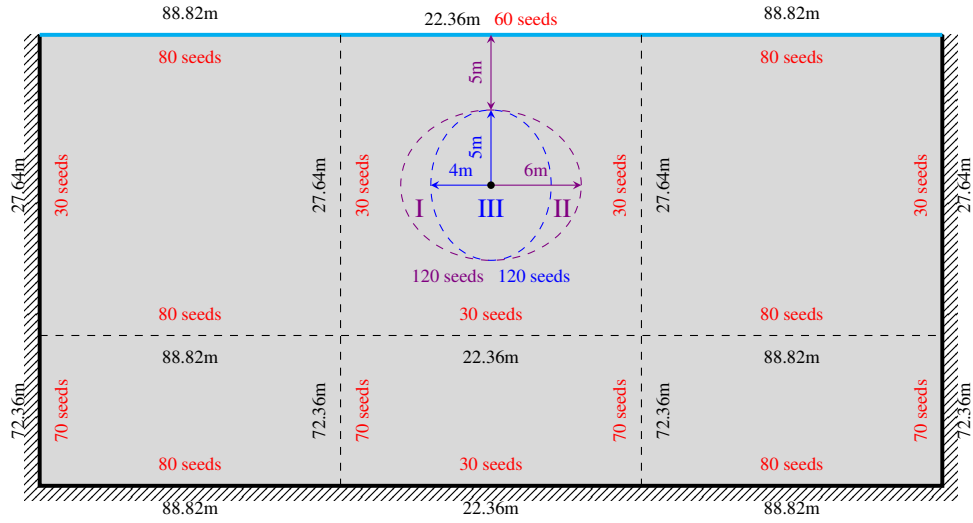

Both of the Lanczos filtering and solution convergence discussed above could not necessarily illustrate the correctness of the proposed solution, and the finite element solution should be conducted for comparisons to fortify the verification. The FEM solution is conducted on FEM software {\texttt{ABAQUS 2016}}. The mechanical parameters of the geomaterial in Case 1 ($ E $, $ \nu $, $ \gamma $, and $ k_{0} $ in Table \ref{tab:1}) are used in the corresponding FEM solution. The geometry of the FEM model is briefly shown in Fig. \ref{fig:12}. The left, right, and bottom boundaries of the geomaterial are constrained, while the top boundary is left free. Two ellipses are sketched using {\texttt{Partition Face}} of {\texttt{Part}} module to facilitate model construction procedure, and three regions marked as \uppercase\expandafter{\romannumeral1}, \uppercase\expandafter{\romannumeral2}, and \uppercase\expandafter{\romannumeral3} are partitioned, as shown in Fig. \ref{fig:12}. The meshing strategy is also shown in Fig. \ref{fig:12}, and 22990 standard elements of plane strain are generated. The three steps of the FEM solution are listed in Table \ref{tab:2}. Meanwhile, the proposed solution takes $ x_{0} = 500 $, $ N_{0} = 80 $, and $ M = 360 $, as verified in previous sections.

\begin{figure}[htb]
  \centering
  \includegraphics[scale = 0.4]{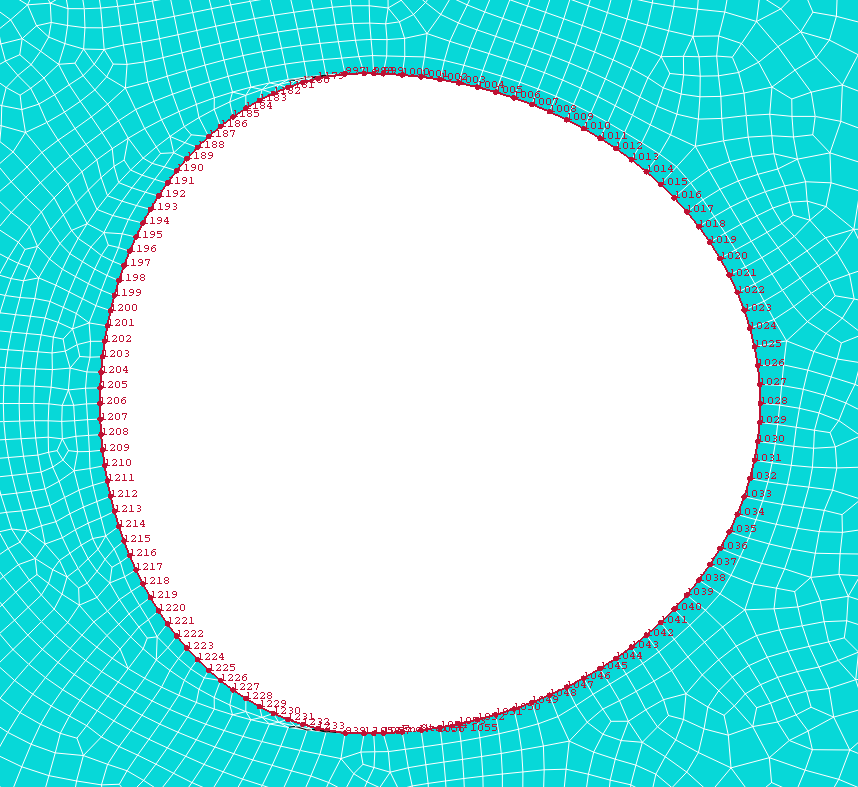}
  \caption{Meshing near cavity after Step 3 (see Table \ref{tab:2}) and manually selected {\texttt{Path}} along cavity boundary}
  \label{fig:13}
\end{figure}

\begin{figure}[htb]
  \centering
  \includegraphics{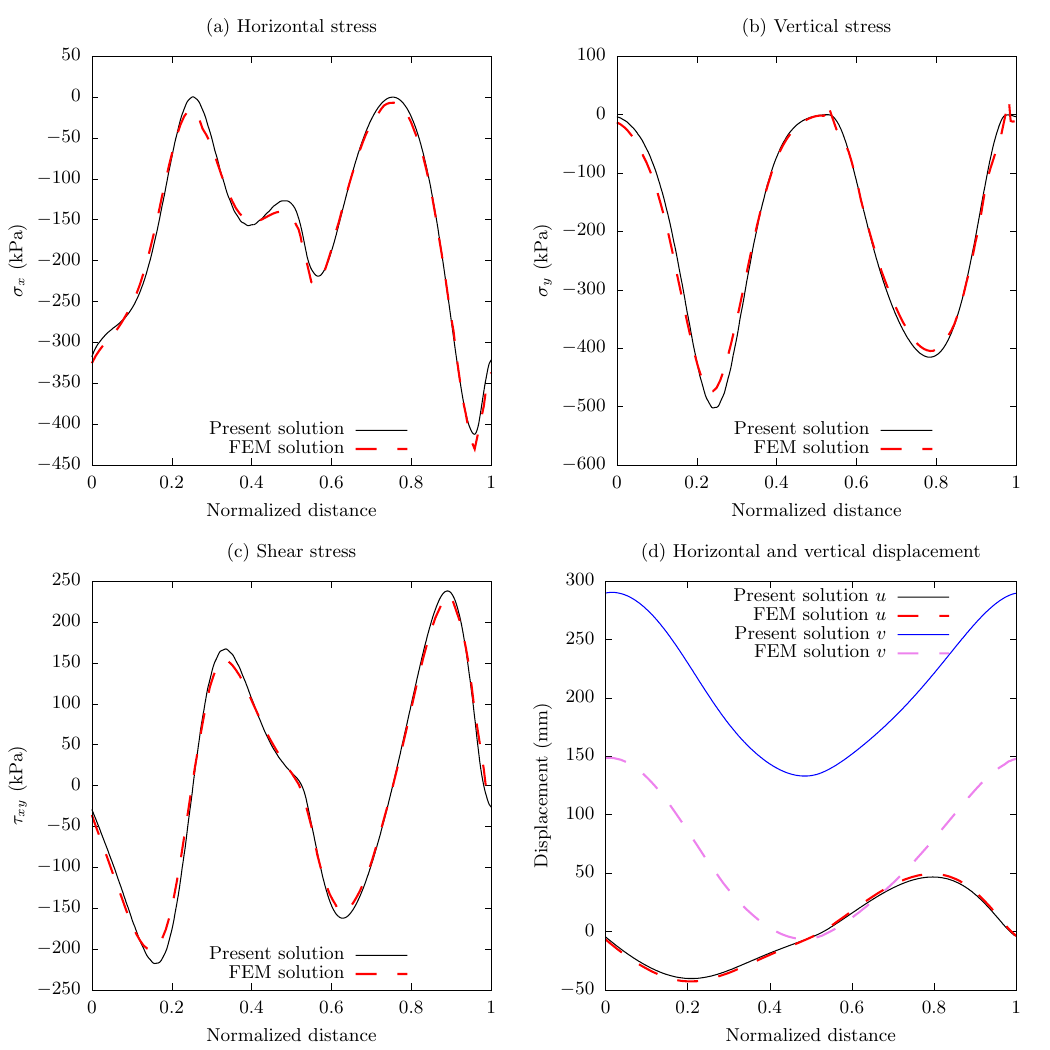}
  \caption{Stress and displacement comparisons along manually selected {\texttt{Path}} (see Fig. \ref{fig:13}) between present solution and FEM solution}
  \label{fig:14}
\end{figure}

The stress and displacement components along cavity boundary are expected to vary greatly, and can be chosen as indexes for comparisons of verification. Since the cavity shape is irregular, the {\texttt{Path}} to extract stress and displacement datum of FEM solution should be manually determined by adding nodes along cavity boundary one by one, as shown in Fig. \ref{fig:13}. The comparisons of the stress and displacement components extracted from the manually selected {\texttt{Path}} of the FEM solution and the ones computed via the proposed solution are illustrated in Fig. \ref{fig:14}. The abscissa axis uses normalized distance to faciliate code programming of analytical solution. Clearly, all three stress components and horizontal displacement components along cavity boundary between two solutions are in good agreements, while the vertical displacement components along cavity boundary between two solutions are similar in shape with a constant difference. Such a constant difference is reasonable, since the geomaterial region of the FEM solution is finite ($ 200{\rm{{m}}}\times 100{\rm{m}} $ in size), while the proposed solution assumes that the geomaterial is infinite. The FEM verification in Ref \cite[]{LIN2024396} has demonstrated that size increase of geomaterial region would make the vertical displacement obtained by the FEM solution gradually approach the one obtained by the analytical solution. In other words, the vertical displacement along cavity boundary obtained by the proposed solution is the limit of the one obtained by the FEM solution, when the geomaterial region is infinitely large. Therefore, the stress and displacement comparisons in Fig. \ref{fig:14} verify the proposed solution.

\begin{table}[htb]
  \centering
  \caption{Steps of finite element solution}
  \label{tab:2}
  \begin{tabular}{cccc}
    \toprule
    Step & Step 1 & Step 2 &  Step 3 \\
    \midrule
    Procedure & (Initial) & Geostatic & Static, General \\
    Load &  \makecell[c]{Applying constraints \\ and geostatic stress} & Applying gravity & - \\
    Interaction & - & - & \makecell[c]{Deactivating regions \\ \uppercase\expandafter{\romannumeral2} and \uppercase\expandafter{\romannumeral3} in Fig. \ref{fig:12}} \\
    \bottomrule
  \end{tabular} 
\end{table}

\subsection{Comparisons with Lin's solution}
\label{sec:comp-with-lins}

The case verifications above only take a nonaxisymmetrical shallow cavity of Case 1 as examples, and no other cavity shapes are considered. In this section, more complicated and axisymmetrical cavity shapes studied by Lin et al. in Ref \cite[]{LIN2024106008} are considered, and the general formation of the backward conformal mapping that maps a unit annulus onto a lower half plane containing an axisymmetrical cavity \cite[]{Zengguisen2019,lu2021complex} can be expressed as
\begin{equation}
  \label{eq:7.4}
  z(\zeta) = -{\rm{i}}a\frac{1+\zeta}{1-\zeta} + {\rm{i}}\sum\limits_{k = 1}^{K} b_{k}(\zeta^{k}-\zeta^{-k})
\end{equation}
where $ a $ and $ b_{k} $ are both real coefficients that control the cavity shape. The inner radius of the unit annulus in Eq. (\ref{eq:7.3}) is also denoted by $ \alpha $. Thus, the {\emph{collocation points}} of the first-step forward conformal mapping in Appendix \ref{sec:first-step-composite} are selected as
\begin{equation}
  \label{eq:7.5}
  z_{i} = -{\rm{i}}a\frac{1+\alpha\sigma_{i}}{1-\alpha\sigma_{i}} + {\rm{i}}\sum\limits_{k = 1}^{K} b_{k}(\alpha^{k}\sigma_{i}^{k} - \alpha^{-k}\sigma_{i}^{-k}), \quad i = 1,2,3,\cdots,2N
\end{equation}
where
\begin{equation}
  \label{eq:7.6}
  \sigma_{i} = \exp \left( {\rm{i}}\frac{i-1}{N}\pi \right)
\end{equation}
The {\emph{assignment factor}} takes $ k_{2} = 0.8 $, $ N = 30 $ as that in Eq. (\ref{eq:7.3}).

\begin{figure}[htb]
  \centering
  \includegraphics{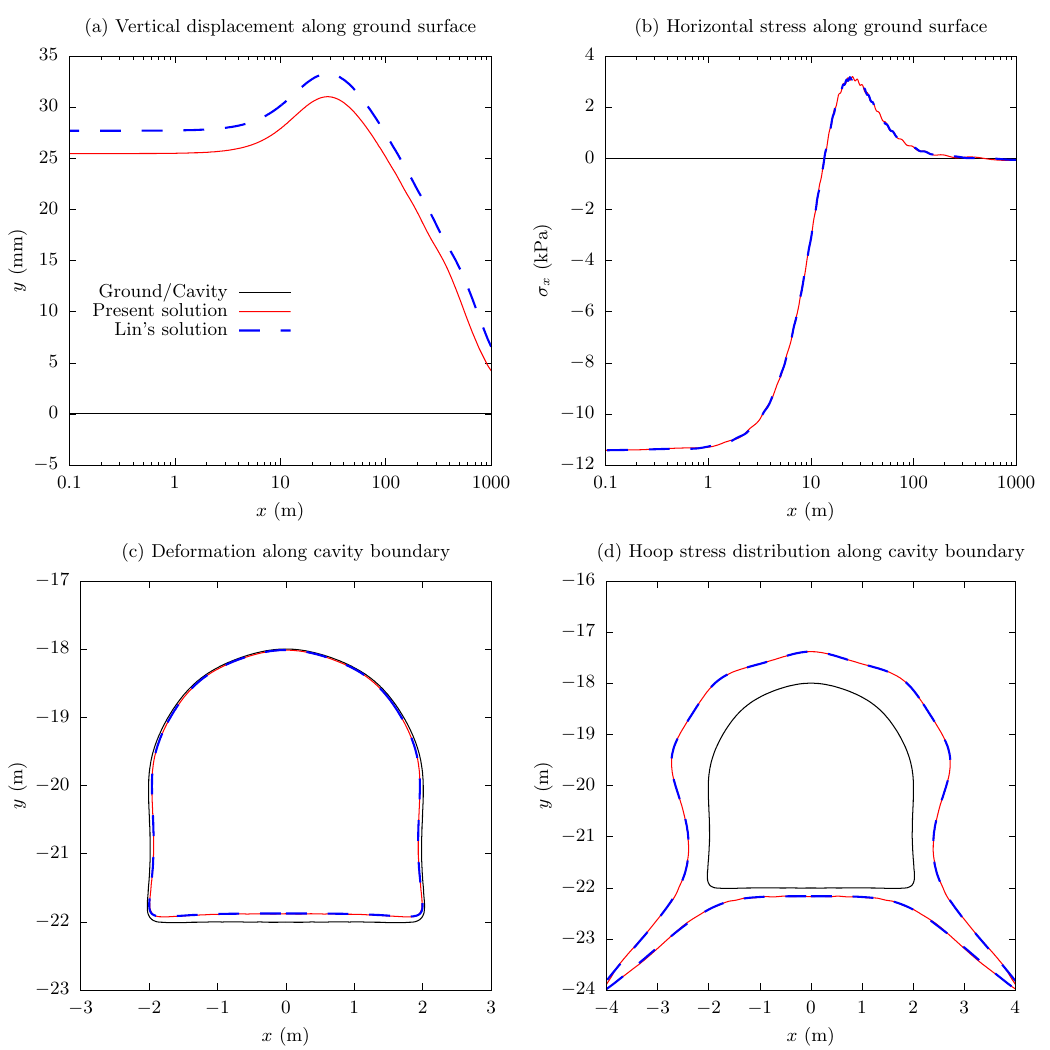}
  \caption{Stress and displacement comparisons between present solution and Lin's solution on cavity shape of Case 2 \cite[]{LIN2024106008}}
  \label{fig:15}
\end{figure}

\begin{figure}[htb]
  \centering
  \includegraphics{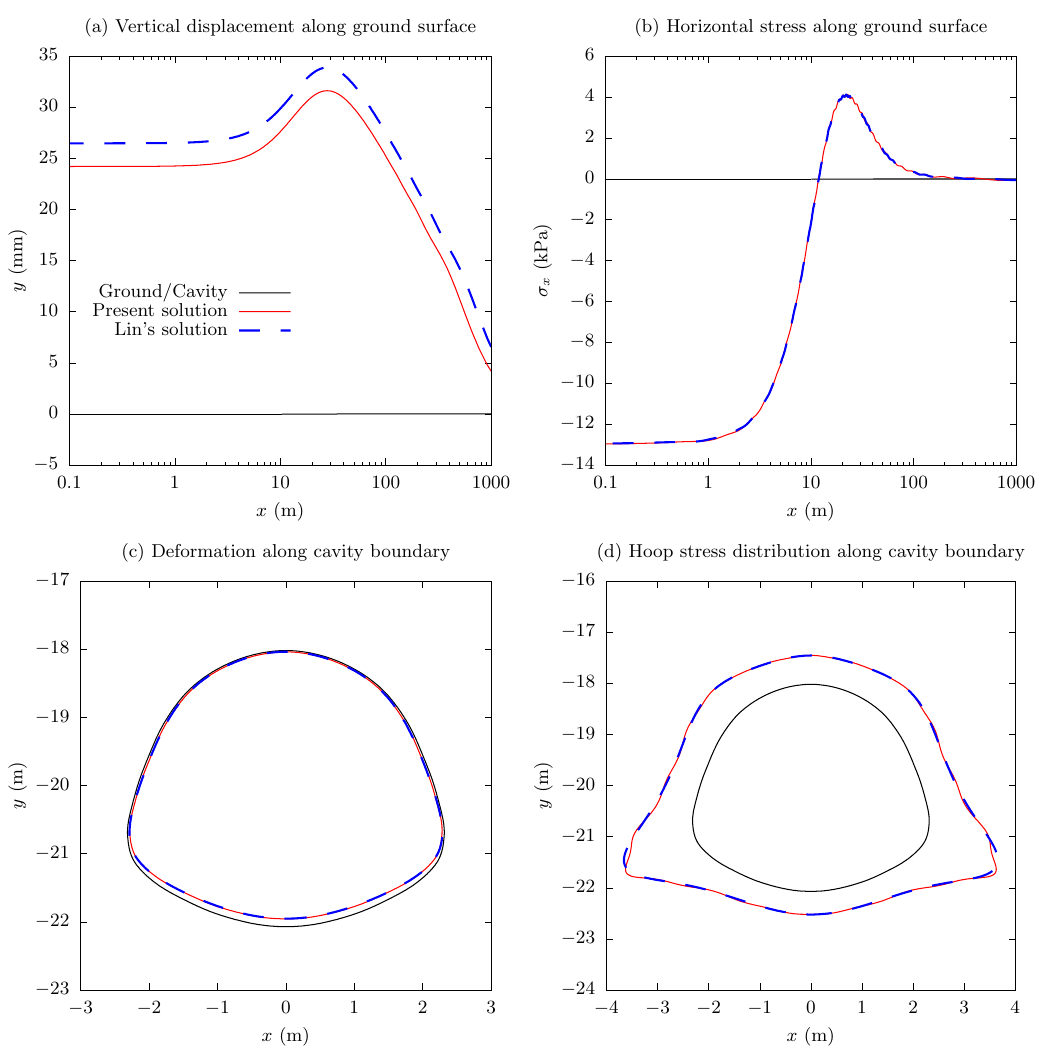}
  \caption{Stress and displacement comparisons between present solution and Lin's solution on cavity shape of Case 4 \cite[]{LIN2024106008}}
  \label{fig:16}
\end{figure}

In this section, we choose cavity shapes of Cases 2 and 4 in Table 1 of Ref \cite[]{LIN2024106008}, so that the real coefficients in Eq. (\ref{eq:7.4}) are determined, as well as the {\emph{collocation points}} in Eq. (\ref{eq:7.5}). The mechanical parameters of the geomaterial are the same to those in Table \ref{tab:1}, and we take $ x_{0} = 500 $, $ N_{0} = 80 $, and $ M = 360 $, just as in Section \ref{sec:comp-with-finite}. Substituting all necessary parameters into the proposed solution and Lin's solution in Ref \cite[]{LIN2024106008} (we only need to change the input parameters of corresponding open sourced {\texttt{.f90}} files released in author Luobin Lin's \texttt{Github} repository {\underline{github.com/luobinlin987/noncircular-shallow-tunnelling-reasonable-displacement}}) gives Figs. \ref{fig:15} and \ref{fig:16}, where the cavity boundary is obtained using the bidirectional composite conformal mapping in Eq. (\ref{eq:3.3}). Figs. \ref{fig:15} and \ref{fig:16} show that the both stress and displacement components along ground surface and cavity boundary are in good agreements between two analytical solutions, except for a small constant difference of vertical displacement. Therefore, Figs. \ref{fig:15} and \ref{fig:16} indicate that the proposed solution are capable of dealing with more complicated shapes in shallow cavity excavation in gravitational geomaterial.

\section{Further discussions}
\label{sec:further-discussion}

Three acute defects have been put forward in our previous solution of Ref \cite[]{LIN2024106008} that (1) no bidirectional conformal mapping is found; (2) the solution is not suitable for asymmetrical cavity, which is commonly seen in real-world engineering; (3) high cost of solving the coefficients in the backward conformal mapping in Eq. (\ref{eq:7.4}). In this paper, all three problems are solved by constructing a conformal mapping in Eq. (\ref{eq:3.3}). Firstly, the conformal mapping in Eq. (\ref{eq:3.3}) is bidirectional. Secondly, such a conformal mapping is suitable for asymmetrical cavity in a composite manner by introducing the Charge Simulation Method (CSM) in Appendix \ref{sec:first-step-composite}. Finally, the computation in both steps of the bidrectional composite conformal mapping (forwardly or backwardly) only relies on solving very simple linear systems, which is low-cost in both computation method and computer infrastructures. Therefore, the proposed solution depending on the bidrectional composite conformal mapping in Eq. (\ref{eq:3.3}) is highly efficient with good accuracy as shown in the verification section.

Despite the many advantages of the proposed solution, its certain defects should also be disclosed. The most obvious defect of the proposed solution is the bold geomaterial simulation discussed in Section \ref{sec:geom-resel} that we use a real unit annulus $ {\bm{\omega}}_{0} $ to simulate the originally mapped geomaterial in a quasi unit annulus region $ {\bm{\omega}} $. The geomaterial simulation scheme {\emph{changes the definition domain of the mapped geomaterial}}, which has been mentioned in the third paragraph in Section \ref{sec:geom-resel}. Subsequently, the original boundary conditions are also simulated by corresponding new ones, as discussed in Section \ref{sec:simul-bound-cond}. To guarantee the accuracy of the proposed solution, the cavity depth should not be too small to ensure that the quasi unit annulus region $ {\bm{\omega}} $ is as approximate as possible to the real unit annulus $ {\bm{\omega}}_{0} $, so that the geometry and boundary conditions in simulating geomaterial in Fig. \ref{fig:3}b-1 would be as approximate as possible to those in the original geomaterial in Fig. \ref{fig:3}a-1. In the perspective of practical application, the depth of cavity excavation is generally not too small in real-world engineering to guarantee the stability of the above geomaterial and ground surface, thus, the geomaterial simulation would not greatly affect the accuracy of the proposed solution.

Comparing to the formal and analytical defect discussed above, the other defect of the proposed solution is latent and numerical that the proposed bidirectional composite conformal mapping is not capable of dealing with cavity with sharp corners. To be specific, the CSM in the first-step forward conformal in Eq. (\ref{eqa:1}) can not deal with cavity with sharp corners, so that no case verification of strict rectangular or vertical-wall cavity shape with sharp corners is conducted. In practical application, sharp corners of a cavity are generally avoided in real-world engineering to reduce severe stress concentration, thus, the CSM in Eq. (\ref{eqa:1}) that is uncapable of dealing with sharp corners is still adequate in the proposed solution.

\section{Conclusions}
\label{sec:conclusions}

In this paper, a new complex variable solution of over-/under-break cavity excavation in gravitational geomaterial with reasonable far-field displacement is proposed by integrating a bidirectional composite conformal mapping sequentially consisting of Charge Simulation Method (CSM) and Verruijt's mapping. When cavity depth is not too small, a geomaterial simulation scheme is applied in the proposed solution, and the mixed boundary conditions are correspondingly simulated as well and further turned into a homogenerous Riemann-Hilbert problem. The Riemann-Hilbert problem is then iteratively solved to obtain infinite complex potential series, which lead to the final stress and displacement components in geomaterial. The infinite complex potential series of the solution are truncated to obtain numerical results with rectification of Lanczos filtering. Several numerical cases are conducted to verify the bidrectional conformal mapping and the subsequent complex variable solution, including comparisons with corresponding finite element solution and an existing analytical solution. The results are all in good agreements and indicate that the proposed complex variable solution is capable of dealing with over-/under-break shallow tunnelling in gravitational geomaterial with reasonable far-field displacement.

\clearpage
\section*{Acknowledgement}
\label{sec:acknowledgement}

This study is financially supported by the Natural Science Foundation of Fujian Province, China (Grant No. 2022J05190 and 2023J01938), the Scientific Research Foundation of Fujian University of Technology (Grant No. GY-Z20094 and GY-Z21026), and the National Natural Science Foundation of China (Grant No. 52178318). The authors would like to thank Professor Changjie Zheng and Ph.D. Yiqun Huang for their suggestions on this study.

\clearpage
\appendix

\section{Bidirectional composite conformal mapping}
\label{sec:bidir-comp-conf}

\subsection{First step of composite conformal mapping}
\label{sec:first-step-composite}

The first-step forward conformal mapping in Eq. (\ref{eq:3.1a}) temporarily ignores the upper boundary $ {\bm{C}}_{1} $ and only focuses on the mapping of the cavity boundary $ {\bm{C}}_{2} $, and laterally maps the simply-connected exterior region outside of cavity boundary $ {\bm{C}}_{2} $ in the global physical plane $ z = x + {\rm{i}} y $ onto its image $ \overline{\bm{\varOmega}}^{\prime} $ in the corresponding local mapping plane $ w = \Re w + {\rm{i}} \Im w $. To obtain such forward conformal mappings, the Charge Simulation Method (CSM) \cite[]{amano1994charge} is applied. The CSM has remarkable features, such as (a) no iteration needed in computation, (b) simplicity of mathematical formation, and (c) low cost of code programming.

The general formula of CSM to map the exterior region $ \overline{\bm{\varOmega}} $ onto its image $ \overline{\bm{\varOmega}}^{\prime} $ can be expressed as
\begin{equation}
  \label{eqa:1}
  w(z) = (z-z_{c}) \cdot \exp \left[ \varGamma + \sum\limits_{k=1}^{2N} Q_{k}\ln(z -Z_{k}) \right], \quad z \in \overline{\bm{\varOmega}}
\end{equation}
where $ \varGamma $ is the Robin constant of region $ \overline{\bm{\varOmega}} $ \cite[]{symm1967numerical}, and $ \varGamma = \lim\limits_{z\rightarrow\infty} \ln|w^{\prime}(z)| $; $ Z_{k} $ are called {\em{charge points}}, which are located outside of and near boundary $ {\bm{C}}_{2} $; $ 2N+1 $ denotes the quantity of charge points; $ Q_{k} $ are called {\emph{charges}}, which are real and should be determined by {\emph{collocation conditions}}; the natural logarithmic items can be separated as
\begin{equation}
  \label{eqa:2}
  \ln(z-Z_{k}) = \ln|z-Z_{k}| + {\rm{i}}\arg(z-Z_{k})
\end{equation}
wherein $ \arg $ denotes the principal value of argument. Different from the general formation of the exterior conformal mapping in Ref \cite[]{amano1994charge}, which is a linear simplification of the integral conformal mappings proposed by Gaier \cite[]{gaier1976integralgleichungen,hough1983integral}, we use the linear simplification of the formation proposed by Symm \cite[]{symm1967numerical} in Eq. (\ref{eqa:1}) instead.

With the separation in Eq. (\ref{eqa:2}), the real parts of Eq. (\ref{eqa:1}) would form the {\emph{collocation condition}}, which determines the {\emph{charges}} $ Q_{k} $ in Eq. (\ref{eqa:1}) and can be expressed as
\begin{equation}
  \label{eqa:3}
  \sum\limits_{k=1}^{2N} Q_{k}\ln|z_{i}-Z_{k}| + \varGamma = -\ln|z_{i}-z_{c}|, \quad i = 1,2,3,\cdots,2N
\end{equation}
where $ z_{i} $ are called {\emph{collocation points}}, which are artificially selected and sequentially distributed in the induced orientation (to always keep the exterior region left) along boundary $ {\bm{C}}_{2} $. To guarantee the single-valuedness of the argument, the {\emph{charges}} should meet the following requirement as
\begin{equation}
  \label{eqa:4}
  \sum\limits_{k=1}^{2N} Q_{k} = 0
\end{equation}
The {\emph{charge points}} are given in Amano's manner \cite[]{amano1994charge} as
\begin{equation}
  \label{eqa:5}
  \left\{
    \begin{aligned}
      & Z_{k} = (z_{k}-z_{c}) + k_{2} \cdot h_{k} \cdot {\rm{e}}^{{\rm{i}}\varTheta_{k}} \\
      & h_{k} = \frac{1}{2}( |z_{k+1}-z_{k}| + |z_{k}-z_{k-1}| ) \\
      & \varTheta_{k} = \arg(z_{k+1}-z_{k-1}) - \frac{\pi}{2} \\
    \end{aligned}
  \right.
  , \quad
  \begin{aligned}
    & k = 1,2,3,\cdots,2N
  \end{aligned}
\end{equation}
where $ k_{2} $ are positive and real constants called {\emph{assignment factor}}. Eqs. (\ref{eqa:3}) and (\ref{eqa:4}) complete the simultaneous linear system to seek the forward conformal mappings of region $ \overline{\bm{\varOmega}} $ onto its image $ \overline{\bm{\varOmega}}^{\prime} $ in Eq. (\ref{eq:3.4}).

In pratical computation, a branch cut discontinuity of argument $ \arg(z-Z_{k}) $ in the imaginary part of Eq. (\ref{eqa:2}) exists along the negative real axis of coordinate system $ z = x+{\rm{i}}y $ in many coding languages, such as {\tt{Fortran}}, {\tt{Python}}, {\tt{C/C++}}, {\tt{Matlab}}, and {\tt{Mathematica}}. To cancel such a discontinuity of argument in coding, Eq. (\ref{eqa:1}) is modified by using Eq. (\ref{eqa:4}) as
\begin{equation}
  \label{eqa:1'}
  \tag{A.1'}
  w(z) = (z-z_{c}) \cdot \exp \left( \varGamma + \sum\limits_{k=1}^{2N} Q_{k} \ln\frac{z-Z_{k}}{z-z_{0}} \right), \quad z \in \overline{\bm{\varOmega}}
\end{equation}
where $ z_{0} $ is an arbitrary point within boundary $ {\bm{C}}_{2} $, and $ z_{0} = z_{c} $ can generally be used.

The backward conformal mapping can be obtained using the following series as
\begin{equation}
  \label{eqa:6}
  z(w) = \sum\limits_{k=-1}^{2N-2} q_{k}w^{-k}, \quad w \in \overline{\bm{\varOmega}}^{\prime}
\end{equation}
where $ q_{k} $ denote the complex coefficients to be determined, which can be determined using point correspondence as 
\begin{equation}
  \label{eqa:7}
  \sum\limits_{k=-1}^{2N-2} q_{k}w_{i}^{-k} = z_{i}, \quad i = 1,2,3,\cdots,2N
\end{equation}
where mapping points of corresponding {\emph{collocation points}} can be given by Eq. (\ref{eqa:1'}) as
\begin{equation}
  \label{eqa:8}
  w_{i} = (z_{i}-z_{c}) \cdot \exp \left( \varGamma + \sum\limits_{k=1}^{2N} Q_{k} \ln\frac{z_{i}-Z_{k}}{z_{i}-z_{0}} \right), \quad
  \begin{aligned}
    & i = 1,2,3,\cdots,2N
  \end{aligned}
\end{equation}
Eqs. (\ref{eqa:6})-(\ref{eqa:8}) complete the backward conformal mapping of mapping region $ \overline{\bm{\varOmega}}^{\prime} $ onto its preimage $ \overline{\bm{\varOmega}} $.

\subsection{Second step of composite conformal mapping}
\label{sec:second-step-comp}

The second-step bidirectional conformal mapping is based on the Verruijt's conformal mapping \cite[]{Verruijt1997displacement,Verruijt1997traction}.

\begin{subequations}
  \label{eqa:9}
  \begin{equation}
    \label{eqa:9a}
    \zeta(w) = \frac{w+{\rm{i}}h+{\rm{i}}a}{w+{\rm{i}}h-{\rm{i}}a}, \quad w \in \overline{\bm{\varOmega}}^{\prime}
  \end{equation}
  \begin{equation}
    \label{eqa:9b}
    w(\zeta) = -{\rm{i}}a\frac{1+\zeta}{1-\zeta}+{\rm{i}}h, \quad \zeta \in \overline{\bm{\omega}}
  \end{equation}
  with
  \begin{equation}
    \label{eqa:9c}
    \left\{
      \begin{aligned}
        & a = h\frac{1-\alpha^{2}}{1+\alpha^{2}} \\
        & \alpha = \frac{r}{h+\sqrt{h^{2}-r^{2}}} \\
      \end{aligned}
    \right.
  \end{equation}
\end{subequations}
where $ r = 1 $ is the radius of the unit circular cavity in the mapping plane $ w = u + {\rm{i}} v $. The value of $ h $ is no longer the depth of the upper boundary of the mapping geomaterial $ \overline{\bm{\varOmega}}^{\prime} $, since the upper boundary is slightly curved. To establish the bidirectional conformal mapping in Eq. (\ref{eqa:9}), we choose $ h = \Im T_{2}^{\prime} $, as long as the joint point $ T_{2}^{\prime} $ is far away from the origin of coordinate system $ w = \Re w + {\rm{i}}\Im w $. Therefore, Eq. (\ref{eqa:9}) bidirectionally maps the quasi lower half plane containing a unit circular cavity $ \overline{\bm{\varOmega}}^{\prime} $ in mapping plane $ w = \Re w + {\rm{i}}\Im w $ in Fig. \ref{fig:3}a-2 and the quasi unit annulus $ \overline{\bm{\omega}} $ in mapping plane $ \zeta = \rho \cdot {\rm{e}}^{{\rm{i}}\theta} $ in Fig. \ref{fig:3}a-3 onto each other.

\clearpage

\clearpage

\clearpage
\begin{figure}[htb]
  \centering
  \begin{tikzpicture}
    \tikzstyle{every node} = [scale = 0.7]
    \fill [gray!30] (0,0) rectangle (8,-6);
    \draw [cyan, line width = 1pt] (0.5,0) -- (7.5,0);
    \fill [pattern = north east lines] (0,0) rectangle (0.5,0.2);
    \fill [pattern = north east lines] (8,0) rectangle (7.5,0.2);
    \draw [line width = 1pt] (0,0) -- (0.5,0);
    \draw [line width = 1pt] (7.5,0) -- (8,0);
    \draw [->] (4,0) -- (4.5,0) node [above] {$x$};
    \draw [->] (4,0) -- (4,0.5) node [right] {$y$};
    \node at (4,0) [above right] {$O$};
    \fill [white] (4,-4) arc [start angle = -90, end angle = 90, x radius = 1.5, y radius = 1] arc [start angle = 90, end angle = 270, x radius = 1.2, y radius = 1];
    \draw [violet, line width = 1pt] (4,-4) arc [start angle = -90, end angle = 90, x radius = 1.5, y radius = 1] arc [start angle = 90, end angle = 270, x radius = 1.2, y radius = 1];
    \draw [blue, dashed, line width = 1pt] (4,-3) circle [radius = 1];
    \draw [<-, blue] (5,-3) -- (6,-3) node [right, align = left] {Expected symmetrical \\ cavity boundary};
    \node at (4,-4) [below right, violet] {Real asymmetrical cavity boundary};
    \node at (4,-4) [below left, violet] {$ {\bm{C}}_{2} $};
    \draw [dash dot, blue] (4,0) -- (4,-6);
    \node at (4,-5) [right, blue] {Expected symmetrical line};
    \fill [violet] (4,-3) circle [radius = 0.05];
    \node at (4,-3) [right, violet] {$ z_{c} $};
    \fill [red] (0.5,0) circle [radius = 0.05];
    \fill [red] (7.5,0) circle [radius = 0.05];
    \node at (0.5,0) [above right, red] {$ T_{1} $};
    \node at (7.5,0) [above left, red] {$ T_{2} $};
    \node at (0,-6) [above right] {$ \overline{\bm{\varOmega}} $};
    \node at (8,-6) [above left] {$ z = x + {\rm{i}}y $};
    \node at (4,0) [below right, cyan] {Straight ground surface};
    \node at (4,0) [below left, cyan] {$ {\bm{C}}_{1f} $};
    \node at (0,0) [below right] {$ {\bm{C}}_{1c} $};
    \node at (8,0) [below left] {$ {\bm{C}}_{1c} $};
    \node at (4,-6) [below] {(a-1) Geomaterial in lower half plane containing an asymmetrical cavity in physical plane $ z $};

    \draw [->, line width = 1.5pt] (3.5,-6.6) -- (3.5,-7.6);
    \node at (3.5,-7.1) [left] {$ w = w(z) $};
    \draw [<-, line width = 1.5pt] (4.5,-6.6) -- (4.5,-7.6);
    \node at (4.5,-7.1) [right] {$ z = z(w) $};

    \fill [gray!30] (0,-8) rectangle (8,-14);
    \draw [cyan, line width = 1pt] (1,-8) -- (6.5,-8);
    \fill [pattern = north east lines] (0,-8) rectangle (1,-7.8);
    \fill [pattern = north east lines] (8,-8) rectangle (6.5,-7.8);
    \draw [line width = 1pt] (0,-8) -- (1,-8);
    \draw [line width = 1pt] (8,-8) -- (6.5,-8);
    \fill [white] (4,-10) circle [radius = 0.8];
    \draw [violet, line width = 1pt] (4,-10) circle [radius = 0.8];
    \draw [->] (4,-10) -- (4.5,-10) node [above] {$ \Re w $};
    \draw [->] (4,-10) -- (4,-9.5) node [right] {$ \Im w $};
    \node at (4,-10) [above right] {$ o $};
    \fill [red] (1,-8) circle [radius = 0.05];
    \fill [red] (6.5,-8) circle [radius = 0.05];
    \node at (1,-8) [below, red] {$ T_{1}^{\prime} $};
    \node at (6.5,-8) [below, red] {$ T_{2}^{\prime} $};
    \node at (0,-14) [above right] {$ \overline{\bm{\varOmega}}^{\prime} $};
    \node at (8,-14) [above left] {$ w = \Re w + {\rm{i}}\Im w $};
    \draw [dashed, red] (4,-10) -- (7.5,-10);
    \draw [<->, red] (7,-10) -- (7,-8);
    \node at (7,-9) [above, rotate = 90, red] {$ h = \Im T_{2} $};
    \draw [->, violet] (4,-10) -- ({4+0.8*cos(135)},{-10+0.8*sin(135)}) node [above left, draw, shape = rectangle] {$ r = 1 $};
    \fill [violet] (4,-10) circle [radius = 0.05];
    \node at (4,-8) [below, cyan] {Slightly curved ground surface $ {\bm{C}}_{1f}^{\prime} $};
    \node at (0,-8) [below right] {$ {\bm{C}}_{1c}^{\prime} $};
    \node at (8,-8) [below left] {$ {\bm{C}}_{1c}^{\prime} $};
    \node at (4,-10.8) [below, violet, align = center] {$ {\bm{C}}_{2}^{\prime} $ \\ Unit circle};
    \node at (4,-14) [below, align = left] {(a-2) Geomaterial in quasi lower half plane containing a unit cavity in mapping plane $ w $};

    \draw [->, line width = 1.5pt] (3.5,-14.6) -- (3.5,-15.8);
    \node at (3.5,-15.2) [left] {$ \zeta = \zeta(w) $};
    \draw [<-, line width = 1.5pt] (4.5,-14.6) -- (4.5,-15.8);
    \node at (4.5,-15.2) [right] {$ w = w(\zeta) $};

    \fill [gray!30] (4,-19) circle [radius = 3];
    \fill [pattern = north east lines] ({4+3.2*cos(10)},{-19+3.2*sin(10)}) arc [start angle = 10, end angle = -20, radius = 3.2] -- ({4+3*cos(-20)},{-19+3*sin(-20)}) arc [start angle = -20, end angle = 10, radius = 3] -- ({4+3.2*cos(10)},{-19+3.2*sin(10)});
    \draw [cyan, line width = 1pt] (4,-19) circle [radius = 3];
    \draw [line width = 1pt] ({4+3*cos(-20)},{-19+3*sin(-20)}) arc [start angle = -20, end angle = 10, radius = 3];
    \fill [white] (4,-19) circle [radius = 1];
    \draw [violet, line width = 1pt] (4,-19) circle [radius = 1];
    \draw [->] (4,-19) -- (4.5,-19) node [above] {$ \rho $};
    \draw [->] (4.3,-19) arc [start angle = 0, end angle = 360, radius = 0.3];
    \node at (4.2,-19) [below right] {$ \theta $};
    \node at (4,-19) [above] {$o$};
    \fill (4,-19) circle [radius = 0.05];
    \fill [red] ({4+3*cos(10)},{-19+3*sin(10)}) circle [radius = 0.05];
    \fill [red] ({4+3*cos(-20)},{-19+3*sin(-20)}) circle [radius = 0.05];
    \node at ({4+3*cos(-20)},{-19+3*sin(-20)}) [left, red] {$ t_{1} $};
    \node at ({4+3*cos(10)},{-19+3*sin(10)}) [left, red] {$ t_{2} $};
    \node at (5,-22) [right] {$ \zeta = \rho \cdot {\rm{e}}^{{\rm{i}}\theta} $};
    \node at (4,-22) [above] {$ \overline{\bm{\omega}} $};
    \node at (1,-19) [left, cyan, align = right] {$ {\bm{c}}_{1f} $ \\ Quasi \\ circular \\ outer \\ periphery};
    \draw [->, violet] (4,-19) -- ({4+1*cos(135)},{-19+1*sin(135)}) node [above left, draw, shape = rectangle] {$ r_{i} = \alpha $};
    \draw [->, cyan] (4,-19) -- ({4+3*cos(-135)},{-19+3*sin(-135)}) node [below left, draw, shape = rectangle, cyan] {$ r_{o} \approx 1 $};
    \node at (4,-20) [below, violet] {$ {\bm{c}}_{2} $};
    \node at (7,-19) [left] {$ {\bm{c}}_{1c} $};
    \node at (4,-22.5) [below] {(a-3) Geomaterial in quasi unit annulus in mapping plane $ \zeta $};

    \draw [cyan, ->, line width = 1.5pt] (8.5,-18.5) -- (9.5,-18.5);
    \node at (9,-18.3) [above, cyan, align = center, draw, shape = rectangle] {Geomaterial \\ reselection};
    \draw [brown, <-, line width = 1.5pt] (8.5,-19.5) -- (9.5,-19.5);
    \node at (9,-19.7) [below, brown, align = center, draw, shape = rectangle] {Geomaterial \\ simulation};

    \fill [gray!30] (14,-19) circle [radius = 3];
    \fill [pattern = north east lines] ({14+3.2*cos(10)},{-19+3.2*sin(10)}) arc [start angle = 10, end angle = -20, radius = 3.2] -- ({14+3*cos(-20)},{-19+3*sin(-20)}) arc [start angle = -20, end angle = 10, radius = 3] -- ({14+3.2*cos(10)},{-19+3.2*sin(10)});
    \draw [brown, line width = 1pt] (14,-19) circle [radius = 3];
    \draw [line width = 1pt] ({14+3*cos(-20)},{-19+3*sin(-20)}) arc [start angle = -20, end angle = 10, radius = 3];
    \fill [white] (14,-19) circle [radius = 1];
    \draw [violet, line width = 1pt] (14,-19) circle [radius = 1];
    \draw [->] (14,-19) -- (14.5,-19) node [above] {$ \rho $};
    \draw [->] (14.3,-19) arc [start angle = 0, end angle = 360, radius = 0.3];
    \node at (14.2,-19) [below right] {$ \theta $};
    \node at (14,-19) [above] {$o$};
    \fill (14,-19) circle [radius = 0.05];
    \fill [red] ({14+3*cos(10)},{-19+3*sin(10)}) circle [radius = 0.05];
    \fill [red] ({14+3*cos(-20)},{-19+3*sin(-20)}) circle [radius = 0.05];
    \node at ({14+3*cos(-20)},{-19+3*sin(-20)}) [left, red] {$ t_{10} $};
    \node at ({14+3*cos(10)},{-19+3*sin(10)}) [left, red] {$ t_{20} $};
    \node at (15,-22) [right] {$ \zeta = \rho \cdot {\rm{e}}^{{\rm{i}}\theta} $};
    \node at (14,-22) [above] {$ \overline{\bm{\omega}}_{0} $};
    \node at (14,-22) [below] {$ \bm{\omega}_{0}^{-} $};
    \node at (11,-19) [left, brown, align = right] {$ {\bm{c}}_{1f0} $};
    \draw [->, violet] (14,-19) -- ({14+1*cos(135)},{-19+1*sin(135)}) node [above left, draw, rectangle] {$ r_{i} = \alpha $};
    \draw [->, brown] (14,-19) -- ({14+3*cos(-135)},{-19+3*sin(-135)}) node [below left, draw, shape = rectangle] {$ r_{o} = 1 $};
    \node at (17,-19) [left] {$ {\bm{c}}_{1c0} $};
    \node at (14,-20) [below, violet] {$ {\bm{c}}_{20} $};
    \node at (14,-22.5) [below, align = left] {(b-3) Reselected geomaterial in mapping plane $ \zeta $ \\ within corresponding unit annulus};

    \draw [<-, line width = 1.5pt] (14,-14.8) -- (14,-15.8);
    \node at (14,-15.3) [left] {$ w = w(\zeta) $};

    \fill [gray!30] (10,-8) rectangle (18,-14);
    \draw [brown, line width = 1pt] (11,-8) -- (16.5,-8);
    \fill [pattern = north east lines] (10,-8) rectangle (11,-7.8);
    \fill [pattern = north east lines] (18,-8) rectangle (16.5,-7.8);
    \draw [line width = 1pt] (10,-8) -- (11,-8);
    \draw [line width = 1pt] (18,-8) -- (16.5,-8);
    \fill [white] (14,-10) circle [radius = 0.8];
    \draw [violet, line width = 1pt] (14,-10) circle [radius = 0.8];
    \draw [->] (14,-10) -- (14.5,-10) node [above] {$ \Re w $};
    \draw [->] (14,-10) -- (14,-9.5) node [right] {$ \Im w $};
    \node at (14,-10) [above right] {$ o $};
    \fill [red] (11,-8) circle [radius = 0.05];
    \fill [red] (16.5,-8) circle [radius = 0.05];
    \node at (11,-8) [below, red] {$ T_{10}^{\prime} $};
    \node at (16.5,-8) [below, red] {$ T_{20}^{\prime} $};
    \node at (10,-14) [above right] {$ \overline{\bm{\varOmega}}_{0}^{\prime} $};
    \node at (18,-14) [above left] {$ w = \Re w + {\rm{i}}\Im w $};
    \draw [dashed, red] (14,-10) -- (17.5,-10);
    \draw [<->, red] (17,-10) -- (17,-8);
    \node at (17,-9) [above, rotate = 90, red] {$ h = \Im T_{2} $};
    \draw [->, violet] (14,-10) -- ({14+0.8*cos(135)},{-10+0.8*sin(135)}) node [above left, draw, shape = rectangle] {$ r = 1 $};
    \fill [violet] (14,-10) circle [radius = 0.05];
    \node at (14,-8) [below, brown] {Straight ground surface $ {\bm{C}}_{1f0}^{\prime} $};
    \node at (10,-8) [below right] {$ {\bm{C}}_{1c0}^{\prime} $};
    \node at (18,-8) [below left] {$ {\bm{C}}_{1c0}^{\prime} $};
    \node at (14,-10.8) [below, violet] {$ {\bm{C}}_{20}^{\prime} $};
    \node at (14,-14) [below, align = left] {(b-2) Reselected geomaterial in mapping plane $ w $ \\ after first-step backward conformal mapping $ w = w(\zeta) $};

    \draw [<-, line width = 1.5pt, brown] (8.5,-11) -- (9.5,-11);
    \node at (9,-11.2) [below, brown, align = center, draw, shape = rectangle] {Geomaterial \\ simulation};

    \draw [<-, line width = 1.5pt] (14,-6.8) -- (14,-7.6);
    \node at (14,-7.2) [left] {$ z = z(w) $};

    \fill [gray!30] (10,0) rectangle (18,-6);
    \draw [brown, line width = 1pt] (10.5,0) -- (17.5,0);
    \fill [pattern = north east lines] (10,0) rectangle (10.5,0.2);
    \fill [pattern = north east lines] (18,0) rectangle (17.5,0.2);
    \draw [line width = 1pt] (10,0) -- (10.5,0);
    \draw [line width = 1pt] (17.5,0) -- (18,0);
    \draw [->] (14,0) -- (14.5,0) node [above] {$x$};
    \draw [->] (14,0) -- (14,0.5) node [right] {$y$};
    \node at (14,0) [above right] {$O$};
    \fill [white] (14,-4) arc [start angle = -90, end angle = 90, x radius = 1.5, y radius = 1] arc [start angle = 90, end angle = 270, x radius = 1.2, y radius = 1];
    \draw [violet, line width = 1pt] (14,-4) arc [start angle = -90, end angle = 90, x radius = 1.5, y radius = 1] arc [start angle = 90, end angle = 270, x radius = 1.2, y radius = 1]; 
    \fill [red] (10.5,0) circle [radius = 0.05];
    \fill [red] (17.5,0) circle [radius = 0.05];
    \node at (10.5,0) [above right, red] {$ T_{10} $};
    \node at (17.5,0) [above left, red] {$ T_{20} $};
    \node at (10,-6) [above right] {$ \overline{\bm{\varOmega}}_{0} $};
    \node at (18,-6) [above left] {$ z = x + {\rm{i}}y $};
    \node at (14,0) [below, brown] {Slightly curved ground surface $ {\bm{C}}_{1f0} $};
    \node at (10,0) [below right] {$ {\bm{C}}_{1c0} $};
    \node at (18,0) [below left] {$ {\bm{C}}_{1c0} $};
    \node at (14,-4) [below, violet, align = center] { $ {\bm{C}}_{20} $ \\ Original asymmetrical cavity boundary};
    \node at (14,-3) [right, violet] {$ z_{c0} $};
    \draw [->, red] (14,-3) -- (14,-2.5) node [right] {$ R_{y} $};
    \fill [violet] (14,-3) circle [radius = 0.05];
    \node at (14,-6) [below, align = left] {(b-1) Simulating geomaterial in mapping plane $ z $ \\ after second-step backward conformal mapping $ z = z(w) $};

    \draw [<-, line width = 1.5pt, brown] (8.5,-3) -- (9.5,-3);
    \node at (9,-3.2) [below, brown, align = center, draw, shape = rectangle] {Geomaterial \\ simulation};
  \end{tikzpicture}
  \caption{Schematic diagram of bidirectional composite conformal mapping and geomaterial simulation}
  \label{fig:3}
\end{figure}

\clearpage

\clearpage

\clearpage

\clearpage

\clearpage

\clearpage

\clearpage

\clearpage

\clearpage

\clearpage

\clearpage

\clearpage

\clearpage

\clearpage
\bibliographystyle{plainnat}
\bibliography{/home/lin/my_bibtex}  

\end{document}